\documentclass[leqno,12pt]{book}

\usepackage{amssymb}
\usepackage{makeidx}
\usepackage{graphicx}
\usepackage{fancyhdr}
\usepackage{setspace}
\usepackage{xypic}
\input xy
\xyoption{all}

\newtheorem{theorem}{Theorem}[section]
\newtheorem{definition}{Definition}[section]
\newtheorem{proposition}{Proposition}[section]
\newtheorem{corollary}{Corollary}[section]
\newtheorem{lemma}{Lemma}[section]
\newtheorem{remark}{Remark}[section]

\renewcommand \theequation{\thesection.\arabic{equation}}
\renewcommand{\bibname}{References}

\makeindex

\def\v{\varphi}
\def\o{\omega}

\def\a{\alpha}
\def\b{\beta}
\def\l{\lambda}
\def\d{\delta}
\def\p{\partial}
\def\t{\theta}
\def\vp{\varepsilon}
\def\s{\sigma}
\def\D{\Delta}
\def\k{\kappa}
\def\G{\Gamma}
\def\O{\Omega}
\def\L{\Lambda}
\def\dfrac{\displaystyle\frac}
\def\({\left(}
\def\){\right)}

\def\spa{\vspace*{-.3cm}\\ }
\def\cal#1{\mathcal{#1}}

\makeatletter
\renewcommand*\l@part[2]{%
  \ifnum \c@tocdepth >\m@ne
    \addpenalty{-\@highpenalty}%
    \vskip 1.0em \@plus\p@
    \setlength\@tempdima{1.5em}%
    \begingroup
      \parindent \z@ \rightskip \@pnumwidth
      \parfillskip -\@pnumwidth
      \leavevmode \bfseries
      \advance\leftskip\@tempdima
      \hskip -\leftskip
      #1\nobreak\
       \leaders\hbox{$\m@th
        \mkern \@dotsep mu\hbox{.}\mkern \@dotsep
        mu$}\hfil\nobreak\hb@xt@\@pnumwidth{\hss #2}\par
      \penalty\@highpenalty
    \endgroup
  \fi}
\makeatother
\makeatletter
\renewcommand*\l@chapter[2]{%
  \ifnum \c@tocdepth >\m@ne
    \addpenalty{-\@highpenalty}%
    \vskip 1.0em \@plus\p@
    \setlength\@tempdima{1.5em}%
    \begingroup
      \parindent \z@ \rightskip \@pnumwidth
      \parfillskip -\@pnumwidth
      \leavevmode \bfseries
      \advance\leftskip\@tempdima
      \hskip -\leftskip
      #1\nobreak\
       \leaders\hbox{$\m@th
        \mkern \@dotsep mu\hbox{.}\mkern \@dotsep
        mu$}\hfil\nobreak\hb@xt@\@pnumwidth{\hss #2}\par
      \penalty\@highpenalty
    \endgroup
  \fi}
\makeatother
\makeatletter
\renewcommand\ps@plain{\ps@empty}
\makeatother
\makeatletter
\def\l@part#1#2{\addpenalty{\@secpenalty}%
   \addvspace{2em \@plus\p@}%
   \begingroup
     \parindent \z@
     \rightskip \z@ \@plus 5em
     \hrule\vskip5\p@
     \bfseries\boldmath
     \leavevmode
     #1\par
     \vskip5\p@
     \hrule
     \vskip\p@
     \nobreak
   \endgroup}\makeatother
\makeatletter
\def\@part[#1]#2{%
     \ifnum \c@secnumdepth >-2\relax
       \refstepcounter{part}%
       \addcontentsline{toc}{part}{\partname~\thepart\hspace{1em}#1}%
     \else
       \addcontentsline{toc}{part}{#1}%
     \fi
     \markboth{}{}%
     {\centering
      \interlinepenalty \@M
      \normalfont
      \ifnum \c@secnumdepth >-2\relax
        \huge\bfseries \partname~\thepart
        \par
        \vskip 20\p@
      \fi
      \Huge \bfseries #2\par}%
     \@endpart}
\makeatother

\begin{document}
\begin{titlepage}
\setstretch{1.4}
\centerline{\large{\textbf{RADU MIRON}}}
\vspace*{4cm}

\begin{center}{\Large{\textbf{THE GEOMETRY OF MYLLER CONFIGURATIONS. \vspace*{2mm}\\ Applications to Theory of Surfaces and Nonholonomic Manifolds}}}
\end{center}

\vspace*{8cm}\centerline{{\bf Editura Academiei Rom\^{a}ne}}\vspace*{-2mm}
\centerline{\footnotesize{\bf Bucure\c{s}ti, 2010}}
\end{titlepage}
\setstretch{1.2}

\pagestyle{fancy}

\fancyhf{} \fancyhead[LE,RO]{\slshape\thepage}

\newcommand{\helv}{\fontsize{9}{11}\selectfont}
\fancyhf{} \fancyhead[LE,RO]{\helv\thepage} \fancyhead[LO]{\helv
\nouppercase\rightmark} \fancyhead[RE]{\helv\nouppercase\leftmark}

% rid of headers
   %\renewcommand{\headrulewidth}{0.4pt} % and the line

%\pagestyle{myheadings}
%\thispagestyle{empty}
\markboth{}{}

\vspace*{5cm}

\begin{flushright}
{Dedicated to academicians\\ {\bf Alexandru Myller} and {\bf Octav Mayer},\\ two great Romanian mathematicians}
\end{flushright}

\cleardoublepage

\newpage

\vspace*{5cm}

\hspace{3cm}{\bf MOTTO}

\medskip

\hspace{3cm}{\it{It is possible to bring for the great disappeared\\ \hspace*{3.4cm} scientists
varied homages. Sometimes, the strong\\ \hspace*{3.4cm} wind of progress erases the
trace of their steps. To\\ \hspace*{3.4cm} not forget them means to continue their
works,\\ \hspace*{3.4cm} connecting them to the living present.}}

\bigskip

\hspace*{2.85cm} Octav Mayer

\addcontentsline{toc}{chapter}{Preface}

\chapter*{Preface}

In 2010, the Mathematical Seminar of the ``Alexandru Ioan Cuza''
University of Ia\c{s}i comes to its 100th anniversary of
prodigious existence.

The establishing of the Mathematical Seminar by Alexandru Myller
also marked the beginning of the School of Geometry in Ia\c{s}i,
developed in time by prestigious mathematicians of international
fame, Octav Mayer, Gheorghe Vr\^{a}nceanu, Grigore Moisil, Mendel
Haimovici, Ilie Popa, Dimitrie Mangeron,  Gheorghe Gheorghiev and
many others.

Among the first paper works, published by Al. Myller and O. Mayer,
those concerning the generalization of the Levi-Civita parallelism
must be specified, because of the relevance for the international
recognition of the academic School of Ia\c{s}i. Later on, through
the effort of two generations of Ia\c{s}i based mathematicians,
these led to a body of theories in the area of differential
geometry of Euclidian, affine and projective spaces.

At the half-centenary of the Mathematical Seminary, in 1960, the
author of the present opuscule synthesized the field$'$s results
and laid it out in the form of a ``whole, superb theory'', as
mentioned by Al. Myller. In the same time period, the book {\it The
geometry of the Myller configurations} was published by the
Technical Publishing House. It represents, as mentioned by Octav
Mayer in the Foreword, ``the most precious tribute ever offered to
Alexandru Myller''. Nowadays, at the 100th anniversary of the
Mathematical Seminary ``Alexandru Myller'' and 150 years after the
enactment, made by Alexandru Ioan Cuza, to set up the University
that carries his name, we are going to pay homage to those two
historical acts in the Romanian culture, through the publishing,
in  English, of the volume {\it The Geometry of the Myller
Configurations}, completed with an ample chapter, containing new
results of the Romanian geometers, regarding applications in the
studying of non-holonomic manifolds in the Euclidean space. On
this occasion, one can better notice the undeniable value of the
achievements in the field made by some great Romanian
mathematicians, such as Al. Myller, O. Mayer, Gh. Vr\u{a}nceanu,
Gr. Moisil, M. Haimovici, I. Popa, I. Creang\u{a} and Gh. Gheorghiev.

The initiative for the re-publishing of the book belongs to the
leadership of the ``Alexandru Ioan Cuza'' University of Ia\c{s}i, to the local branch of the Romanian Academy, as
well as NGO Formare Studia Ia\c{s}i. A competent
assistance was offered to me, in order to complete this work, by
my former students and present collaborators, Professors Mihai
Anastasiei and Ioan Buc\u{a}taru, to whom I express my utmost gratitude.

\vspace*{1cm}

\,\,\, Ia\c{s}i, 2010 \hfill Acad. Radu Miron

\cleardoublepage 
\addcontentsline{toc}{chapter}{Preface of the book  {\it Geometria configura\c{t}iilor Myller}
written in 1966 by Octav Mayer, former member of the Romanian Academy}
\addcontentsline{toc}{chapter}{A short biography of Al. Myller}

\chapter*{Preface of the book {\it Geometria configura\c{t}iilor Myller}
written in 1966 by Octav Mayer, former member of the Romanian
Academy {\Large{(Translated from Romanian)}}}

One can bring to the great scientists, who passed away, various homages. Some times, the strong wind
of the progress wipes out the trace of their steps. To not forget them
means to continue their work, connecting them to the living
present. In this sense, this scientific work is the most precious homage
which can be dedicated to Alexandru Myller.

Initiator of a modern education of Mathematics at the University
of Iassy, founder of the Geometry School, which is still
flourishing today, in the third generation, Alexandru Myller was
also a hardworking researcher, well known inside the country and
also abroad due to his papers concerning Integral Equations and
Differential Geometry.

Our Academy elected him as a member and published his scientific
work. The ``A. Humboldt'' University from Berlin awarded him the
title of ``doctor Honoris Causa'' for ``special efforts in creating
an independent Romanian Mathematical School''.

Some of Alexandru Myller$'$s discoveries have penetrated fruitfully the
impetuous torrent of ideas, which have changed the Science of Geometry
during the first decades of the century. Among others, it is the case
of ``Myller configurations''. The reader would be perhaps interested
to find out some more details about how he discovered these configurations.

Let us imagine a surface $S$, on which there is drawn a curve $C$
and, along it, let us circumscribe to $S$ a developing surface
$\Sigma $; then we apply (develop) the surface $\Sigma$ on a
plane $\pi $. The points of the curve $C$, connected to the
respective tangent planes, which after are developed overlap each
other on the plane $\pi $, are going to represent a curve $C'$
into this plane.
In the same way, a series $(d)$ of  tangent directions to the
surface $S$ at the points of the curve $C$ becomes developing, on
the plane $\pi $, a series $(d')$ of directions getting out from
the points of the curve $C'$. The directions $(d)$ are parallel on
the surface $S$ if the directions $(d')$ are parallel in the
common sense.

Going from this definition of T. Levi-Civita parallelism (valid in
the Euclidian space), Alexandru Myller arrived to a more general
concept in a sensible process of abstraction. Of course, it was
not possible to leave the curve $C$ aside, neither the
directions $(d)$ whose parallelism was going to be defined
in a more general sense. What was left aside was the surface $S$.

For the surface $\Sigma$ was considered, in a natural way, the enveloping
of the family of planes constructed from the points of the curve
$C$, planes in which are given the directions $(d)$. Keeping
unchanged the remainder of the definition one gets to what Alexandru
Myller called ``parallelism into a family of planes''.

A curve $C$, together with a family of planes on its points and a
family of given directions (in an arbitrary way) in these planes constitutes
what the author called a ``Myller configuration''.

It was considered that this new introduced notion had a central
place into the classical theory of surfaces, giving the
possibility to interpret and link among them many particular
facts. It was obvious that this notion can be successfully applied
in the Geometry of other spaces different from the Euclidian one.
Therefore the foreworded study worthed all the efforts made by the
author, a valuable mathematician from the third generation of the
Ia\c si Geometry School.

By recommending this work, we believe that the reader (who needs
only basic Differential Geometry) will be attracted by the clear
lecture of Radu Miron and also by the beauty of the subject, which
can still be developed further on.

\bigskip

\hfill Octav Mayer, 1966

\newpage

\chapter*{A short biography of\\ Al. Myller}

%\vspace*{-2cm}
\begin{figure}[h]
\begin{center}
\includegraphics[width=6cm,height=7cm]{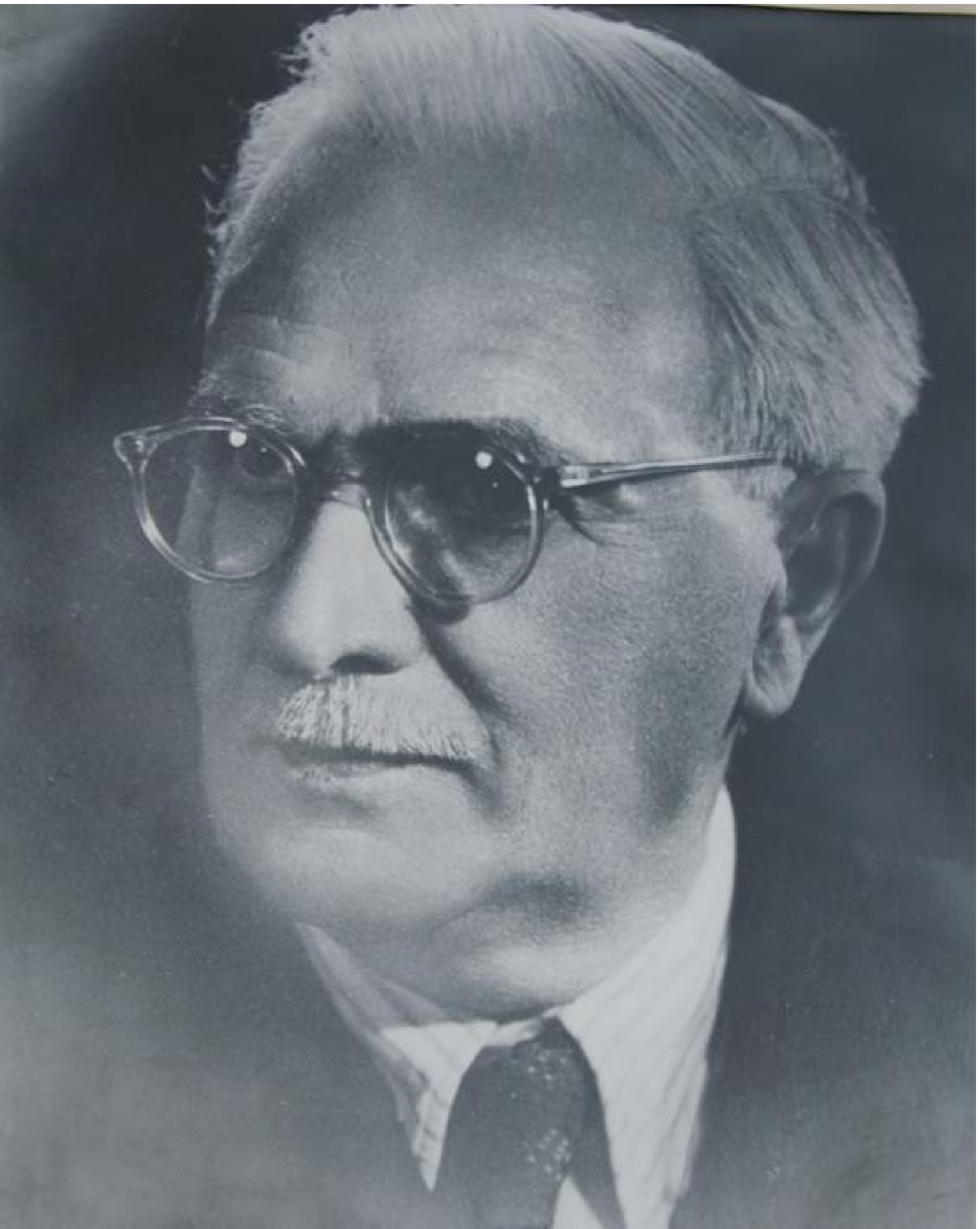}
\end{center}
\end{figure}

{\it ALEXANDRU MYLLER} was born in Bucharest in 1879, and died in
Ia\c{s}i on the 4th of July 1965. Romanian mathematician. Honorary
Member (27 May 1938) and Honorific Member (12 August 1948) of the
Romanian Academy. High School and University studies (Faculty of
Science) in Bucharest, finished with a bachelor degree in
mathematics.

He was a Professor of Mathematics at the Pedagogical Seminar; he
sustained his PhD thesis {\it Gewohnliche Differential Gleichungen
Hoherer Ordnung}, in G\"{o}ttingen, under the scientific
supervision of David Hilbert.

He worked as a Professor at the Pedagogical Seminar and the School
of Post and Telegraphy (1907 - 1908), then as a lecturer at the
University of Bucharest (1908 - 1910). He is appointed Professor
of Analytical Geometry at the University of Ia\c{s}i (1910 -
1947); from 1947 onward - consultant Professor. In 1910, he sets
up the Mathematical Seminar at the ``Alexandru Ioan Cuza``
University of Ia\c{s}i, which he endows with a library full of
didactic works and specialty journals, and which, nowadays, bears
his name. Creator of the mathematical school of Ia\c{s}i, through
which many reputable Romanian mathematicians have passed, he
conveyed to us research works in fields such as integral
equations, differential geometry and history of mathematics. He
started his scientific activity with papers on the integral
equations theory, including extensions of Hilbert$'$s results, and
then he studied integral equations and self- adjoint of even and
odd order linear differential equations, being the first
mathematician to introduce integral equations with skew symmetric
kernel.

He was the first to apply integral equations to solve problems for
partial differential equations of hyperbolic type. He was also
interested in the differential geometry, discovering a
generalization of the notion of parallelism in the Levi-Civita
sense, and introducing the notion today known as concurrence in
the Myller sense. All these led to ``the differential geometry of
Myller configurations''(R. Miron, 1960). Al. Myller, Gh.
\c{T}i\c{t}eica and O. Mayer have created ``the differential
centro-affine geometry``, which the history of mathematics refers
to as ``purely Romanian creation$!$`` He was also concerned with
problems from the geometry of curves and surfaces in Euclidian
spaces. His research outcomes can be found in the numerous
memoirs, papers and studies, published in Ia\c{s}i and abroad:
Development of an arbitrary function after Bessel$'$s functions
(1909); Parallelism in the Levi-Civita sense in a plane system
(1924); The differential centro-affine geometry of the plane
curves (1933) etc. - within the volume ``Mathematical Writings``
(1959), Academy Publishing House. Alexandru Myller has been an
Honorary Member of the Romanian Academy since 1938. Also, he was a
Doctor Honoris Causa of the Humboldt University of Berlin.

Excerpt from the volume ``Members of the Romanian Academy 1866 -
1999. Dictionary.`` - Dorina N. Rusu, page 360.

\cleardoublepage 
 %\setcounter{page}{7}
%\clearpage%\thispagestyle{empty}
%\input{pref}
%\clearpage%\thispagestyle{empty}
%\mainmatter%\thispagestyle{empty}
\addcontentsline{toc}{chapter}{Introduction}

\chapter*{Introduction}

In the differential geometry of curves in the Euclidean space $E_{3}$
one introduces, along a curve $C$, some versor fields, as tangent,
principal normal or binormal, as well as some plane fields as
osculating, normal or rectifying planes.

More generally, we can consider a versor field $(C,\overline{\xi})$
or a plane field $(C,\pi)$.
A pair $\{(C,\overline{\xi}), (C,\pi)\}$ for which
$\overline{\xi}\in \pi$, has been called in 1960 [23] by the present
author, a Myller configuration in the space $E_{3}$, denoted by
$\mathfrak{M}(C, \overline{\xi},\pi)$. When the planes $\pi$ are tangent to $C$ then we have a tangent
Myller configuration $\mathfrak{M}_{t}(C,\overline{\xi},\pi)$.

Academician Alexandru Myller studied in 1922 the notion of
parallelism of $(C,\overline{\xi})$ in the plane field $(C,\pi)$
obtaining an interesting generalization of the famous parallelism
of Levi-Civita on the curved surfaces. These investigations have
been continued by Octav Mayer which introduced new fundamental
invariants for $\mathfrak{M}(C, \overline{\xi}, \pi)$. The
importance of these studies was underlined by Levi Civita in
Addendum to his book {\it Lezioni di calcolo differentiale assoluto},
1925.

Now, we try to make a systematic presentation of the geometry of
Myller configurations $\mathfrak{M}(C,\overline{\xi},\pi)$ and
$\mathfrak{M}_{t}(C,\overline{\xi}, \pi)$ with applications to the
differential geometry of surfaces and to the geometry of
nonholonomic manifolds in the Euclidean space $E_{3}$.

Indeed, if $C$ is a curve on the surface $S\subset E_{3}$, $s$ is
the natural parameter of curve $C$ and $\overline{\xi}(s)$ is a
tangent versor field to $S$ along $C$ and $\pi(s)$ is tangent
planes field to $S$ along $C$, we have a tangent Myller
configuration $\mathfrak{M}_{t}(C,\overline{\xi},\pi)$ intrinsic
associated to the geometric objects $S,C,\overline{\xi}$.
Consequently, the geometry of the field $(C, \overline{\xi})$ on
surface $S$ is the geometry of the associated Myller
configurations $\mathfrak{M}_{t}(C,\overline{\xi},\pi)$. It is
remarkable that the geometric theory of $\mathfrak{M}_{t}$ is a
particular case of that of general Myller configuration
$\mathfrak{M}(C,\overline{\xi},\pi)$.

For a Myller configuration $\mathfrak{M}(C,\overline{\xi},\pi)$ we
determine a Darboux frame, the fundamental equations and a
complete system of invariants $G,K,T$ called, geodesic curvature,
normal curvature and geodesic torsion, respectively, of the versor
field $(C,\overline{\xi})$ in Myller configuration
$\mathfrak{M}(C,\overline{\xi}, \pi)$. A fundamental theorem, when
the functions $(G(s),K(s),$ $T(s))$ are given, can be
proven.

The invariant $G(s)$ was discovered by Al. Myller (and named by
him the deviation of parallelism). $G(s) = 0$ on curve $C$
characterizes the parallelism of versor field $(C,\overline{\xi})$
in $\mathfrak{M},$ [23], [24], [31], [32], [33]. The second invariant
$K(s)$ was introduced by O. Mayer (it was called the curvature of
parallelism). Third invariant $T(s)$ was found by E. Bortolotti
[3].

In the particular case, when $\mathfrak{M}$ is a tangent Myller
configuration $\mathfrak{M}_{t}(C,\overline{\xi},\pi)$ associated
to a tangent versor field $(C,\overline{\xi})$ on a surface $S$,
the $G(s)$ is an intrinsic invariant and $G(s)=0$ along $C$ leads
to the Levi Civita parallelism.

In configurations $\mathfrak{M}_{t}(C,\overline{\xi},\pi)$ there
exists a natural versor field $(C,\overline{\a})$, where
$\overline{\a}(s)$ are the tangent versors to curve $C$. The
versor field $(C,\overline{\a})$ has a Darboux frame $\cal{R} =
(P(s); \overline{\a}, \overline{\mu}^{*}, \overline{\nu})$ in
$\mathfrak{M}_{t}$ where $\overline{\nu}(s)$ is normal to plane
$\pi(s)$ and $\overline{\mu}^{*}(s) = \overline{\nu}(s)\times
\overline{\a}(s).$ The moving equations of $\cal{R}$ are:
\begin{eqnarray*} \dfrac{d\overline{r}}{ds} &=& \overline{\a}(s),\;\;
s\in (s_{1},s_{2})\\\dfrac{d\overline{\a}}{ds}&=&
\k_{g}(s)\overline{\mu}^{*} +\k_{n}(s)\overline{\nu},\\
\dfrac{d\overline{\mu}^{*}}{ds}&=& -\k_{g}(s)\overline{\a} +
\tau_{g}(s)\overline{\nu},\\\dfrac{d\overline{\nu}}{ds}&=&
-\k_{n}(s)\overline{\a} - \tau_{g}(s)\overline{\mu}^{*}.
\end{eqnarray*}
The functions $\k_{g}(s), \k_{n}(s)$ and $\tau_{g}(s)$ form a
complete system of invariants of the curve in $\mathfrak{M}_{t}$.

A theorem of existence and uniqueness for the versor fields
$(C,\overline{\a})$ in $\mathfrak{M}_{t}(C,$ $\overline{\a},\pi)$,
when the invariant $\k_{g}(s)$, $k_n$ and $\tau_{g}(s)$ are given, is proved. The function $\k_{g}(s)$ is called the geodesic
curvature of the curve $C$ in $\mathfrak{M}_{t}$; $\k_{n}(s)$ is
the normal curvature and $\tau_{g}(s)$ is the geodesic torsion of
$C$ in $\mathfrak{M}_{t}$.

The condition $\k_{g}(s) = 0$, $\forall s\in (s_{1},s_{2})$
characterizes the geodesic (autoparallel lines) in
$\mathfrak{M}_{t}$; $\k_{n}(s) = 0,$ $\forall s\in (s_{1},s_{2})$
give us the asymptotic lines and $\tau_{g}(s) = 0,$ $\forall s\in
(s_{1},s_{2})$ characterizes the curvature lines $C$ in
$\mathfrak{M}_{t}$

One can remark that in the case when $\mathfrak{M}_{t}$ is the
associated Myller confi\-gu\-ration to a curve $C$ on a surface $S$ we
obtain the classical theory of curves on surface $S$. It is
important to remark that the Mark Krien$'$s formula (2.10.3) leads
to the integral formula (2.10.9) of Gauss-Bonnet for surface $S$, studied by R. Miron in the book [62].

Also, if $C$ is a curve of a nonholonomic manifold $E_{3}^{2}$ in
$E_{3}$, we uniquely determine a Myller configuration
$\mathfrak{M}_{t}(C,\overline{\a},\pi)$ in which
$(C,\overline{\a})$ is the tangent versor field to $C$ and
$(C,\pi)$ is the tangent plane field to $E_{3}^{2}$ along $C$, [54], [64], [65].

In this case the geometry of $\mathfrak{M}_{t}$ is the geometry of
curves $C$ in the nonholonomic manifolds $E_{3}^{2}$. Some new
notions can be introduced as: concurrence in Myller sense of
versor fields $(C,\overline{\xi})$ on $E_{3}^{2}$, extremal
geodesic torsion, the mean torsion and total torsion of
$E_{3}^{2}$ at a point, a remarkable formula for geodesic torsion
and an indicatrix of Bonnet for geodesic torsion, which is not
reducible to a pair of equilateral hyperbolas, as in the case of
surfaces.

The nonholonomic planes, nonholonomic spheres of Gr. Moisil, can
be studied by means of techniques from the geometry of Myller
configurations $\mathfrak{M}_{t}$.

We finish the present introduction pointing out some important
developments of the geometry of Myller configurations:

The author extended the notion of Myller configuration in
Riemannian Geometry [24], [25], [26], [27]. Izu Vaisman [76] has studied
the Myller configurations in the symplectic geometry.

Mircea Craioveanu realized a nice theory of Myller configurations
in infinit dimensional Riemannian manifolds. Gheorghe Gheorghiev
developed the configuration $\mathfrak{M}$, [11], [55] in the
geometry of versor fields in the Euclidean space and applied it in
hydromechanics.

N.N. Mihalieanu studied the Myller configurations in Minkowski
spaces [61]. For Myller configurations in a Finsler, Lagrange or Hamilton spaces we refer to the paper [58] and to the books of R. Miron and M. Anastasiei [68], R. Miron, D. Hrimiuc, H. Shimada and S. Sab\u{a}u [70].

All these investigations underline the usefulness of geometry of
Myller configurations in differential geometry and its
applications.

\newpage

\cleardoublepage

\chapter{Versor fields. Plane fields in $E_{3}$}

First of all we investigate the geometry of a versor field
$(C,\overline{\xi})$ in the Euclidean space $E_{3}$, introducing an
invariant frame of Frenet type, the moving equations of this
frame, invariants and proving a fundamental theorem. The
invariants are called $K_{1}$-curvature and $K_{2}$-torsion of
$(C,\overline{\xi})$. Geometric interpretations for $K_{1}$ and
$K_{2}$ are pointed out. The parallelism of $(C,\overline{\xi})$,
concurrence of $(C,\overline{\xi})$ and the enveloping of versor
field $(C,\overline{\xi})$ are studied, too.

A similar study is made for the plane fields $(C,\pi)$ taking into
account the normal versor field $(C,\overline{\nu})$, with
$\overline{\nu}(s)$ a normal versor to the plane $\pi(s)$.

\section{Versor fields $(C,\overline{\xi})$}

\setcounter{theorem}{0}\setcounter{equation}{0}

In the Euclidian space a versor field $(C,\overline{\xi})$ can be
analytical represented in an orthonormal frame $\cal{R} =
(0;\overline{i}_{1}, \overline{i}_{2},\overline{i}_{3})$, by
\begin{equation}
\overline{r} = \overline{r}(s),\;\;\;\;\overline{\xi} = \overline{\xi}(s), \;\; s\in
(s_{1},s_{2})
\end{equation}
where $s$ is the arc length on curve $C$,
$$
\overline{r}(s) = \overline{OP}(s) = x(s)\overline{i}_{1} +
y(s)\overline{i}_{2} + z(s)\overline{i}_{3}\textrm{ and }
$$  $\overline{\xi}(s) =
\overline{\xi}(\overline{r}(s))\hspace{-0.5mm}=\hspace{-0.5mm} \overrightarrow{PQ} =
\xi^{1}(s)\overline{i}_{1}\hspace{-0.5mm} +\hspace{-0.5mm} \xi^{2}(s)\overline{i}_{2}\hspace{-0.5mm} +\hspace{-0.5mm}
\xi^{3}(s)\overline{i}_{3}$, $\|\overline{\xi}(s)\|^2\hspace{-0.5mm} =\hspace{-0.5mm}
\langle\overline{\xi}(s),\overline{\xi}(s)\rangle =1.$

\smallskip

All geometric objects considered in this book are assumed to be of
class $C^{k}$, $k\geq 3$, and sometimes of class $C^{\infty}$. The
pair $(C,\overline{\xi})$ has a geometrical meaning. It follows
that the pair $\(C,\dfrac{d\overline{\xi}}{ds}\)$ has a geometric
meaning, too. Therefore, the norm:
\begin{equation}
K_{1}(s) = \left\|\dfrac{d\overline{\xi}(s)}{ds}\right\|
\end{equation}
is an invariant of the field $(C,\overline{\xi})$.

We denote
\begin{equation}
\overline{\xi}_1(s) = \overline{\xi}(s)
\end{equation}
and let $\overline{\xi}_{2}(s)$ be the versor of vector
$\dfrac{d\overline{\xi}_{1}}{ds}$. Thus we can write
$$
\dfrac{d\overline{\xi}_{1}(s)}{ds} =K_{1}(s)\overline{\xi}_{2}(s).
$$
Evidently, $\overline{\xi}_{2}(s)$ is orthogonal to
$\overline{\xi}_{1}(s)$.

It follows that the frame
\begin{equation}
\cal{R}_{F} = (P(s); \overline{\xi}_{1}, \overline{\xi}_{2},
\overline{\xi}_{3}),\;\; \overline{\xi}_{3}(s) =
\overline{\xi}_{1}\times \overline{\xi}_{2}
\end{equation}
is orthonormal, positively oriented and has a geometrical meaning.
$\cal{R}_{F}$ is called the Frenet frame of the versor field
$(C,\overline{\xi})$.

We have:

\begin{theorem}
The moving equations of the Frenet frame $\cal{R}_{F}$ are:
\begin{equation}
\hspace*{5mm}\dfrac{d\overline{r}}{ds} =
a_{1}(s)\overline{\xi}_{1}+a_{2}(s)\overline{\xi}_{2} +
a_{3}(s)\overline{\xi}_{3}\;\;\;\; a_{1}^{2}(s)+a_{2}^{2}(s) +
a_{3}^{2}(s) =1
\end{equation}
and
\begin{eqnarray}
\dfrac{d\overline{\xi}_{1}}{ds} &=&
K_{1}(s)\overline{\xi}_{2},\nonumber\\\dfrac{d\overline{\xi}_{2}}{ds}
&=& -K_{1}(s)\overline{\xi}_{1}(s) +
K_{2}(s)\overline{\xi}_{3},\\\dfrac{d\overline{\xi}_{3}}{ds}& =&
-K_{2}(s)\overline{\xi}_{2}(s),\nonumber
\end{eqnarray}
where $K_{1}(s)>0.$ The functions $K_{1}(s), K_{2}(s), a_{1}(s),
a_{2}(s)$, $a_{3}(s),$ $s\in (s_{1},s_{2})$ are invariants of the
versor field $(C,\overline{\xi})$.\end{theorem}

The proof does not present difficulties.

The invariant $K_{1}(s)$ is called the curvature of
$(C,\overline{\xi})$ and has the same geometric interpretation as
the curvature of a curve in $E_{3}$. $K_{2}(s)$ is called the
torsion and has the same geometrical interpretation as the torsion
of a curve in $E_{3}$.

The equations (1.1.5), (1.1.6) will be called the {\it fundamental or
Frenet equations} of the versor field $(C,\overline{\xi})$.

In the case $a_{1}(s)= 1,$ $a_{2}(s) = 0,$ $a_{3}(s) =0$ the
tangent versor $\dfrac{d\overline{r}}{ds}$ is denoted by
$$
\overline{\a}(s) =
\dfrac{d\overline{r}}{ds}(s)\leqno{(1.1.5')}
$$
The equations (1.1.5), (1.1.6) are then the Frenet equations of a
curve in the Euclidian space $E_{3}$.

For the versor field $(C,\overline{\xi})$ we can formulate a
fundamental theorem:

\begin{theorem}
If the functions $K_{1}(s)>0,$ $K_{2}(s)$, $a_{1}(s), a_{2}(s),
a_{3}(s)$, $(a_{1}^{2} + a_{2}^{2} + a_{3}^{2} =1)$ of class
$C^{\infty}$ are apriori given, $s\in [a,b]$ there exists a curve
$C:[a,b]\to E_{3}$ parametrized by arclengths and a versor field
$\overline{\xi}(s), s\in [a,b]$, whose the curvature, torsion and
the functions $a_{i}(s)$ are $K_{1}(s), K_{2}(s)$ and
$a_{i}(s)$. Any two such versor fields $(C,\overline{\xi})$ differ
by a proper Euclidean motion.
\end{theorem}

For the proof one applies the same technique like in the proof of
Theorem 11, p. 45 from [75], [76].

\begin{remark}\rm
1. If $K_{1}(s) = 0,$ $s\in (s_{1}, s_{2})$ the versors
$\overline{\xi}_{1}(s)$ are parallel in $E_{3}$ along the curve
$C.$

2. The versor field $(C,\overline{\xi})$ determines a ruled
surface $S(C,\overline{\xi})$.

3. The surface $S(C,\overline{\xi})$ is a cylinder iff the
invariant $K_{1}(s)$ vanishes.

4. The surface $S(C,\overline{\xi})$ is with director plane iff
$K_{2}(s) = 0$.

5. The surface $S(C,\overline{\xi})$ is developing iff the
invariant $a_{3}(s)$ vanishes.
\end{remark}

If the surfaces $S(C,\overline{\xi})$ is a cone we say that the
versors field $(C,\overline{\xi})$ is {\it concurrent}.

\begin{theorem}
A necessary and sufficient condition for the versor field $(C,\overline{\xi})$, $(K_{1}(s)\neq 0)$,
to be concurrent is the following
$$
\dfrac{d}{ds}\left(\dfrac{a_{2}(s)}{K_{1}(s)}\right)-a_{1}(s) =
0,\;\; a_{3}(s) = 0,\;\; \forall s\in (s_{1},s_{2}).
$$
\end{theorem}

\medskip

\section{Spherical image of a versor field $(C,\overline{\xi})$}

\setcounter{theorem}{0}\setcounter{equation}{0}

Consider the sphere $\Sigma$ with center a fix point $O\in E_{3}$
and radius 1.

\begin{definition}
The spherical image of the versor field $(C,\overline{\xi})$ is the
curve $C^{*}$ on sphere $\Sigma$ given by $\overline{\xi}^{*}(s) =
\overrightarrow{OP^{*}}(s) = \overline{\xi}(s),$ $\forall s\in
(s_{1},s_{2})$.
\end{definition}

From this definition we have:
\begin{equation}
d\overline{\xi}^{*}(s) =\overline{\xi}_{2}(s) K_{1}(s)ds.
\end{equation}
Some immediate properties:

1. It follows
\begin{equation}
ds^{*} = K_{1}(s)ds.
\end{equation}

Therefore:

The arc length of the curve $C^{*}$ is
\begin{equation}
s^{*} = s_{0}^{*} +\int_{s_{0}}^{s}K_{1}(\sigma) d\sigma
\end{equation}
where $s_{0}^{*}$ is a constant and $[s_{0},s]\subset
(s_{1},s_{2})$.

2. The curvature $K_{1}(s)$ of $(C,\overline{\xi})$ at a
point $P^{*}(s)$ is expressed by
\begin{equation}
K_{1}(s) = \dfrac{ds^*}{ds}
\end{equation}

3. $C^{*}$ is reducible to a point iff $(C,\overline{\xi})$ is a
parallel versor field in $E_{3}$.

4. The tangent line at point $P^{*}\in C^{*}$ is parallel with
principal normal line of $(C,\overline{\xi})$.

5. Since $\overline{\xi}_{1}(s)\times \overline{\xi}_{2}(s) =
\overline{\xi}_{3}(s)$ it follows that the direction of binormal
versor field $\overline{\xi}_{3}$ is the direction tangent to
$\Sigma$ orthogonal to tangent line of $C^{*}$ at point $P^{*}$.

6. The geodesic curvature $\k_{g}$ of the curve $C^{*}$ at a point
$P^{*}$ verifies the equation
\begin{equation}
\k_{g} ds^{*} =K_{2}ds.
\end{equation}

7. The versor field $(C,\overline{\xi})$ is of null torsion (i.e.
$K_{2}(s)=0$) iff $\k_{g} = 0.$ In this case $C^{*}$ is an arc of
a great circle on $\Sigma$.

\section{Plane fields $(C,\pi)$}

\setcounter{theorem}{0}\setcounter{equation}{0}

A plane field $(C,\pi)$ is defined by the versor field
$(C,\overline{\nu}(s))$ where $\overline{\nu}(s)$ is normal to
$\pi(s)$ in every point $P(s)\in C.$ We assume the $\pi(s)$ is
oriented. Consequently, $\overline{\nu}(s)$ is well determined.

Let $\cal{R}_{\pi} = (P(s), \overline{\nu}_{1}(s),
\overline{\nu}_{2}(s), \overline{\nu}_{3}(s))$ be the Frenet frame,
$\overline{\nu}(s) = \overline{\nu}_{1}(s)$, of the versor field
$(C, \overline{\nu})$.

It follows that the fundamental equations of the plane field
$(C,\pi)$ are:
\begin{equation}
\dfrac{d\overline{r}}{ds}(s) = b_{1}(s){\overline{\nu}}_{1} +
b_{2}(s){\overline{\nu}}_{2} + b_{3}(s)\nu_{3},\;\;
b_{1}^{2}+b_{2}^{2}+b_{3}^{2}=1.
\end{equation}
\begin{eqnarray}
\dfrac{d\overline{\nu}_{1}}{ds} &=& \chi_{1}(s)\overline{\nu}_{2},\nonumber\\
\dfrac{d\overline{\nu}_{2}}{ds} &=& -\chi_{1}\overline{\nu}_{1} +
\chi_{2}(s)\overline{\nu}_{3}\\
\dfrac{d\overline{\nu}_{3}}{ds}&=&
-\chi_{2}(s)\overline{\nu}_{2}.\nonumber
\end{eqnarray}

The invariant $\chi_{1}(s) = \left\|\dfrac{d\overline{\nu}_{1}}{dr}\right\|$
is called the curvature and $\chi_{2}(s)$ is the torsion of the
plane field $(C,\pi(s))$.

The following properties hold:

1. The characteristic straight lines of the field of planes
$(C,\pi)$ cross through the corresponding point $P(s)\in C$ iff the
invariant $b_{1}(s) = 0$, $\forall s\in (s_{1},s_{2})$.

2. The planes $\pi(s)$ are parallel along the curve $C$ iff the
invariant $\chi_{1}(s)$ vanishes, $\forall s\in (s_{1},s_{2})$.

3. The characteristic lines of the plane field $(C,\pi)$ are
parallel iff the invariant $\chi_{2}(s) = 0$, $\forall s\in
(s_{1},s_{2})$.

4. The versor field $(C,\overline{\nu}_{3})$ determines the
directions of the characteristic line of $(C,\pi)$.

5. The versor field $(C,\overline{\nu}_{3})$, with
$\chi_{2}(s)\neq 0$, is concurrent iff:
$$
b_{1}(s) = 0,\;\; b_{3}(s) +
\dfrac{d}{dt}\left(\dfrac{b_{2}(s)}{\chi_{2}(s)}\right)=0.
$$

6. The curve $C$ is on orthogonal trajectory of the generatrices of
the ruled surface $\cal{R}(C,\overline{\nu}_{3})$ if $b_{3}(s) =
0,$ $\forall s\in (s_{1},s_{2})$.

7. By means of equations (1.3.1), (1.3.2) we can prove a
fundamental theorem for the plane field $(C,\pi).$

\chapter{Myller configurations $\mathfrak{M}(C,\overline{\xi},\pi)$}

The notions of versor field $(C,\overline{\xi})$ and the plane
field $(C,\pi)$ along to the same curve $C$ lead to a more general
concept named Myller configuration
$\mathfrak{M}(C,\overline{\xi},\pi)$, in which every versor
$\overline{\xi}(s)$ belongs to the corresponding plane $\pi(s)$ at
point $P(s)\in C.$ The geometry of $\mathfrak{M} =
\mathfrak{M}(C,\overline{\xi},\pi)$ is much more rich as the
geometries of $(C,\overline{\xi})$ and $(C,\pi)$ separately taken.
For $\mathfrak{M}$ one can define its geometric invariants, a Darboux frame
and introduce a new idea of parallelism or concurrence of versor
$(C,\overline{\xi})$ in $\mathfrak{M}$. The geometry of
$\mathfrak{M}$ is totally based on the fundamental equations of
$\mathfrak{M}$. The basic idea of this construction belongs to Al.
Myller [31], [32], [33], [34] and it was considerable developed by O. Mayer
[20], [21], R. Miron [23], [24], [62] (who proposed the name of Myller
Configuration and studied its complete system of invariants).

\section{Fundamental equations of Myller configuration}

\setcounter{theorem}{0}\setcounter{equation}{0}

\begin{definition}
A Myller configuration $\mathfrak{M} =
\mathfrak{M}(C,\overline{\xi},\pi)$ in the Euclidean space $E_{3}$
is a pair $(C,\overline{\xi})$, $(C,\pi)$ of versor field and
plane field, having the property: every $\overline{\xi}(s)$
belongs to the plane $\pi(s)$.
\end{definition}

Let $\overline{\nu}(s)$ be the normal versor to plane $\pi(s)$.
Evidently $\overline{\nu}(s)$ is uniquely determined if $\pi(s)$
is an oriented plane for $\forall s\in (s_{1},s_{2})$.

By means of versors $\overline{\xi}(s), \overline{\nu}(s)$ we can
determine the {\it Darboux frame} of $\mathfrak{M}:$
\begin{equation}
\cal{R}_{D} = (P(s); \overline{\xi}(s), \overline{\mu}(s),
\overline{\nu}(s)),
\end{equation}
where
\begin{equation}
\overline{\mu}(s) = \overline{\nu}(s)\times \overline{\xi}(s).
\end{equation}

$\cal{R}_{D}$ is geometrically associated to
$\mathfrak{M}$. It is orthonormal and positively oriented.

Since the versors $\overline{\xi}(s), \overline{\mu}(s),
\overline{\nu}(s)$ have a geometric meaning, the same pro\-per\-ties
have the vectors $\dfrac{d\overline{\xi}}{ds},$
$\dfrac{d\overline{\mu}}{ds}$ and $\dfrac{d\overline{\nu}}{ds}$.

Therefore, we can prove, without difficulties:

\begin{theorem}
The moving equations of the Darboux frame of $\mathfrak{M}$ are as
follows:
\begin{equation}
\dfrac{d\overline{r}}{ds} = \overline{\a}(s) =
c_{1}(s)\overline{\xi} + c_{2}(s)\overline{\mu} +
c_{3}(s)\overline{\nu};\;\; c_{1}^{2} + c_{2}^{2}+c_{3}^{2} =1
\end{equation}
and
\begin{eqnarray}
\dfrac{d\overline{\xi}}{ds} &=& G(s)\overline{\mu} +
K(s)\overline{\nu},\nonumber\\
\dfrac{d\overline{\mu}}{ds} &=& -G(s)\overline{\xi} +
T(s)\overline{\nu},\\\dfrac{d\overline{\nu}}{ds} &=&
-K(s)\overline{\xi} - T(s)\overline{\mu}\nonumber
\end{eqnarray}
and $c_{1}(s), c_{2}(s), c_{3}(s); G(s), K(s)$ and $T(s)$ are uniquely determined and are invariants.
\end{theorem}

The previous equations are called {\it the fundamental equations}
of the Myller configurations $\mathfrak{M}$.

Terms:

$G(s)$ -- is the geodesic curvature, $K(s)$ -- is the normal
curvature, $T(s)$ -- is the geodesic torsion of the versor field $(C,\overline{\xi})$ in Myller configuration $\mathfrak{M}$.

For $\mathfrak{M}$ a fundamental theorem can be stated:

\begin{theorem}
Let be a priori given $C^{\infty}$ functions $c_{1}(s), c_{2}(s),
c_{3}(s)$, $[c_{1}^{2} + c_{2}^{2} + c_{3}^{2}=1]$, $G(s),$
$K(s)$, $T(s)$, $s\in [a,b]$. Then there is a Myller configuration
$\mathfrak{M}(C,\overline{\xi},\pi)$ for which $s$ is the arclength of
curve $C$ and the given functions are its invariants. Two such
configuration differ by a proper Euclidean motion.
\end{theorem}

\noindent{\it Proof.}
By means of given functions $c_{1}(s), \ldots, G(s), \ldots$ we
can write the system of differential equations (2.1.3), (2.1.4).
Let the initial conditions $\overline{r}_{0} =
\overrightarrow{OP}_{0}$, $(\overline{\xi}_{0},
\overline{\mu}_{0}, \overline{\nu}_{0})$, an orthonormal,
positively oriented frame in $E_{3}$.

From (2.1.4) we find an unique solution $(\overline{\xi}(s),
\overline{\mu}(s),\overline{\nu}(s))$, $s\in [a,b]$ with the
property
$$
\overline{\xi}(s_{0}) = \overline{\xi}_{0},\;\;
\overline{\mu}(s_{0}) = \overline{\mu}_{0},\;\;
\overline{\nu}(s_{0}) = \overline{\nu}_{0},
$$
with $s_{0}\in [a,b]$ and $(\overline{\xi}(s), \overline{\mu}(s),
\overline{\nu} (s))$ being an orthonormal, positively oriented
frame.

Then consider the following solution of (2.1.3)
$$
\overline{r}(s) = \overline{r}_{0} +
\int_{s_{0}}^{s}[c_{1}(\s)\overline{\xi}(\s) +
c_{2}(\s)\overline{\mu}(\s) +c_{3}(\s)\overline{\nu}(\s)]d\s,
$$
which has the property $\overline{r}(s_{0}) = \overline{r}_{0}$,
and $\left\|\dfrac{d\overline{r}}{ds}\right\| =1.$

Thus $s$ is the arc length on the curve $\overline{r} = r(s)$.

Now, consider the configuration
$\mathfrak{M}(C,\overline{\xi},\pi(s))$, $\pi(s)$ being the plane
orthogonal to versor $\overline{\nu}(s)$ at point $P(s)$.

We can prove that $\mathfrak{M}(C,\overline{\xi},\pi)$ has as
invariants just $c_{1},c_{2},c_{3},G,K,T.$

The fact that two configurations $\mathfrak{M}$ and
$\mathfrak{M}^{\prime}$, obtained by changing the initial
conditions $(\overline{r}_{0}; \overline{\xi}_{0},
\overline{\mu}_{0}, \overline{\nu}_{0})$ to
$(\overline{r}_{0}^{\prime}, \overline{\xi}_{0}^{\prime},
\overline{\mu}_{0}^{\prime}, \overline{\nu}_{0}^{\prime})$ differ
by a proper Euclidian motions follows from the property that the
exists an unique Euclidean motion which apply $(\overline{r}_{0};
\overline{\xi}_{0}, \overline{\mu}_{0},\overline{\nu}_{0})$ to
$(\overline{r}^{\prime}_{0}, \overline{\xi}_{0}^{\prime},
\overline{\mu}_{0}^{\prime},
\overline{\nu}_{0}^{\prime})$.

\section{Geometric interpretations of invariants}

\setcounter{theorem}{0}
\setcounter{equation}{0}

The invariants $c_{1}(s), c_{2}(s), c_{3}(s)$ have simple
geometric interpretations:
$$
c_{1}(s) = \cos \sphericalangle (\overline{\a}, \overline{\xi}),\;
c_{2}(s) =\cos \sphericalangle (\overline{\a}, \overline{\mu}),\;
c_{3} = \cos \sphericalangle (\overline{\a}, \overline{\nu}).
$$

We can find some interpretation of the invariants $G(s), K(s)$ and
$T(s)$ considering a variation of Darboux frame
$$
R_{D}(P(s); \overline{\xi}(s),
\overline{\mu}(s),\overline{\nu}(s))$$$$\to
R^{\prime}_{D}(P^{\prime}(s+\D s), \overline{\xi}(s+\D s),
\overline{\mu}(s+\D s), \overline{\nu}(s+\D s)),
$$
obtained by the Taylor expansion
\begin{equation}
\begin{array}{lll}
\overline{r}(s+\D s)&=& \overline{r}(s)+\dfrac{\D s
}{1!}\dfrac{d\overline{r}}{ds} + \dfrac{(\D s)^{2}}{2!}
\dfrac{d^{2}r}{ds^{2}} + \ldots +\vspace*{2mm}\\ & & +\dfrac{(\D
s)^{n}}{n!}\left(\dfrac{d^{n}\overline{r}}{ds^{n}}
+ \overline{\o}_{0}(s,\D
s)\right),\\\overline{\xi}(s+\D s)&=&
\overline{\xi}(s)+\dfrac{\D s }{1!}\dfrac{d\overline{\xi}}{ds} +
\dfrac{(\D s)^{2}}{2!} \dfrac{d^{2}\overline{\xi}}{ds^{2}} + \ldots
+\\ & & +\dfrac{(\D s)^{n}}{n!}\left(\dfrac{d^{n}\overline{\xi}}{ds^{n}}
+ \overline{\o}_{1}(s,\D
s)\right),\vspace*{2mm}\\\overline{\mu}(s+\D s)&=&
\overline{\mu}(s)+\dfrac{\D s }{1!}\dfrac{d\overline{\mu}}{ds} +
\dfrac{(\D s)^{2}}{2!} \dfrac{d^{2}\overline{\mu}}{ds^{2}} + \ldots
+\\ & & +\dfrac{(\D
s)^{n}}{n!}\left(\dfrac{d^{n}\overline{\mu}}{ds^{n}}+ \overline{\o}_{2}(s,\D
s)\right),\vspace*{2mm}\\\overline{\nu}(s+\D s)&=& \overline{\nu}(s)+\dfrac{\D s
}{1!}\dfrac{d\overline{\nu}}{ds} + \dfrac{(\D s)^{2}}{2!}
\dfrac{d^{2}\overline{\nu}}{ds^{2}} + \ldots +\vspace*{2mm}\\ & &\dfrac{(\D
s)^{n}}{n!}\left(\dfrac{d^{n}\overline{\nu}}{ds^{n}}+ \overline{\o}_{3}(s,\D s)\right),
\end{array}
\end{equation}
where
\begin{eqnarray}
\lim\limits_{\D s\to 0}\overline{\o}_{i}(s,\D s) = 0,\;\;
(i=0,1,2,3).
\end{eqnarray}

Using the fundamental formulas (2.1.3), (2.1.4) we can write, for
$n=1:$
\begin{eqnarray}
\overline{r}(s+\D s) &=& \overline{r}(s) + \D s
(c_{1}\overline{\xi} + c_{2} \overline{\mu} + c_{3}\overline{\nu}
+ \overline{\o}_{0}(s,\D s)),\nonumber\\\overline{\xi}(s+\D s) &=&
\overline{\xi}(s) + \D s(G(s) \overline{\mu} + K(s)\overline{\nu}
+ \overline{\o}_{1}(s,\D s) ),\\\overline{\mu}(s+\D s) &=&
\overline{\mu}(s) + \D s (-G(s)\overline{\xi} + T(s)\overline{\nu}
+ \overline{\o}_{2}(s,\D s) )\nonumber\\\overline{\nu}(s)+\D s &=&
\overline{\nu}(s) + \D s(-K(s)\overline{\xi} - T(s)\overline{\mu}
+ \overline{\o}_{3}(s,\D s)).\nonumber
\end{eqnarray} and (2.2.2) being verified.

Let $\overline{\xi}^{*}(s+\D s)$ be the orthogonal projection of
the versor $\overline{\xi}(s+\D s)$ on the plane $\pi(s)$ at point
$P(s)$ and let $\D \psi_{1}$ be the oriented angle of the versors
$\overline{\xi}(s), \overline{\xi}^{*}(s+\D s)$. Thus, we have

\begin{theorem}
The invariant $G(s)$ of versor field $(C,\xi)$ in Myller
configuration $\mathfrak{M}(C,\overline{\xi},\pi)$ is given by
$$
G(s) = \lim_{\D s\to 0}\dfrac{\D\psi_{1}}{\D s}.
$$
\end{theorem}

By means of second formula (2.2.3), this Theorem can be proved
without difficulties.

Therefore the name of {\it geodesic curvature} of $(C,\overline{\xi})$
in $\mathfrak{M}$ is justified.

Consider the plane $(P(s);
\overline{\xi}(s), \overline{\nu}(s))$-called the normal plan of
$\mathfrak{M}$ which contains the versor $\overline{\xi}(s)$.

Let be the vector $\overline{\xi}^{**}(s+\D s)$ the orthogonal
projection of versor $\overline{\xi}(s+\D s)$ on the normal plan
$(P(s); \overline{\xi}(s), \overline{\nu}(s))$. The angle $\D
\psi_{2} = \sphericalangle (\overline{\xi}(s),
\overline{\xi}^{**}(s+\D s))$ is given by the forma
$$
\sin \D \psi_{2} = \dfrac{\langle\overline{\xi}(s),
\overline{\xi}(s+\D s),
\overline{\mu}(s)\rangle}{\|\overline{\xi}^{**}(s+\D s)\|}.
$$
By (2.2.3) we obtain
$$
\sin \D \psi_{2} = \dfrac{K(s) +\langle \overline{\xi}(s),
\overline{\o}_{1}(s, \D s), \overline{\mu}(s)
\rangle}{\|\overline{\xi}^{**}(s+\D s)\|}\D s.
$$
Consequently, we have:

\begin{theorem}
The invariant $K(s)$ has the following geometric interpretation
$$
K(s) = \lim_{\D s\to 0}\dfrac{\D \psi_{2}}{\D s}.
$$\end{theorem}

Based on the previous result we can call {\it $K(s)$ the normal
curvature} of $(C,\overline{\xi})$ in $\mathfrak{M}$.

A similar interpretation can be done for the invariant $T(s)$.

\begin{theorem}
The function $T(s)$ has the interpretation:
$$
T(s) = \lim\limits_{\D s\to 0}\dfrac{\D \psi_{3}}{\D s},
$$
where $\D \psi_{3}$ is the oriented angle between
$\overline{\mu}(s)$ and $\overline{\mu}^{*}(s+\D s)$-which is the
orthogonal projection of $\overline{\mu}(s+\D s)$ on the normal
plane $(P(s); \overline{\mu}(s),$
$\overline{\nu}(s))$.\end{theorem}

This geometric interpretation allows to give the name {\it
geodesic torsion} for the invariant $T(s)$.

\section{The calculus of invariants $G, K, T$}
\setcounter{theorem}{0}
\setcounter{equation}{0}

The fundamental formulae (2.1.3), (2.1.4) allow to calculate the
expressions of second derivatives of the versors of Darboux frame
$\cal{R}_{D}$. We have:
\begin{equation}\hspace*{8mm}
\begin{array}{lll}
\dfrac{d^{2}\overline{r}}{ds^{2}}&=& \left(\dfrac{dc_{1}}{ds} -
Gc_{2}-K c_{3} \right)\overline{\xi}+\left(\dfrac{dc_{2}}{ds} +
Gc_{1} - Tc_{3} \right)\overline{\mu} +\vspace*{2mm}\\&&+\left(\dfrac{dc_{3}}{ds}
+ K c_{1} + Tc_{2}\right)\overline{\nu}
\end{array}
\end{equation}
and
\begin{equation}
\hspace*{8mm}\begin{array}{lll}
\dfrac{d^{2}\overline{\xi}}{ds^{2}} &=& -(G^{2} +
T^{2})\overline{\xi} + \left(\dfrac{dG}{ds} -
KT\right)\overline{\mu} +\left(\dfrac{dK}{ds} +
GT\right)\overline{\nu},\vspace*{2mm}\\ \dfrac{d^{2}\overline{\mu}}{ds^{2}}
&=& -\left(\dfrac{dG}{ds} +KT \right)\overline{\xi} -
(G^{2}+T^{2})\overline{\mu} + \left(\dfrac{dT}{ds} -
GT\right)\overline{\nu},\vspace*{2mm}\\ \dfrac{d^{2}\nu}{ds^{2}} &=&
\left(-\dfrac{dK}{ds}+GT\right)\overline{\xi} - \left(\dfrac{dT}{ds}
+GT \right)\overline{\mu} - (K^{2}+T^{2})\overline{\nu}.
\end{array}
\end{equation}
These formulae will be useful in the next part of the book.

From the fundamental equations (2.1.3), (2.1.4) we get

\begin{theorem}
The following formulae for invariants $G(s), K(s)$ and $T(s)$
hold:
\begin{equation}
G(s) = \left\langle \overline{\xi},
\dfrac{d\xi}{ds},\overline{\nu}\right\rangle,
\end{equation}
\begin{equation}
K(s) = \left\langle \dfrac{d\overline{\xi}}{ds},
\overline{\nu}\right\rangle = -\left\langle \overline{\xi},
\dfrac{d\overline{\nu}}{ds}\right\rangle, T(s)=\left\langle\overline{\xi},\overline{\nu},\dfrac{d\overline{\nu}}{ds}\right\rangle.
\end{equation}
\end{theorem}
Evidently, these formulae hold in the case when $s$ is the arclength of
the curve $C$.

\section{Relations between the invariants of the field $(C,\overline{\xi})$ and the invariants of $(C,\overline{\xi})$ in
$\mathfrak{M}(C,\overline{\xi},\pi)$}
\setcounter{theorem}{0}
\setcounter{equation}{0}

The versors field $(C,\overline{\xi})$ in $E_{3}$ has a Frenet
frame $\cal{R}_{F} = (P(s), \overline{\xi}_{1},
\overline{\xi}_{2}, \overline{\xi}_{3})$ and a complete system of
invariants $(a_{1},a_{2},a_{3}; K_{1}, K_{2})$ verifying the
equations (1.1.5) and (1.1.6).

The same field $(C,\overline{\xi})$ in Myller configuration
$\mathfrak{M}(C,\overline{\xi},\pi)$ has a Darboux frame
$\cal{R}_{D} = (P(s), \overline{\xi}, \overline{\mu},
\overline{\nu})$ and a complete system of invariants $(c_{1},
c_{2}, c_{3}; G,K,T)$. If we relate $\cal{R}_{F}$ to $\cal{R}_{D}$
we obtain
\begin{eqnarray*}
\overline{\xi}_{1}(s) &=&
\overline{\xi}(s),\\\overline{\xi}_{2}(s)&=& \overline{\mu}(s)\sin
\v +\overline{\nu}(s)\cos \v,\\\overline{\xi}_{3}(s) &=&-\mu(s)cos
\v + \overline{\nu}(s)\sin \v,
\end{eqnarray*}
with $\v = \sphericalangle (\overline{\xi}_{2},\overline{\nu})$
and
$\langle\overline{\xi}_{1},\overline{\xi}_{2},\overline{\xi}_{3}\rangle
= \langle\overline{\xi},\overline{\mu},\overline{\nu}\rangle =1.$

Then, from (1.1.5) and (2.1.3) we can determine the relations between the two
systems of invariants.

From
$$
\dfrac{d\overline{r}}{ds} = \overline{\a}(s) =
a_{1}\overline{\xi}_{1} + a_{2}\overline{\xi}_{2} +
a_{3}\overline{\xi}_{3} = c_{1}\overline{\xi} +
c_{2}\overline{\mu}+c_{3}\overline{\nu}
$$
it follows
\begin{eqnarray}
c_{1}(s)&=&a_{1}(s)\nonumber\\c_{2}(s) &= &a_{1}(s)\sin \v
-a_{3}(s)\cos \v\\c_{3}(s)&=& a_{2}(s)\cos \v + a_{3}(s)\sin
\v.\nonumber
\end{eqnarray}
And, (1.1.6), (2.1.4) we obtain
\begin{eqnarray}
G(s)&=& K_{1}(s) \sin \v\nonumber\\K(s) &=& K_{1}(s)\cos \v\\T(s)
&=& K_{2}(s) + \dfrac{d\v}{ds}.\nonumber
\end{eqnarray}

These formulae, allow to investigate some important properties of
$(C,\overline{\xi})$ in $\mathfrak{M}$ when some invariants
$G,K,T$ vanish.

\section{Relations between invariants of normal field $(C,\overline{\nu})$ and invariants $G,K,T$}
\setcounter{theorem}{0}
\setcounter{equation}{0}

The plane field $(C,\pi)$ is characterized by the normal versor
field $(C,\overline{\nu})$, which has as Frenet frame $\cal{R}_{F}
= (P(s); \overline{\nu}_{1}, \overline{\nu}_{2},
\overline{\nu}_{3})$ with $\overline{\nu}_{1} = \overline{\nu}$
and has $(b_{1}, b_{2}, b_{3}, \chi_{1},\chi_{2})$ as a complete
system of invariants. They satisfy the formulae (1.3.1), (1.3.2).
But the frame $\cal{R}_{F}$ is related to Darboux frame
$\cal{R}_{D}$ of $(C,\overline{\xi})$ in $\mathfrak{M}$ by the
formulae
\begin{eqnarray}
\overline{\nu}_{1} &=&
\overline{\nu}(s)\nonumber\\-\overline{\nu}_{2}&=& \sin \s
\overline{\xi} +\cos \s \overline{\mu}\\\overline{\nu}_{3}&=&
-\cos \s \overline{\xi}+\sin \s \overline{\mu}\nonumber
\end{eqnarray}
where $\s = \sphericalangle (\overline{\xi}(s),
\overline{\nu}_{3}(s))$.

Proceeding as in the previous section we deduce

\begin{theorem}
The following relations hold:
\begin{eqnarray}
c_{1}&=& -b_{2}\sin \s + b_{3}\cos \s\nonumber\\-c_{2}&=&
b_{2}\cos \s + b_{3}\sin \s\\c_{3}&=& b_{2}\nonumber
\end{eqnarray}
and
\begin{eqnarray}
K&=& \chi_{1}\sin \s\nonumber\\T&=& \chi_{1}\cos \s\\G&=&
\chi_{2}+\dfrac{d\s}{ds}.\nonumber
\end{eqnarray}
\end{theorem}

A first consequence of previous formulae is given by

\begin{theorem}
The invariant $K^{2}+T^{2}$ depends only on the plane field
$(C,\pi)$. We have
\begin{equation}
K^{2}+T^{2} = \chi_{1}.
\end{equation}
\end{theorem}

The proof is immediate, from (2.5.3).

\section{Meusnier$^{\prime}$s theorem. Versor fields $(C,\overline{\xi})$ conjugated with tangent versor $(C,\overline{\a})$}
\setcounter{theorem}{0}
\setcounter{equation}{0}

Consider the vector field $\overline{\xi}^{**}(s+\D s)$, $(|\D
s|<\vp, \vp>0)$, the orthogonal projection of versor
$\overline{\xi}(s+\D s)$ on the normal plane $(P(s);
\overline{\xi}(s), \overline{\nu}(s))$. Since, up to terms of
second order in $\D s$, we have
$$
\overline{\xi}^{**}(s+\D s) = \overline{\xi}(s) + \dfrac{\D
s}{1!}K(s)\overline{\nu}(s) + \overline{\t}(s,\D s)\dfrac{(\D
s)^{2}}{2!},
$$
for $\D s\to 0$ one gets:
\begin{equation}
\dfrac{d\overline{\xi}^{**}}{ds} = K(s)\overline{\nu}(s).
\end{equation}

Assuming $K(s)\neq 0$ we consider the point $P_{c}^{**}$-called
the {\it center of curvature} of the vector field
$(C,\overline{\xi}^{**})$, given by
$$
\overrightarrow{PP}_{c}^{**} = \dfrac{1}{K(s)}\overline{\nu}(s).
$$
On the other hand the field of versors $(C,\overline{\xi})$ have a
center of curvature $P_{c}$ given by $\overrightarrow{PP}_{c} =
\dfrac{1}{K_{1}(s)}\overline{\xi}_{2}$.

The formula (2.4.2), i.e., $K(s) = K_{1}(s)\cos \v,$ shows that
the orthogonal projection of vector $\overrightarrow{PP}_{c}^{**}$
on the (osculating) plane $(P(s); \overline{\xi}_{1},
\overline{\xi}_{2})$ is the vector $\overrightarrow{PP}_{c}$.

Indeed, we have
\begin{equation}
\dfrac{\cos \v}{K} = \dfrac{1}{K_{1}}
\end{equation}
As a consequence we obtain  a theorem of Meusnier type:

\begin{theorem}
The curvature center of the field $(C,\overline{\xi})$ in
$\mathfrak{M}$ is the orthogonal projection on the osculating
plane $(M; \overline{\xi}_{1}, \overline{\xi}_{2})$ of the
curvature center $P_{c}^{**}$.
\end{theorem}

\begin{definition}
The versor field $(C,\overline{\xi})$ is called conjugated with
tangent versor field $(C,\overline{\a})$ in the Myller
configuration $\mathfrak{M}(C,\overline{\xi},\pi)$ if the
invariant $K(s)$ vanishes.
\end{definition}

Some immediate consequences:

1. $(C,\overline{\xi})$ is conjugated with $(C,\overline{\a})$ in
$\mathfrak{M}$ iff the line $(P;\overline{\xi})$ is parallel in
$E_{3}$ to the characteristic line of the planes $\pi(s)$, $s\in
(s_{1},s_{2})$.

2. $(C,\overline{\xi})$ is conjugated with $(C,\overline{\a})$ in
$\mathfrak{M}$ iff $|T(s)| = \chi_{1}(s)$.

3. $(C,\overline{\xi})$ is conjugated with $(C,\overline{\a})$ in
$\mathfrak{M}$ iff the osculating planes
$(P;\overline{\xi}_{1},\overline{\xi}_{2})$ coincide to the planes
$\pi(s)$ of $\mathfrak{M}$.

4. $(C,\overline{\xi})$ is conjugated with $(C,\overline{\xi})$
iff the asimptotic planes of the ruled surface
$\cal{R}(C,\overline{\xi})$ coincide with the planes $\pi(s)$ of
$\mathfrak{M}$.

\section{Versor field $(C,\overline{\xi})$ with null geodesic torsion}
\setcounter{theorem}{0}
\setcounter{equation}{0}

A new relation of conjugation of versor field $(C,\overline{\xi})$
with the tangent versor field $(C,\overline{\a})$ is obtained in
the case $T(s) = 0.$

\begin{definition}
The versor field $(C,\overline{\xi})$ is called orthogonal
conjugated with the tangent versor field $(C,\overline{\a})$ in
$\mathfrak{M}$ if its geodesic torsion $T(s) = 0,$ $\forall s\in
(s_{1},s_{2})$.

Some properties
\begin{itemize}
\item[1.] $(C,\overline{\xi})$ is orthogonal conjugated with
$(C,\overline{\a})$ in $\mathfrak{M}$ iff $\overline{\mu}(s)$ are
parallel with the characteristic line of planes $\pi(s)$ along the
curve $C$.

\item[2.] $(C,\overline{\xi})$ is orthogonal conjugated with
$(C,\overline{\a})$ in $\mathfrak{M}$ if $|K_{1}(s)| =
\chi_{1}(s)$, along $C$.
\end{itemize}
\end{definition}

\begin{theorem}
Assuming that the versor field $(C,\overline{\xi})$ in the
configuration $\mathfrak{M}(C,\overline{\xi},\pi)$ has two of the
following three properties then it has the third one, too:
\begin{itemize}
\item[a.] The osculating planes $(P; \overline{\xi}_{1},
\overline{\xi}_{2})$ are parallel in $E_{3}$ along $C$.

\item[b.] The osculating planes
$(P;\overline{\xi}_{1},\overline{\xi}_{2})$ have constant angle
with the plans $\pi(s)$ on $C$.

\item[c.] The geodesic torsion $T(s)$ vanishes on $C$.
\end{itemize}
\end{theorem}

The proof is based on the formula $T(s) = K_{2}(s) +
\dfrac{d\v}{ds}$, $\v = \sphericalangle
(\overline{\xi}_{2},\overline{\nu})$.

Consider two Myller configurations
$\mathfrak{M}(C,\overline{\xi},\pi)$ and
$\mathfrak{M}^{\prime}(C,\overline{\xi},\pi^{\prime})$ which have
in common the versor field $(C,\overline{\xi})$. Denote by $\v =
\sphericalangle (\overline{\xi}_{2},\overline{\nu})$, $\v^{\prime}
= \sphericalangle(\overline{\xi}_{2},\overline{\nu}^{\prime})$.
Then the geodesic torsions of $(C,\overline{\xi})$ in
$\mathfrak{M}$ and $\mathfrak{M}^{\prime}$ are follows:
$$
T(s) = K_{2}(s) + \dfrac{d\v}{ds},\;\; T^{\prime}(s) = K_{2}(s) +
\dfrac{d\v^{\prime}}{ds}.
$$

Evidently, we have $\v-\v^{\prime} = \sphericalangle
(\overline{\nu},\overline{\nu}^{\prime})$.

By means of these relations we can prove, without difficulties:

\begin{theorem}
If the Myller configurations $\mathfrak{M}(C,\overline{\xi},\pi)$
and $\mathfrak{M}^{\prime}(C,\overline{\xi},\pi^{\prime})$ have two
of the following properties:

\begin{itemize}
\item[a)] $(C,\overline{\xi})$ has the null geodesic torsion,
$T(s)=0$ in $\mathfrak{M}$.

\item[b)] $(C,\overline{\xi})$ has the null geodesic torsion
$T^{\prime}(s)$in $\mathfrak{M}^{\prime}$.

\item[c)] The angle
$\sphericalangle(\overline{\nu},\overline{\nu}^{\prime})$is
constant along $C$, then $\mathfrak{M}$ and
$\mathfrak{M}^{\prime}$ have the third property.
\end{itemize}
\end{theorem}

\begin{remark}\rm
The versor field $(C,\overline{\nu}_{2})$ is orthogonally conjugated
with the tangent versors $(C,\overline{\a})$ in the configuration
$\mathfrak{M}(C,\overline{\xi},\pi)$.
\end{remark}

\section{The vector field parallel in Myller sense in configurations $\mathfrak{M}$}
\setcounter{theorem}{0}
\setcounter{equation}{0}

Consider $(C,\overline{V})$ a vector field, along the curve $C$.
We denote $\overline{V}(s) = \overline{V}(\overline{r}(s))$ and
say that $(C,\overline{V})$ is a vector field in the configuration
$\mathfrak{M} = \mathfrak{M}(C,\overline{\xi},\pi)$ if the vector
$\overline{V}(s)$ belongs to plane $\pi(s),$ $\forall s\in
(s_{1},s_{2})$.

\begin{definition}
The vectors field $(C,\overline{V})$ in
$\mathfrak{M}(C,\overline{\xi},\pi)$ is parallel in Myller sense
if the vector field $\dfrac{d\overline{V}}{ds}$ is normal to
$\mathfrak{M}$, i.e. $\dfrac{d\overline{V}}{ds} =
\l(s)\overline{\nu}(s),$ $\forall s\in
(s_{1},s_{2})$.\end{definition}

The parallelism in Myller sense is a direct generalization of
Levi-Civita parallelism of tangent vector fields along a curve
$C$ of a surface $S$.

It is not difficult to prove that the vector field
$\overline{V}(s)$ is parallel in Myller sense if the vector field
$$
\overline{V}^{\prime}(s+\D s) = pr_{\pi(s)}\overline{V}(s+\D s)
$$
is parallel in ordinary sens in $E_{3}$-up to terms of second
order in $\D s$.

In Darboux frame, $\overline{V}(s)$ can be represented by its
coordinate as follows:
\begin{equation}
\overline{V}(s) =V^{1}(s)\overline{\xi}(s) +
V^{2}(s)\overline{\mu}(s).
\end{equation}

By virtue of fundamental equations (2.1.4) we find:
\begin{equation}
\hspace*{8mm}\dfrac{d\overline{V}}{ds} = \left(\dfrac{dV^{1}}{ds} - G
V^{2}\right)\overline{\xi} +\left(\dfrac{dV^{2}}{ds} +
GV^{1}\right)\overline{\mu} + (KV^{1} +TV^{2})\overline{\nu}.
\end{equation}
Taking into account the Definition 2.8.1, one proves:

\begin{theorem}
The vector field $\overline{V}(s)$, $(2.8.1)$ is parallel in Myller
sense in configuration $\mathfrak{M}(C, \overline{\xi}, \pi)$ iff
coordinates $V^{1}(s), V^{2}(s)$ are solutions of the system of
differential equations:
\begin{equation}
\dfrac{dV^{1}}{ds}-GV^{2} = 0,\;\; \dfrac{dV^{2}}{ds}+GV^{1} = 0.
\end{equation}
\end{theorem}

In particular, for $\overline{V}(s) = \overline{\xi}(s)$, we
obtain

\begin{theorem}
The versor field $\overline{\xi}(s)$ is parallel in Myller sense
in $\mathfrak{M}(C,\overline{\xi},\pi)$ iff the geodesic curvature
$G(s)$ of $(C,\overline{\xi})$ in $\mathfrak{M}$ vanishes.
\end{theorem}

This is a reason that Al. Myller says that $G$ is {\it the
deviation of parallelism} [31]. Later we will see that $G(s)$ is
an intrinsec invariant in the geometry of surfaces in $E_{3}$.

By means of (2.8.3) we have

\begin{theorem}
There exists an unique vector field $\overline{V}(s)$, $s\in
(s_{1}^{\prime}, s_{2}^{\prime})\subset (s_{1},s_{2})$ parallel in
Myller sense in the configuration
$\mathfrak{M}(C,\overline{\xi},\pi)$ which satisfy the initial
condition $\overline{V}(s_{0}) = \overline{V}_{0}$, $s_{0}\in
(s_{1}^{\prime}, s_{2}^{\prime})$ and $\langle \overline{V}_{0},
\overline{\nu}(s_{0}) \rangle = 0.$
\end{theorem}

Evidently, theorem of existence and uniqueness of solutions of
system (2.8.3), is applied in this case.

In particular, if $G(s)$ = constant, then the general solutions of
(2.8.3) can be obtained by algebric operations.

An important property of parallelism in Myller sense is expressed
in the next theorem.

\begin{theorem}
The Myller parallelism of vectors in $\mathfrak{M}$ preserves the
lengths and angles of vectors.
\end{theorem}

\begin{proof}\rm If $\dfrac{d\overline{V}}{ds} =
\l(s)\overline{\nu}$, then $\dfrac{d}{ds}\langle \overline{V},
\overline{V} \rangle=0.$ Also,
$\dfrac{d\overline{V}^{\prime}}{ds}=\l(s)\overline{\nu}$,
$\dfrac{d\overline{U}}{ds} = \l^{\prime}(s)\overline{\nu}$, then
$\dfrac{d}{ds}\langle \overline{V}(s),
\overline{U}^{\prime}(s)\rangle = 0$
\end{proof}

\section{Adjoint point, adjoint curve and concurrence in Myller sense}
\setcounter{theorem}{0}
\setcounter{equation}{0}

The notions of adjoint point, adjoint curve and concurrence in
Myller sense in a configuration $\mathfrak{M}$ have been
introduced and studied by O. Mayer [20] and Gh. Gheorghiev
[11], [55]. They applied these notions, to the theory of surfaces,
nonholomorphic manifolds and in the geometry of versor fields in
Euclidean space $E_{3}$.

In the present book we introduce these notions in a different way.

Consider the vector field
$$
\overline{\xi}^{*}(s+\D s) = pr_{\pi(s)}\overline{\xi}(s+\D s).
$$
Taking into account the formula (2.2.1)$^{\prime}$ we can write up
to terms of second order in $\D s:$
\begin{equation}
\overline{\xi}^{*}(s+\D s) = \overline{\xi}(s) + \D s  (G
\overline{\mu}(s) + \overline{\o}^{*}(s,\D s))
\end{equation}
with
$$
\overline{\o}^{*}(s, \D s) \to 0, (\D s\to 0).
$$

Let $C^{\prime}$ be the orthogonal projection of the curve $C$ on
the plane $\pi(s)$. A neighbor point $P^{\prime}(s+\D s)$ is
projected on plane $\pi(s)$ in the point $P^{*}(s+\D s)$ given by
\begin{equation}
\begin{array}{l}
\overline{r}^{*}(s+\D s) = \overline{r}(s) + \D s
(c_{1}\overline{\xi} + c_{2}\overline{\mu} +
\overline{\o}_{0}^{*}(s, \D s)),\\ \overline{\o}_{0}^{*}(s, \D s)
\to 0, (\D s\to 0).
\end{array}
\end{equation}

\begin{definition}
The adjoint point of the point $P(s)$ with respect to
$\overline{\xi}(s)$ in $\mathfrak{M}$ is the characteristic point
$P_{a}$ on the line $(P;\overline{\xi})$ of the plane ruled
surface $R(C^{*}, \overline{\xi}^{*})$.
\end{definition}

One proves that the position vector $\overline{R}(s)$ of adjoint
point $P_{a}$ for $G\neq 0$, is as follows:
\begin{equation}
\overline{R}(s) = \overline{r}(s) -
\dfrac{c_{2}}{G}\overline{\xi}(s).
\end{equation}
The vector field $(C^{*}, \overline{\xi}^{*})$ from (2.9.3) is
called {\it geodesic field}. A result established by O. Mayer [20]
holds:

\begin{theorem}
If the versor field $(C, \overline{\xi})$ is enveloping in space
$E_{3}$, then the adjoint point $P_{a}$ of the point $P(s)$ in
$\mathfrak{M}$ is the contact point of the line $(P,
\overline{\xi})$ with the cuspidale line.
\end{theorem}

\begin{definition}
The geometric locus of the adjoint points corresponding to the
versor field $(C, \overline{\xi})$ in $\mathfrak{M}$ is the
adjoint curve $C_{a}$ of the curve $C$ in $\mathfrak{M}$.
\end{definition}

The adjoint curve $C_{a}$ has the vector equations (2.9.3) for
$\forall s\in (s_{1},s_{2})$.

Now, we can introduce

\begin{definition}
The versor field $(C, \overline{\xi})$ is concurrent in Myller
sense in $\mathfrak{M}(C, \overline{\xi}, \pi)$ if, at every point
$P(s)\in C$ the geodesic vector field $(C^{*},
\overline{\xi}^{*})$ is concurrent.
\end{definition}

For $G(s)\neq 0$, we have

\begin{theorem}
The versor field $(C,\xi)$ is concurrent in Myller sense in
$\mathfrak{M}$ iff the following equation hold
\begin{equation}
\dfrac{d}{ds}\left(\dfrac{c_{2}}{G}\right) = c_{1}.
\end{equation}
\end{theorem}

For the proof see Section 1.1, Chapter 1, Theorem 1.1.3.

\section{Spherical image of a configuration $\mathfrak{M}$}
\setcounter{theorem}{0}
\setcounter{equation}{0}

In the Section 2, Chapter 1 we defined the spherical image of a
versor field $(C,\overline{\xi})$. Applying this idea to the
normal vectors field $(C,\overline{\nu})$ to a Myller
configuration $\mathfrak{M}(C, \overline{\xi}, \pi)$ we define the
notion of spherical image $C^{*}$ of $\mathfrak{M}$ as being
\begin{equation}
\overline{\nu}^{*}(s) = \overline{OP}^{*}(s) = \overline{\nu}(s)
\end{equation}
Thus, the relations between the curvature $\chi_{1}$ and torsion
$\chi_{2}$ of $(C,\overline{\nu})$ and geodesic curvature
$\k_{g}^{*}$ of $C^{*}$ at a point $P^{*}\in C^{*}$ and arclength
$s^{*}$ are as follows
$$
\chi_{1} = \dfrac{ds^{*}}{ds},\; \chi_{2}ds = \k_{g}^{*}ds^{*}.
$$
The properties enumerated in Section 2, ch 1, can be obtained for the
spheric image $C^{*}$ of configuration $\mathfrak{M}$.

Consider the versor field $\overrightarrow{P^{*}P^{*}_{1}} =
\xi^{*}(s) = \xi(s)$ and the angle $\t = \sphericalangle
(\overline{\nu}_{2}, \overline{\xi})$. It follows
\begin{equation}
K = -\dfrac{ds^{*}}{ds}\cos \t,\;\; (for \chi_{2}\neq 0).
\end{equation}

Thus, for $\t = \pm \dfrac{\pi}{2}$ one gets $K=0,$ which leads to
a new interpretation of the fact that $\overline{\xi}(s)$ is
conjugated to $\overline{\a}(s)$ in Myller configurations.

Analogous one can obtain $T = -\dfrac{ds^{*}}{ds}\cos
\widetilde{\t},$ $\widetilde{\t} =
\sphericalangle(\overline{\nu}_{2}, \overline{\mu}^{*})$ (with
$\overline{\mu}^{*}(s) = \overline{\mu}(s)$, applied at the point
$P^{*}(s)\in C^{*}$). For $\widetilde{\t} = \pm \dfrac{\pi}{2}$ it
follows that $\overline{\xi}(s)$ are orthogonally conjugated with
$\overline{\a}(s) .$

The problem is to see if the Gauss-Bomet formula can be extended
to Myller configurations $\mathfrak{M}(C, \overline{\xi}, \pi)$.
In the case $(\a(s)\in \pi(s)$, (i.e. $\overline{\a}(s)\bot
\overline{\nu}(s)$) such a problem was suggested by Thomson [18] and it
has been solved by Mark Krein in 1926, [18].

Here, we study this problem in the general case of Myller
configuration, when $\langle \overline{\a}(s),
\overline{\nu}(s)\rangle\neq 0.$

First of all we prove

\begin{lemma} Assume that we have:

1. $\mathfrak{M}(C, \overline{\xi},\pi)$ a Myller configuration of
class $C^{3}$, in which $C$ is a closed curve, having $s$ as
arclength.

2. The spherical image $C^{*}$ of $\mathfrak{M}$ determines on the
support sphere $\Sigma$ a simply connected domain of area $\o$.

In this hypothesis we have the formula
\begin{equation}
\o = 2\pi - \int_{C}\left(\overline{\nu}, \dfrac{d\overline{\nu}}{ds},
\dfrac{d^{2}\overline{\nu}}{ds^{2}}\right)/
\left\|\dfrac{d\overline{\nu}}{ds}\right\|^{2}ds.
\end{equation}
\end{lemma}

\begin{proof}\rm Let $\Sigma$ be the unitary sphere of
center $O\in E_{3}$ and a simply connected domain $D$, delimited
by $C^{*}$ on $\Sigma.$ Assume that $D$ remains to left with
respect to an observer looking in the sense of versor
$\overline{\nu}(s)$, when he is going along $C^{*}$ in the
positive sense.

Thus, we can take the following representation of $\Sigma$:
\begin{eqnarray}
x^{1}&=&\cos \v \sin \t\nonumber\\
x^{2}&=&\sin \v \sin \t\\x^{3}&=&\cos \v,\;\; \v\in [0, 2\pi), \t
\in \(-\dfrac{\pi}{2}, \dfrac{\pi}{2}\).\nonumber
\end{eqnarray}
The curve $C^{*}$ can be given by
\begin{equation}
\v = \v(s),\; \t = \t(s),\;\; s\in [0,s_{1}]
\end{equation}
with $\v(s), \t(s)$ of class $C^{3}$ and $C^{*}$ being closed:
$\v(0) = \v(s_{1})$, $\t(0) = \t(s_1)$.

The area $\o$ of the domain $D$ is
\begin{equation}
\begin{array}{l}
\o = \int_{0}^{s_{1}}\int_{0}^{\t(s_{1})}\sin \t d\t =
\int_{0}^{s_{1}}(1-\cos \t)\dfrac{d\v}{ds}ds =\\ = 2\pi
-\int_{0}^{s_{1}}\cos \t \dfrac{d\v}{ds}ds.
\end{array}
\end{equation}
Noticing that the versor $\overline{\nu}^{*} = \overline{\nu}(s)$
has the coordinate (2.10.4) a straightforward calculus leads to
$$
\left\langle \overline{\nu}, \dfrac{d\overline{\nu}}{ds},
\dfrac{d^{2}\overline{\nu}}{ds^{2}} \right\rangle/
\left\|\dfrac{d\overline{\nu}}{ds}\right\|^{2} = \cos \t
\dfrac{d\v}{ds} +\dfrac{d}{ds} \mbox{arctg}\, \dfrac{\sin
\t\dfrac{d\v}{ds}}{\dfrac{d\t}{ds}}.
$$
Denoting by $\psi$ the angle between the meridian $\v = \v_{0}$ and
curve $C^{*}$, oriented with respect to the versor
$\overline{\nu}$ we have
\begin{equation}
\mbox{tg} \psi =\(\sin \t \dfrac{d\v}{ds}\)/ \dfrac{d\t}{ds}.
\end{equation}
The previous formulae lead to
\begin{equation}
\left\langle \overline{\nu}, \dfrac{d\overline{\nu}}{ds},
\dfrac{d^{2}\overline{\nu}}{ds^{2}} \right\rangle/
\left\|\dfrac{d\overline{\nu}}{ds}\right\|^{2} = \cos \t
\dfrac{d\v}{ds} + \dfrac{d\psi}{ds}.
\end{equation}
But in our conditions of regularity $\displaystyle\int_{C}\dfrac{d\psi}{ds}ds =
0.$ Thus (2.10.7), (2.10.8) implies the formula (2.10.3)
\end{proof}

It is not difficult to see that the formula (2.10.3) can be
generalized in the case when the curve $C$ of the configuration
$\mathfrak{M}$ has a finite number of angular points. The second
member of the formula (2.10.3) will be additive modified with the
total of variations of angle $\psi$ at the angular points
corresponding to the curve $C$.

Now, one can prove the generalization of Mark Krein formula.

\begin{theorem}
Assume that we have

1. $\mathfrak{M}(C, \overline{\xi}, \pi)$ a Myller configuration
of class $C^{3}$ $($i.e. $C$ is of class $C^{3}$ and $\xi(s)$,
$\overline{\nu}(s)$ are the class $C^{2}$$)$ in which $C$ is a
closed curve, having $s$ as natural parameter.

2. The spherical image $C^{*}$ of $\mathfrak{M}$ determine on the
support sphere $\Sigma$ a simply connected domain of area $\o.$

3. $\s$ the oriented angle between the versors
$\overline{\nu}_{3}(s)$ and $\overline{\xi}(s)$. In these
conditions the following formula hold:
\begin{equation}
\o = 2\pi -\int_{C}G(s)ds +\int_{C}d\s.
\end{equation}
\end{theorem}

\noindent {\bf Proof.} The first two conditions allow to apply the
Lemma 2.10.1. The fundamental equations (1.3.1), (1.3.2) of
$(C,\overline{\nu})$ give us for $\overline{\nu} =
\overline{\nu}_{1}$:
$$
\dfrac{d\overline{\nu}_{1}}{ds} = \chi_{1}\overline{\nu}_{2},\;
\dfrac{d_{2}\overline{\nu}_{1}}{ds^{2}} =
\dfrac{d\chi_{1}}{ds}\overline{\nu}_{2} +
\chi_{1}(-\chi_{1}\overline{\nu}_{1} +
\chi_{2}\overline{\nu}_{3}).
$$
So,
$$
\left\langle \overline{\nu}_{1}, \dfrac{d\overline{\nu}_{1}}{ds},
\dfrac{d^{2}\overline{\nu}_{1}}{ds^{2}} \right\rangle =
\chi_{1}^{2}\chi_{2}.
$$

Thus, the formula (2.10.3), leads to the following formula
$$
\o = 2\pi -\int_{C}\chi_{2}(s)ds.
$$
But, we have $G(s) =\chi_{2}(s)  +\dfrac{d\s}{ds}$, $G(s)$ being
the geodesic curvature of $(C,\overline{\xi})$ in $\mathfrak{M}$.
Then the last formula is exactly (2.10.9).

If $G=0$ for $\mathfrak{M}$, then we have $\o = 2\pi. $

Indeed $G(s) = 0$ along the curve $C$ imply $\o = 2\pi +2k \pi$,
$k\in \mathbb{N}$. But we have $0\leq \o <4\pi,$ so $k=0.$

A particular case of Theorem 2.10.1 is the famous result of Jacobi:

{\it The area $\o$ of the domain $D$ determined on the sphere $\Sigma$
by the closed curve $C^{*}$-spherical image of the principal normals
of a closed curve $C$ in $E_{3}$, assuming $D$ a simply connected
domain, is a half of area of sphere $\Sigma$.}

In this case we consider the Myller configuration $\mathfrak{M}(C,
\overline{\a}, \pi)$, $\pi(s)$ being the rectifying planes of $C$.
\newpage
\thispagestyle{empty}

\def\tg{{\rm{tg}}}

\chapter{Tangent Myller configurations $\mathfrak{M}_{t}$}

The theory of Myller configurations
$\mathfrak{M}(C,\overline{\xi},\pi)$ presented in the Chapter 2
has an important particular case when the tangent versor fields
$\overline{\a}(s)$, $\forall s\in (s_{1},s_{2})$ belong to the
corresponding planes $\pi(s)$. These Myller configurations will be
denoted by $\mathfrak{M}_{t} = \mathfrak{M}_{t}(C,
\overline{\xi},\pi)$ and named {\bf tangent} {\it Myller configuration}.

The geometry of $\mathfrak{M}_{t}$ is much more rich that the
geometrical theory of $\mathfrak{M}$ because in $\mathfrak{M}_{t}$
the tangent field has some special properties. So,
$(C,\overline{\a})$ in $\mathfrak{M}_{t}$ has only three
invariants $\k_{g}, \k_{n}$ and $\tau_{g}$ called {\it geodesic
curvature}, {\it normal curvature} and {\it geodesic torsion}, respectively, of
the curve $C$ in $\mathfrak{M}_{t}$.

The curves $C$ with $\k_{g} = 0$ are {\it geodesic lines} of
$\mathfrak{M}_{t}$; the curves $C$ with the property $\k_{n}=0$
are the {\it asymptotic lines} of $\mathfrak{M}_{t}$ and the curve $C$
for which $\tau_{g} = 0$ are the {\it curvature lines} for
$\mathfrak{M}_{t}$. The mentioned invariants have some geometric
interpretations as the geodesic curvature, normal curvature and
geodesic torsion of a curve $\cal{C}$ on a surfaces $S$ in the
Euclidean space $E_{3}$.

\section{The fundamental equations of $\mathfrak{M}_{t}$}
\setcounter{theorem}{0}\setcounter{equation}{0}

Consider a tangent Myller configuration $\mathfrak{M}_{t} =
(C,\overline{\xi},\pi)$. Thus we have
\begin{equation}
\langle \a(s), \nu(s)\rangle = 0,\; \forall s\in (s_{1},s_{2}).
\end{equation}
The Darboux frame $\cal{R}_{D}$ is $\cal{R}_{D}=(P(s);
\overline{\xi}(s), \overline{\mu}(s),\overline{\nu}(s))$ with
$\overline{\mu}(s) = \overline{\nu}(s)\times \overline{\xi}(s)$.

The fundamental equations of $\mathfrak{M}_{t}$ are obtained by
the fundamental equations (2.1.3), (2.1.4), Chapter 2 of a general
Myller configuration $\mathfrak{M}$ for which the invariant
$c_{3}(s)$ vanishes.

\begin{theorem}
The fundamental equations of the tangent Myller configuration
$\mathfrak{M}_{t}(C, \overline{\xi}, \pi)$ are given by the
following system of differential equations:
\begin{equation}
\dfrac{d\overline{r}}{ds} = c_{1}(s)\overline{\xi}(s) +
c_{2}(s)\overline{\mu}(s),\;\; (c_{1}^{2}+c_{2}^{2} = 1),
\end{equation}
\begin{eqnarray}
\dfrac{d\overline{\xi}}{ds} &=&
G(s)\overline{\mu}(s)+K(s)\overline{\nu}(s)\nonumber\\
\dfrac{d\overline{\mu}}{ds}& =& -G(s)\overline{\xi}(s) +
T(s)\overline{\nu}(s)\\\dfrac{d\overline{\nu}}{ds} &=&
-K(s)\overline{\xi}(s) -T(s)\overline{\mu}(s).\nonumber
\end{eqnarray}
\end{theorem}

The invariants $c_{1},c_{2},G, K,T$ have the same geometric
interpretations and the same denomination as in Chapter 2. So,
$G(s)$ is {\it the geodesic curvature} of the field
$(\overline{C},\overline{\xi})$ in $\mathfrak{M}_{t}$, $K(s)$ is
{\it the normal curvature} and $T(s)$ is {\it the geodesic
torsion} of $(C,\overline{\xi})$ in $\mathfrak{M}_{t}$.

The cases when some invariants $G, K, T$ vanish can be investigate
exactly as in the Chapter 2.

In this respect, denoting $\v = \sphericalangle
(\overline{\xi}_{2}(s),\overline{\nu}(s))$ and using the Frenet
formulae of the versor field $(C,\overline{\xi})$ we obtain the
formulae
\begin{equation}
G = K_{1}\sin \v,\;\; K = K_{1}\cos \v, \;\; T = K_{2}
+\dfrac{d\v}{ds}.
\end{equation}
In \S 5, Chapter 2 we get the relations between the invariants of
$(C,\overline{\xi})$ in $\mathfrak{M}_{t}$ and the invariants of
normal versor field $(C,\overline{\nu})$, (Theorem 2.5.1, Ch 2.) For
$\s = \sphericalangle (\overline{\xi}(s), \overline{\nu}_{3}(s))$
we have
\begin{equation}
K = \chi_{1}\sin \s,\;\; T = \chi_{1}\cos \s,\;\; G = \chi_{2}
+\dfrac{d\s}{ds}.
\end{equation}

Others results concerning $\mathfrak{M}_{t}$ can be deduced from
those of $\mathfrak{M}$.

For instance

\begin{theorem}
(Mark Krein) Assuming that we have:

$1.$ $\mathfrak{M}_{t}(C, \overline{\xi}, \pi)$ a tangent Myller
configuration of class $C^{3}$ in which $C$ is a closed curve,
having $s$ as natural parameter.

$2.$ The spherical image $C^{*}$ of $\mathfrak{M}_{t}$ determines
on the support sphere $\Sigma$ a simply connected domain of area
$\o$.

$3.$ $\s = \sphericalangle (\overline{\xi}, \overline{\nu}_{3})$.
\end{theorem}

{\it In these conditions the following Mark-Krein$'$s formula}
holds:
\begin{equation}
\o = 2\pi - \int_{C}G(s)ds +\int_{C}d\s.
\end{equation}

\section{The invariants of the curve $C$ in $\mathfrak{M}_{t}$}
\setcounter{theorem}{0}\setcounter{equation}{0}

A smooth curve $C$ having $s$ as arclength determines the tangent
versor field $(C,\overline{\a})$ with $\overline{\a}(s) =
\dfrac{d\overline{r}}{ds}$. Consequently we can consider a
particular tangent Myller configuration
$\mathfrak{M}_{t}(C,\overline{\a},\pi)$ defined only by curve $C$
and tangent planes $\pi(s)$.

In this case the geometry of Myller configurations
$\mathfrak{M}_{t}(C,\overline{\a},\pi)$ is called the geometry of
curve $C$ in $\mathfrak{M}_{t}$. The Darboux frame of curve $C$ in
$\mathfrak{M}_{t}$ is $\cal{R}_{D} = (P(s); \overline{\a}(s),
\overline{\mu}^{*}(s), \overline{\nu}(s))$, $\overline{\mu}^{*}(s)
= \overline{\nu}(s)\times \overline{\a}(s)$. It will be called
{\it the Darboux frame} of the curve $C$ in $\mathfrak{M}_{t}$.

\begin{theorem}
The fundamental equations of the curve $C$ in the Myller
configuration $\mathfrak{M}_{t}(C,\overline{\a},\pi)$ are given by
the following system of differential equations:
\begin{equation}
\dfrac{d\overline{r}}{ds} = \overline{\a}(s),
\end{equation}
\begin{eqnarray}
\dfrac{d\overline{\a}}{ds}&=&
\k_{g}(s)\overline{\mu}^{*}(s)+\k_{n}(s)\overline{\nu}(s),\nonumber\\
\dfrac{d\overline{\mu}^{*}}{ds}&=&
-\k_{g}(s)\overline{\a}(s)+\tau_{g}(s)\overline{\nu}(s),\\\dfrac{d\overline{\nu}}{ds}&=&
- \k_{n}(s)\overline{\a}(s) -
\tau_{g}(s)\overline{\mu}^{*}(s).\nonumber
\end{eqnarray}\end{theorem}
Of course $(3.2.1), (3.2.2)$ are the moving equations of the
Darboux frame $\cal{R}_{D}$ of the curve $C$.

The invariants $\k_{g}, \k_{n}$ and $\tau_{g}$ are called: the
{\it geodesic curvature}, {\it normal curvature} and {\it geodesic torsion} of the
curve $C$ in $\mathfrak{M}_{t}.$

Of course, we can prove a fundamental theorem for the geometry of
curves $C$ in $\mathfrak{M}_{t}:$

\begin{theorem}
A priori given $C^{\infty}$-functions $\k_{g}(s),$ $\k_{n}(s)$,
$\tau_{g}(s)$, $s\in [a,b]$, there exists a Myller configuration
$\mathfrak{M}_{t}(C,\overline{\a},\pi)$ for which $s$ is arclength
on the curve $C$ and given functions are its invariants. Two such
configurations differ by a proper Euclidean motion.\end{theorem}

The proof is the same as proof of Theorem 2.1.2, Chapter 2.

\begin{remark}\rm
Let $C$ be a smooth curve immersed in a $C^{\infty}$ surface $S$
in $E_{3}$. Then the tangent plans $\pi$ to $S$ along $C$ uniquely
determines a tangent Myller configuration
$\mathfrak{M}_{t}(C,\overline{\a},\pi)$. Its Darboux frame
$\cal{R}_{D}$ and the invariants $k_{g}, \k_{n},\tau_{g}$ of $C$
in $\mathfrak{M}_{t}$ are just the Darboux frame and geodesic
curvature, normal curvature and geodesic torsion of curve $C$ on
the surface $S$.
\end{remark}

Let $\cal{R}_{F} = (P(s); \overline{\a}_{1}(s),
\overline{\a}_{2}(s), \overline{\a}_{3}(s))$, be the Frenet frame of
curve $C$ with $\overline{\a}_{1}(s) = \overline{\a}(s)$.

The Frenet formulae hold:
\begin{equation}
\dfrac{d\overline{r}}{ds} = \overline{\a}_{1}(s)
\end{equation}
\begin{equation}
\left\{\begin{array}{l}
\dfrac{d\overline{\a}_{1}}{ds} = \k(s)\a_{2}(s);\vspace{1.5mm}\\
\dfrac{d\overline{\a}_{2}}{ds} = -\k(s)\overline{\a}_{1}(s) +
\tau(s)\overline{\a}_{3}(s)\nonumber\vspace{1.5mm}\\
\dfrac{d\overline{\a}_{3}}{ds} = -\tau(s)\overline{\a}_{2}(s)
\end{array}\right.
\end{equation}
where $\k(s)$ is the curvature of $C$ and $\tau(s)$ is the torsion of
$C$.

The relations between the invariants $\k_{g}, \k_{n},\tau_{g}$ and
$\k, \tau$ can be obtained like in Section 4, Chapter 2.

\begin{theorem}
The following formulae hold good
\begin{equation}
\k_{g}(s) = \k \sin \v^{*},\;\; \k_{n}(s) = \k \cos \v^{*},\;\;
\tau_{g}(s) = \tau +\dfrac{d\v^{*}}{ds},
\end{equation}
with $\v^{*} = \sphericalangle(\overline{\a}_{2}(s),
\overline{\nu}(s))$.
\end{theorem}

In the case when we consider the relations between the invariants
$\k_{g}, \k_{n}, \tau_{g}$ and the invariants $\chi_{1}, \chi_{2}$
of the normal versor field $(C,\overline{\nu})$ we have from the
formulae (3.1.5):
\begin{equation}
\k_{n} = \chi_{1}\sin \s, \; \tau_{g} = \chi_{1}\cos \s,\;\;
\k_{g} = \chi_{2} + \dfrac{d\s}{ds},
\end{equation}
where $\s = \sphericalangle (\overline{\a}(s),
\overline{\nu}_{3}(s))$.

It is clear that the second formula (3.2.5) gives us a theorem of
Meusnier type, and for $\v^{*} = 0$ or $\v^{*} = \pm \pi$ we have
$\k_{g} = 0, \k_{n} = \pm \k$, $\tau_{g} = \tau.$

For $\s = 0,$ or $\s = \pm \pi,$ from (3.2.6) we obtain $\k_{n} =
0,$ $\tau_{g} = \pm \chi_{1}$, $\k_{g} = \chi_{2}$.

\section{Geodesic, asymptotic and curvature\, \, \, lines in $\mathfrak{M}_{t}(C,\overline{\a},\pi)$}
\setcounter{theorem}{0}\setcounter{equation}{0}

The notion of parallelism of a versor field $(C,\overline{\a})$ in
$\mathfrak{M}_{t}$ along the curve $C$, investigated in the
Section 8, Chapter 2 for the general case, can be applied now for
the particular case of tangent versor field $(C,\overline{\a})$.
It is defined by the condition $\k_{g}(s) =0,$ $\forall s\in
(s_{1},s_{2})$. The curve $C$ with this property is called
{\it geodesic line} for $\mathfrak{M}_{t}$ (or {\it autoparallel curve}).

The following properties hold:

1. {\it The curve $C$ is a geodesic in the configuration
$\mathfrak{M}_{t}$ iff at every point $P(s)$ of $C$ the osculating
plane of $C$ is normal to $\mathfrak{M}_{t}$.}

2. {\it $C$ is geodesic in $\mathfrak{M}_{t}$ iff the equality
$|\k_{n}| = \k$ holds along $C$.}

3. {\it If $C$ is a straight line in $\mathfrak{M}_{t}$ then $C$ is a
geodesic of $\mathfrak{M}_{t}$.}

\bigskip
{\bf Asymptotics}

The curve $C$ is called {\it asymptotic} in $\mathfrak{M}_{t}$ if
$\k_{n} = 0,$ $\forall s\in (s_{1},s_{2})$. An asymptotic $C$ is
called an {\it asymptotic line}, too. The following properties can be
proved without difficulties:

1. {\it $C$ is asymptotic line in $\mathfrak{M}_{t}$ iff at every point
$P(s)\in C$ the osculating plane of $C$ coincides to the plane
$\pi(s)$.}

2. {\it If $C$ is a straight line then $C$ is asymptotic in
$\mathfrak{M}_{t}$.}

3. {\it $C$ is asymptotic in $\mathfrak{M}_{t}$ iff along $C$,
$|\k_{g}| = \k.$}

4. {\it If $C$ is asymptotic in $\mathfrak{M}_{t}$, then $\tau_{g} =
\tau$ along $C$.}

5. {\it If $\overline{\a}(s)$ is conjugate to $\overline{\a}(s)$ then
$C$ is asymptotic line in $\mathfrak{M}_{t}$.}

Therefore we may say that the asymptotic line in
$\mathfrak{M}_{t}$ are the autoconjugated lines.

\medskip

{\bf Curvature lines}

The curve $C$ is called the {\it curvature line} in the
configurations\linebreak $\mathfrak{M}_{t}(C, \overline{\xi},\pi)$
if the ruled surface $\cal{R}(C,\overline{\nu})$ is a developing
surface.

One knows that $\cal{R}(C,\overline{\nu})$ is a developing surface
iff the following equation holds:
$$
\langle\overline{\a}(s), \overline{\nu}(s),
\dfrac{d\overline{\nu}}{ds}(s)\rangle = 0, \;\; \forall s\in
(s_{1},s_{2}).
$$

Taking into account the fundamental equations of
$\mathfrak{M}_{t}$ one gets:

1. {\it $C$ is a curvature line in $\mathfrak{M}_{t}$ iff its geodesic
torsion $\tau(s) = 0$, $\forall s$.}

2. {\it $C$ is a curvature line in $\mathfrak{M}_{t}$ iff the versors
field $(C,\overline{\mu}^{*})$ are conjugated to the tangent
$\overline{\a}(s)$ in $\mathfrak{M}_{t}$.}

\begin{theorem}
If the curve $C$ of the tangent Myller con\-fi\-gu\-ration
$\mathfrak{M}_{t}(C,$ $\overline{\xi},\pi)$ satisfies two from the
following three conditions

a) $C$ is a plane curve.

b) $C$ is a curvature line in $\mathfrak{M}_{t}$.

c) The angle $\v^{*}(s) = \sphericalangle
(\overline{\a}_{2}(s),\overline{\nu}(s))$ is constant,

\noindent then the curve $C$ verifies the third condition,
too.\end{theorem}

For proof, one can apply the Theorem 3.2.3, Chapter 3.

More general, let us consider two tangent Myller configurations
$\mathfrak{M}_{t}(C,$ $\overline{\a},\pi)$ and
$\mathfrak{M}_{t}^{*}(C, \overline{\a},\pi^{*})$ and $\v^{**}(s) =
\sphericalangle (\overline{\nu}(s), \overline{\nu}^{*}(s))$.

We can prove without difficulties:

\begin{theorem}
Assuming satisfied two from the following three conditions
\begin{itemize}
\item[1.] $C$ is curvature line in $\mathfrak{M}_{t}$.

\item[2.] $C$ is curvature line in $\mathfrak{M}_{t}^{*}$.

\item[3.] The angle $\v^{**}(s)$ is constant for $s\in
(s_{1},s_{2})$.
\noindent Then, the third condition is also verified.
\end{itemize}
\end{theorem}

\section{Mark Krein$^{\prime}$s formula}
\setcounter{theorem}{0}\setcounter{equation}{0}

In the Section 10, Chapter 2, we have defined the spherical image
of a general Myller configuration $\mathfrak{M}$. Of course the
definition applies for tangent configurations $\mathfrak{M}_{t}$,
too.

But now will appear some new special properties.

If in the formulae from Section 10, Chapter 2 we take
$\overline{\xi}(s) = \overline{\a}(s)$, $\forall s$, the
invariants $G,K,T$ reduce to geodesic curvature, normal curvature
and geodesic torsion, respectively, of the curve $C$ in
$\mathfrak{M}_{t}$.

Such that we obtain the formulae:
\begin{equation}
\k_{n} = -\dfrac{ds^{*}}{ds}\cos \t_{1},\;\; \t_{1} =
\sphericalangle (\overline{\nu}_{2},\overline{\a})
\end{equation}
\begin{equation}
\tau_{g} = -\dfrac{ds^{*}}{ds}\cos \t_{2},\;\; \t_{2} =
\sphericalangle (\overline{\nu}_{2},\overline{\mu}^{*})
\end{equation}
which have as consequences:

1. {\it $C$ is asymptotic in $\mathfrak{M}_{t}$ iff $\t_{1} = \pm
\dfrac{\pi}{2}$.}

2. {\it $C$ is curvature line in $\mathfrak{M}_{t}$ iff $\t_{2} = \pm
\dfrac{\pi}{2}$.}

\smallskip

The following Mark Krein$'$s theorem holds:

\begin{theorem}
Assume that we have

\begin{itemize}
\item[$1.$] A Myller configuration $\mathfrak{M}_{t}(C,
\overline{\xi},\pi)$ of class $C^{k}$, $(k\geq 3)$ where $s$ is
the natural parameter on the curve $C$ and $C$ is a closed curve.

\item[$2.$] The spherical image $C^{*}$ of $\mathfrak{M}_{t}$ determine on
the support sphere $\Sigma$ a simply connected domain of area
$\o$.

\item[$3.$] $\s = \sphericalangle (\overline{\nu}_{3}, \overline{\a})$.

In these conditions Mark Krein$'$s formula holds:
\begin{equation}
\o  =2\pi - \int_{C}\k_{g}(s)ds +\int_{C}d\s.
\end{equation}
\end{itemize}
\end{theorem}

Indeed the formula (2.10.9) from Theorem 2.10.1, Chapter 2 is
equivalent to (3.4.3).

The remarks from the end of Section 10, Chapter 2 are valid, too.

\section{Relations between the invariants $G,K,T$ of the versor field
$(C,\overline{\xi})$ in $\mathfrak{M}_{t}(C,\overline{\xi},\pi)$ and the
invariants of tangent versor field $(C,\overline{\a})$ in
$\mathfrak{M}_{t}$}
\setcounter{theorem}{0}\setcounter{equation}{0}

Let $\mathfrak{M}_{t}(C,\overline{\xi},\pi)$ be a tangent Myller
configuration, $\cal{R}_{D} = (P(s), \overline{\xi}(s),$$
\overline{\mu}(s),$ $\overline{\nu}(s))$ its Darboux frame and the
tangent Myller configuration
$\mathfrak{M}_{t}^{\prime}(C,\overline{a},\pi)$ determined by the
plane field $(C,\pi(s))$, with $\overline{\a}
=\dfrac{d\overline{r}}{ds}$-tangent versors field to the oriented
curve $C$. The Darboux frame $\cal{R}_{D}^{\prime} = (P(s);
\overline{\a}(s),$ $ \overline{\mu}^{*}(s),$ $\overline{\nu})$ and
oriented angle $\l = \sphericalangle
(\overline{\a},\overline{\xi})$ allow to determine $\cal{R}_{D}$
by means of formulas:
\begin{eqnarray}
\overline{\xi}&=& \overline{\a}\cos \l + \overline{\mu}^{*}\sin
\l\nonumber\\\overline{\mu}&=& -\overline{\a}\sin \l
+\overline{\mu}^{*}\cos \l\\\overline{\nu}&=&
\overline{\nu}.\nonumber
\end{eqnarray}
The moving equations (3.1.2), (3.1.3) and (3.2.1), (3.2.3) of
$\cal{R}_{D}$ and $\mathcal{R}^{\prime}_{D}$ lead to the following
relations between invariants $\k_{g}, \k_{n}$ and $\tau_{g}$ of
the curve $C$ in $\mathfrak{M}_{t}^{\prime}$ and the invariants
$G,K,T$ of versor field $(C,\overline{\xi})$ in
$\mathfrak{M}_{t}:$
\begin{eqnarray}
K&=& \k_{n}\cos \l + \tau_{g}\sin \l\nonumber\\T&=& -\k_{n}\sin \l
+ \tau_{g}\cos \l\\G&=& \k_{g} + \dfrac{d\l}{ds}.\nonumber
\end{eqnarray}
The two first formulae (3.5.2) imply
\begin{equation}
K^{2} + T^{2} = \k_{n}^{2} + \tau_{g}^{2} =
\left(\dfrac{ds^{*}}{ds}\right)^{2}.
\end{equation}

The formula (3.5.2) has some important consequences:

\begin{theorem}
1. The curve $C$ has $\k_{g}(s)$ as geodesic curvature on the
developing surface $E$ generated by planes $\pi(s)$, $s\in
(s_{1},s_{2})$.

2. $G$ is an intrinsec invariant of the developing surface $E$.
\end{theorem}

\begin{proof} 1. Since the planes $\pi(s)$ pass through tangent line
$(P(s), \overline{\a}(s))$ to curve $C$, the developing surface
$E$, enveloping by planes $\pi(s)$ passes through curve $C$. So,
$\k_{g}$ is geodesic curvature of $C\subset E$ at point $P(s)$.

2. The invariants $\k_{g}(s)$ and $\l(s)$ are the intrinsic
invariants of surface $E$. It follows that the invariant $G$ has
the same property.\end{proof}

Now, we investigate an extension of Bonnet result.

\begin{theorem}
Supposing that the versor field $(C, \overline{\xi})$ in tangent
Myller configuration $\mathfrak{M}_{t}(C,\overline{\xi},\pi)$ has
two from the following three properties:
\begin{itemize}
\item[1.] $\overline{\xi}(s)$ is a parallel versor field in
$\mathfrak{M}_{t}$,

\item[2.] The angle $\l = \sphericalangle (\overline{\a},\overline{\xi})$ is
constant,

\item[3.] The curve $C$ is geodesic in $\mathfrak{M}_{t},$
\end{itemize}
then it has the third property, too.
\end{theorem}

The proof is based on the last formula (3.5.2). Also, it is not
difficult to prove:

\begin{theorem}
The normal curvature and geodesic torsion of the versor field
$(C,\overline{\mu}^{*})$ in tangent Myller configuration
$\mathfrak{M}_{t}(C,\overline{\xi},\pi)$, respectively are the
geodesic torsion and normal curvature with opposite sign of the
curve $C$ in $\mathfrak{M}_{t}$.
\end{theorem}

The same formulae (3.5.2) for $K(s)=0,$ $s\in (s_{1}, s_{2})$ and
$\tau_{g}(s)\neq 0$ imply:
\begin{equation}
\tg \l = -\dfrac{\k_{n}(s)}{\tau_{g}(s)}.
\end{equation}

In the theory of surfaces immersed in $E_{3}$ this formula was
independently established by E. Bortolotti [3] and Al. Myller
[31-34].

\begin{definition}
A curve $C$ is called a geodesic helix in tangent configuration
$\mathfrak{M}_{t}(C,\overline{\xi},\pi)$ if the angle $\l(s) =
\sphericalangle (\overline{\a},\overline{\mu})$ is constant $(\neq
0,\pi)$.
\end{definition}

The following two results can be proved without difficulties.

\begin{theorem}
If $C$ is a geodesic helix in
$\mathfrak{M}_{t}(C,\overline{\xi},\pi)$, then $$\k_{n}/\tau_{g} =
\pm \dfrac{\k}{\tau}.$$
\end{theorem}

Another consequence of the previous theory is given by:

\begin{theorem}
If $C$ is geodesic in $\mathfrak{M}_{t}(C,\overline{\xi},\pi)$, it
is a geodesic helix for $\mathfrak{M}_{t}$ iff $C$ is a cylinder
helix in the space $E_{3}$.
\end{theorem}

We can make analogous consideration, taking $\tau=0$ in the formula
(3.5.2).

We stop here the theory of Myller configurations
$\mathfrak{M}(C,\overline{\xi},\pi)$ in Euclidean space $E_{3}$.
Of course, it can be extended for Myller configurations
$\mathfrak{M}(C,\overline{\xi}, \pi)$ with $\pi$ a $k$-plane in
the space $E^{n}, n\geq 3$, $k<n$ or in more general spaces, as
Riemann spaces, [26], [27].

\newpage
\thispagestyle{empty}

\def\pr{{\rm{pr}}}

\chapter[Applications of the theory of Myller configuration $\mathfrak{M}_{t}$]{Applications of theory of Myller configuration $\mathfrak{M}_{t}$ in the geometry of surfaces in
$E_{3}$}

A first application of the theory of Myller configurations
$\mathfrak{M}(C,\overline{\xi},\pi)$ in the Euclidean space $E_{3}$
can be realized to the geometry of surfaces $S$ embedded in
$E_{3}$. We obtain a more clear study of curves $C\subset S$, a
natural definition of Levi-Civita parallelism of vector field
$\overline{V}$, tangent to $S$ along $C$, as well as the notion of
concurrence in Myller sense of vector field, tangent to $S$ along
$C$. Some new concepts, as those of mean torsion of $S$, the total
torsion of $S$, the Bonnet indicatrix of geodesic torsion and its
relation with the Dupin indicatrix of normal curvatures are
introduced. A new property of Bacaloglu-Sophie-Germain curvature
is proved, too. Namely, it is expressed in terms of the total
torsion of the surface $S$.

\section{The fundamental forms of surfaces in $E_{3}$}
\setcounter{theorem}{0}\setcounter{equation}{0}

Let $S$ be a smooth surface embedded in the Euclidean space
$E_{3}$. Since we use the classical theory of surfaces in $E_{3}$
[see, for instance: Mike Spivak, Diff. Geom. vol. I, Bearkley,
1979], consider the analytical representation of $S$, of class $C^k$, $(k\geq 3,\ {\rm{or\ }}k=\infty)$:
\begin{equation}
\overline{r} = \overline{r}(u,v),\;\; (u,v)\in D.
\end{equation}
$D$ being a simply connected domain in plane of variables $(u,v)$,
and the following condition being verified:
\begin{equation}
\overline{r}_{u}\times \overline{r}_{v}\neq \overline{0}\;\;\text{
for }\forall (u,v)\in D.
\end{equation}

Of course we adopt the vectorial notations
\begin{eqnarray*}
\overline{r}_{u} = \dfrac{\p \overline{r}}{\p u},\;\;
\overline{r}_{v} = \dfrac{\p\overline{r}}{\p v}.
\end{eqnarray*}
Denote by $C$ a smooth curve on $S$, given by the parametric
equations
\begin{equation}
u = u(t),\;\; v = v(t),\;\; t\in (t_{1},t_{2}).
\end{equation}
In this chapter all geometric objects or mappings are considered
of $C^{\infty}$-class, up to contrary hypothesis.

In space $E_{3}$, the curve $C$ is represented by
\begin{equation}
\overline{r} = \overline{r}(u(t), v(t)),\;\; t\in (t_{1},t_{2}).
\end{equation}
Thus, the following vector field
\begin{equation}
d\overline{r} = \overline{r}_{u}du + \overline{r}_{v}dv
\end{equation}
is tangent to $C$ at the points $P(t) = P(\overline{r}(u(t)),
v(t))\in C.$ The vectors $\overline{r}_{u}$ and $\overline{r}_{v}$
are tangent to the parametric lines and $d\overline{r}$ is tangent
vector to $S$ at point $P(t).$

From (4.1.2) it follows that the scalar function
\begin{equation}
\D = \|\overline{r}_{u}\times \overline{r}_{v}\|^{2}
\end{equation}
is different from zero on $D$.

The unit vector $\overline{\nu}$
\begin{equation}
\overline{\nu} = \dfrac{\overline{r}_{u}\times
\overline{r}_{v}}{\sqrt{\D}}
\end{equation}
is normal to surface $S$ at every point $P(t).$

The tangent plane $\pi$ at $P$ to $S$, has the equation
\begin{equation}
\langle \overline{R}-\overline{r}(u,v), \overline{\nu}\rangle=0.
\end{equation}
Assuming that $S$ is orientable, it follows that $\overline{\nu}$
is uniquely determined and the tangent plane $\pi$ is oriented (by
means of versor $\overline{\nu}$), too.

The {\it first fundamental form} $\phi$ of surface $S$ is defined by
\begin{equation}
\phi = \langle d\overline{r}, d\overline{r}\rangle,\;\; \forall
(u,v)\in D
\end{equation}
or by $\phi(du,dv) = \langle d\overline{r}(u,v),
d\overline{r}(u,v)\rangle$, $\forall (u,v)\in D.$

Taking into account the equality (4.1.5) it follows that
$\phi(du,dv)$ is a quadratic form:
\begin{equation}
\phi(du,dv) = E du^{2} + 2F du dv + G dv^{2}.
\end{equation}
The coefficients of $\phi$ are the functions, defined on $D$:
\begin{equation}
\hspace*{1cm}E(u,v) = \langle \overline{r}_{u}, \overline{r}_{u}\rangle, F(u,v)
= \langle \overline{r}_{u}, \overline{r}_{v}\rangle,\; G(u,v) =
\langle \overline{r}_{v}, \overline{r}_{v} \rangle.
\end{equation}
But we have the discriminant $\D$ of $\phi:$
\begin{equation}
\D = EG - F^{2}.
\end{equation}
\begin{equation}
E>0,\; G>0,\; \D >0.
\end{equation}
Consequence: the first fundamental form $\phi$ of $S$ is
positively defined.

Since the vector $d\overline{r}$ does not depend on
parametrization of $S$, it follows that $\phi$ has a geometrical meaning
with respect to a change of coordinates in $E_{3}$ and with
respect to a  change of parameters $(u,v)$ on $S$.

Thus $ds,$ given by:
\begin{equation}
ds^{2} = \phi(du,dv) = E du^{2} + 2F du dv + G dv^{2}
\end{equation}
is called the {\it element of arclength} of the surface $S$.

The arclength of a curve $C$ an $S$ is expressed by
\begin{equation}
s = \int_{t_{0}}^{t}\left\{E\(\dfrac{du}{d\s}\)^{2} + 2F
\dfrac{du}{d\s}\dfrac{dv}{d\s} + G\(\dfrac{dv}{d\s}\)^{2}\right\}^{1/2}d\s.
\end{equation}

The function $s = s(t),$ $t\in [a,b]\subset (t_{1},t_{2})$,
$t_{0}\in [a,b]$ is a diffeomorphism from the interval $[a,b]\to
[0,s(t)]$. The inverse function $t =t(s)$, determines a new
parametrization of curve $C$: $\overline{r} = \overline{r}(u(s),
v(s))$. The tangent vector $\dfrac{d\overline{r}}{ds}$ is a versor:
\begin{equation}
\overline{\a} = \dfrac{d\overline{r}}{ds} =
\overline{r}_{u}\dfrac{du}{ds} + \overline{r}_{v}\dfrac{dv}{ds}.
\end{equation}

For two tangent versors $\overline{\a} =
\dfrac{d\overline{r}}{ds}$, $\overline{\a}_{1} =
\dfrac{\d\overline{r}}{\d s}$ at point $P\in S$ the angle
$\sphericalangle(\overline{\a}, \overline{\a}_{1})$ is expressed
by
\begin{equation}
\hspace*{13mm}\cos\hspace{-1mm} \sphericalangle(\hspace{-0.5mm}\overline{\a}, \overline{\a}_{1}\hspace{-0.5mm})\hspace{-1mm} =\hspace{-1mm}
 \dfrac{E du
\d u\hspace{-0.8mm} +\hspace{-0.8mm} F(du \d v\hspace{-0.8mm} +\hspace{-0.8mm} \d u dv)\hspace{-0.8mm} +\hspace{-0.8mm} Gdv \d v}{\sqrt{Edu^{2}\hspace{-0.8mm} +\hspace{-0.8mm} 2F du dv\hspace{-0.8mm} +\hspace{-0.8mm}
Gdv^{2}}\hspace{-1mm}\cdot\hspace{-1mm} \sqrt{E\d u^{2}\hspace{-0.8mm} +\hspace{-0.8mm} 2F \d u \d v\hspace{-0.8mm} +\hspace{-0.8mm} G\d v^{2}}}.
\end{equation}
{\it The second fundamental form} $\psi$ of $S$ at point $P\in S$ is
defined by
\begin{equation}
\psi = \langle \overline{\nu}, d^{2}\overline{r}\rangle = -
\langle d\overline{r}, d\overline{\nu}\rangle\;\; \forall (u,v)\in
D.
\end{equation}

Also, we adopt the notation $\psi(du,dv) = - \langle
d\overline{r}(u,v),d\nu(u,v)\rangle.$ $\psi$ is a quadratic form:
\begin{equation}
\psi(du,dv) =L du^{2} + 2M du dv  + N dv^{2},
\end{equation}
having the coefficients functions of $(u,v)\in D:$
\begin{equation}
\begin{array}{l}
L(u,v) = \dfrac{1}{\sqrt{\D}} \langle \overline{r}_{u},
\overline{r}_{u}, \overline{r}_{uu}\rangle,\; M(u,v) =
\dfrac{1}{\sqrt{\D}}\langle \overline{r}_{u}, \overline{r}_{v},
\overline{r}_{uv} \rangle\\N(u,v) =
\dfrac{1}{\sqrt{\D}}<\overline{r}_{u}, \overline{r}_{v},
\overline{r}_{vv}>.
\end{array}
\end{equation}
Of course from $\psi =-\langle d\overline{r},
d\overline{\nu}\rangle$ it follows that $\psi$ has geometric
meaning.

The two fundamental forms $\phi$ and $\psi$ are enough to
determine a surface $S$ in Euclidean space $E_{3}$ under proper
Euclidean motions. This property results by the integration of the
Gauss-Weingarten formulae and by their differential consequences
given by Gauss-Codazzi equations.

\section{Gauss and Weingarten formulae}
\setcounter{theorem}{0}\setcounter{equation}{0}

Consider the moving frame$$ \cal{R} = (P(\overline{r});
\overline{r}_{u}, \overline{r}_{v}, \overline{\nu})
$$
on the smooth surface $S$ (i.e. $S$ is of class $C^\infty$).

The moving equations of $R$ are given by the Gauss and Weingarten
formulae.

The Gauss formulae are as follows:
\begin{equation}
\begin{array}{lll}
\overline{r}_{uu}&=& \left\{\begin{array}{l}1\\11\end{array}\right\} \overline{r}_{u} + \left\{\begin{array}{l}2\\ 11\end{array}\right\}
 \overline{r}_{v} +
L\overline{\nu},\\\overline{r}_{uv}&=& \left\{\begin{array}{l}1\\
12\end{array}\right\}\overline{r}_{u} + \left\{\begin{array}{l}2\\ 12\end{array}\right\} \overline{r}_{v} +
M\overline{\nu}\\\overline{r}_{vv}&=&\left\{\begin{array}{l}1\\
22\end{array}\right\}\overline{r}_{u} +\left\{\begin{array}{l}2\\ 22\end{array}\right\} \overline{r}_{v} +
N\overline{\nu}.
\end{array}
\end{equation}

The coefficients $\left\{\begin{array}{cc} i\\jk
\end{array}\right\}$, $(i,j,k=1,2; u=u^{1}, v= u^{2})$ are called the
Christoffel symbols of the first fundamental form $\phi:$
\begin{equation}
\hspace*{5mm}{1\brace 11} =-\dfrac{1}{\sqrt{D}} \langle \overline{\nu},
\overline{r}_{v}, \overline{r}_{uu}\rangle, { 1\brace 21} ={
1\brace 21} = -\dfrac{1}{\sqrt{\D}}\langle \overline{\nu},
\overline{r}_{v}, \overline{r}_{uv}\rangle,
\end{equation}
$$
{ 1\brace 22} = -\dfrac{1}{\sqrt{\D}}\langle \overline{\nu},
\overline{r}_{v}, \overline{r}_{vv} \rangle,
$$
\begin{eqnarray}
\hspace*{14mm}{ 2\brace 11}
 = - \dfrac{1}{\sqrt{\D}}\langle \overline{\nu}, \overline{r}_{u},
\overline{r}_{uu}\rangle, {2\brace 21}
 = {2\brace 12}
 = -\dfrac{1}{\sqrt{\D}}\langle \overline{\nu},\overline{r}_{u},
\overline{r}_{uv}\rangle,\nonumber\\ {2\brace 22} =
-\dfrac{1}{\sqrt{\D}}\langle
\overline{\nu},\overline{r}_{u},\overline{r}_{vv}\rangle.\nonumber
\end{eqnarray}

The Weingarten formulae are given by:
\begin{eqnarray}
\dfrac{\p \overline{\nu}}{\p u}&=& \dfrac{1}{\sqrt{\D}}\{(FM -
GL)\overline{r}_{u} + (FL - EM)\overline{r}_{v}\}\\\dfrac{\p
\overline{\nu}}{\p v}&=& \dfrac{1}{\sqrt{\D}}\{(FN -
GM)\overline{r}_{u} + (FM -EN)\overline{r}_{v}\}.\nonumber
\end{eqnarray}
Of course, the equations (4.2.1) and (4.2.3) express the variation
of the moving frame $R$ on surface $S$.

Using the relations $\overline{r}_{uuv} = \overline{r}_{uvu},$
$\overline{r}_{vuv} = \overline{r}_{vvu}$ and
$\dfrac{\p^{2}\overline{\nu}}{\p uv} =
\dfrac{\p^{2}\overline{\nu}}{\p v u}$ applied to (4.2.1), (4.2.3)
one deduces the so called fundamental equations of surface $S$,
known as the Gauss-Codazzi equations.

A fundamental theorem can be proved, when the first and second
fundamental from $\phi$ and $\psi$ are given and Gauss-Codazzi
equations are verified (Spivak [75]).

\newpage

\section{The tangent Myller configuration $\mathfrak{M}_{t}(C,$ $\overline{\xi},\pi)$ associated to a tangent versor field $(C,\overline{\xi})$ on a surface
$S$}
\setcounter{theorem}{0}\setcounter{equation}{0}

Assuming that $C$ is a curve on surface $S$ in $E_{3}$ and $(C,
\overline{\xi})$ is a versor field tangent to $S$ along $C$, there
is an uniquely determined tangent Myller configuration
$\mathfrak{M}_{t} = \mathfrak{M}_{t}(C, \overline{\xi},\pi)$ for
which $(C,\pi)$ is tangent field planes to $S$ along $C$.

The invariants $(c_{1}, c_{2}, G,K,T)$ of $(C, \overline{\xi})$ in
$\mathfrak{M}_{t}(C, \overline{\xi},\pi)$ will be called the
invariants of tangent versor field $(C, \overline{\xi})$ on
surface $S$.

Evidently, these invariants have the same values on every smooth
surface $S^{\prime}$ which contains the curve $C$ and is tangent
to $S$.

Using the theory of $\mathfrak{M}_{t}$ from Chapter 3, for $(C,
\overline{\xi})$ tangent to $S$ along $C$ we determine:

1. The Darboux frame
\begin{equation}
R_{D} = (P(s); \overline{\xi}(s), \overline{\mu}(s),
\overline{\nu}(s)),
\end{equation}
where $s$ is natural parameter on $C$ and
\begin{equation}
\overline{\nu} = \dfrac{1}{\sqrt{\D}} \overline{r}_{u}\times
\overline{r}_{v},\;\; \overline{\mu} = \overline{\nu}\times
\overline{\xi}.
\end{equation}
Of course $\cal{R}_{D}$ is orthonormal and positively oriented.

The invariants $(c_{1}, c_{2}, G,K,M)$ of $(C, \overline{\xi})$ on
$S$ are given by (3.1.2), (3.1.3) Chapter 3. So we have  the
moving equations of $R.$

\begin{theorem}
The moving equations of the Darboux frame of tangent versor field
$(C,\overline{\xi})$ to $S$ are given by
\begin{equation}
\dfrac{d\overline{r}}{ds} = \overline{\a}(s) =
c_{1}(s)\overline{\xi} + c_{2}(s)\overline{\mu},\;\; (c_{1}^{2} +
c_{2}^{2} =1)
\end{equation}
and
\newpage
\begin{eqnarray}
\dfrac{d\overline{\xi}}{ds} &=& G(s)\overline{\mu}(s) +
K(s)\overline{\nu}(s),\nonumber\\ \dfrac{d\overline{\mu}}{ds} &=&
-G(s)\overline{\xi}(s) +
T(s)\overline{\nu}(s),\\\dfrac{d\overline{\nu}}{ds} &=&
-K(s)\overline{\xi} - T(s)\overline{\mu}(s),\nonumber
\end{eqnarray}
where $c_{1}(s), c_{2}(s), G(s), K(s), T(s)$ are invariants with
respect to the changes of coordinates on $E_{3}$, with respect of transformation of local coordinates on $S$, $(u,v)\to(\widetilde{u},\widetilde{v})$ and with respect
to transformations of natural parameter $s\to s_{0}+s$.
\end{theorem}

A fundamental theorem can be proved exactly as in Chapter 2.

\begin{theorem}
Let be $c_{1}(s),c_{2}(s)$, $[c_{1}^{2} + c_{2}^{2}=1]$, $G(s),
K(s), T(s)$, $s\in [a,b]$, a priori given smooth functions. Then
there exists a tangent Myller configuration
$\mathfrak{M}_{t}(C,\overline{\xi},\pi)$ for which $s$ is the
arclength of curve $C$ and the given functions are its invariants.
$\mathfrak{M}_{t}$ is determined up to a proper Euclidean motion.
The given functions are invariants of the versor field
$(C,\overline{\xi})$ for any smooth surface $S$, which contains
the curve $C$ and is tangent planes $\pi(s)$, $s\in [a,b]$.
\end{theorem}

The geometric interpretations of the invariants $G(s), K(s)$ and
$T(s)$ are those mentioned in the Section 2, Chapter 3.

So, we get
$$
G(s) = \lim_{\D s\to 0}\dfrac{\D\psi_{1}}{\D s}, \;\; \text{
with } \D \psi_{1} = \sphericalangle (\overline{\xi}(s),
\pr_{\pi(s)}\overline{\xi}(s+\D s)),
$$
$$
K(s) = \lim_{\D s\to 0}\dfrac{\D\psi_{2}}{\D s}, \;\; \text{
with } \D \psi_{2} = \sphericalangle (\overline{\mu},
\pr_{(P; \overline{\xi},\overline{\nu})}\overline{\xi}(s+\D
s)),
$$
$$
T(s) = \lim_{\D s\to 0}\dfrac{\D\psi_{3}}{\D s}, \;\; \text{
with } \D \psi_{3} = \sphericalangle (\overline{\mu},
\pr_{(P; \overline{\mu},\overline{\nu})}\overline{\xi}(s+\D
s)).
$$

For this reason $G(s)$ is called the {\it geodesic curvature} of $(C,
\overline{\xi})$ on $S$; $K(s)$ is called {\it the normal curvature} of
$(C, \overline{\xi})$ on $S$ and $T(s)$ is named the {\it geodesic
torsion} of the tangent versor field $(C, \overline{\xi})$ on $S$.

The calculus of invariants $G,K,T$ is exactly the same as it has been done in Chapter
2, (2.3.3), (2.3.4):
\begin{eqnarray}
&&G(s) = \left\langle \overline{\xi}, \dfrac{d\overline{\xi}}{ds},
\overline{\nu} \right\rangle,\;\; K(s) =
\left\langle\dfrac{d\overline{\xi}}{ds},\overline{\nu}\right\rangle = -\left\langle
\overline{\xi},\dfrac{d\overline{\nu}}{ds}
\right\rangle,\nonumber\\&&T(s) = \left\langle \overline{\xi},
\overline{\nu}, \dfrac{d\overline{\nu}}{ds}\right\rangle.
\end{eqnarray}
The analytical expressions of invariants $G,K,T$ on $S$ are given in next
section.

\section{The calculus of invariants $G,K,T$ on a surface}
\setcounter{theorem}{0}\setcounter{equation}{0}

The tangent versor $\overline{\a} = \dfrac{d\overline{r}}{ds}$ to
the curve $C$ in $S$ is
\begin{equation}
\overline{\a}(s) = \dfrac{d\overline{r}}{ds} = \overline{r}_{u}
\dfrac{du}{ds} + \overline{r}_{v}\dfrac{dv}{ds}.
\end{equation}
$\overline{\a}(s)$ is the versor of the tangent vector
$d\overline{r}$, from (4,4)
\begin{eqnarray}
d\overline{r} = \overline{r}_{u}du + \overline{r}_{v}dv.
\end{eqnarray}

The coordinate of a vector field $(C, \overline{V})$ tangent to
surface $S$ along the curve $C\subset S$, with respect to the
moving frame $\cal{R} = (P; \overline{r}_{u}, \overline{r}_{v},
\overline{\nu})$ are the $C^{\infty}$ functions $V^{1}(s),
V^{2}(s)$:
$$
\overline{V}(s)=V^{1}(s)\overline{r}_{u}+V^{2}(s)\overline{r}_{v}.
$$
The square of length of vector $\overline{V}(s)$ is:
$$
\langle \overline{V}(s),\overline{V}(s)\rangle  = E(V^{1})^{2} +
2FV^{1}V^{2} + G(V^{2})^{2}
$$
and the scalar product of two tangent vectors $\overline{V}(s)$
and
$\overline{U}(s)=U^{1}(s)\overline{r}_{u}+U^{2}(s)\overline{r}_{v}$
is as follows
$$
\langle \overline{U}, \overline{V}\rangle = EU^{1}V^{1} +
F(U^{1}V^{2} + V^{1}U^{2}) +GU^{2}V^{2}.
$$
Let $C$ and $C^{\prime}$ two smooth curves on $S$ having $P(s)$ as
common point. Thus the tangent vectors $d\overline{r},
\delta\overline{r}$ at a point $P$ to $C$, respectively to $C^{\prime}$
are
$$
d\overline{r} = \overline{r}_{u}du + \overline{r}_{v}dv,
$$
$$
\d \overline{r} = \overline{r}_{u} \d u +\overline{r}_{v}\d v.
$$
They correspond to tangent directions $(du,dv)$, $(\d u, \d v)$ on
surface $S$.

Evidently, $\overline{\a}(s) = \dfrac{d\overline{r}}{ds}$ and
$\overline{\xi}(s) = \dfrac{\d \overline{r}}{\d s}$ are the tangent
versors to curves $C$ and $C^{\prime}$, respectively at point
$P(s)$, and $(C, \overline{\a})$, $(C, \overline{\xi})$ are the
tangent versor fields along the curve $C$.

The Darboux frame $\cal{R}_{D} = (P(s); \overline{\xi}(s),
\overline{\mu}(s), \overline{\nu}(s))$ has the versors
$\overline{\xi}(s), \overline{\mu}(s),$ $\overline{\nu}(s)$ given
by
\begin{eqnarray}
\overline{\xi}(s) &=& \overline{r}_{u}\dfrac{\d u}{\d s} +
\overline{r}_{v}\dfrac{\d v}{\d s},\nonumber\\\overline{\mu}(s)&=&
\dfrac{1}{\sqrt{\D}}\left[(E \overline{r}_{v} - F
\overline{r}_{u})\dfrac{\d u}{\d s} +(F \overline{r}_{v} - G
\overline{r}_{u})\dfrac{\d v}{\d s}\right],\\\overline{\nu}(s) &=&
\dfrac{1}{\D}(\overline{r}_{u} \times \overline{r}_{v}).\nonumber
\end{eqnarray}
Taking into account (4.3.3) and (4.4.3) it follows:
\begin{eqnarray}
c_{1}&=& \langle \overline{\a}, \overline{\xi}\rangle = \dfrac{E du
\d u +F (du \d v + dv \d u)+G dv \d v}{\sqrt{\phi(du,
dv)}\sqrt{\phi(\d u, \d v)}}\\c_{2}&=& \langle
\overline{\a},\overline{\mu}\rangle=\dfrac{\sqrt{\D}(\d u dv - \d v
du)}{\sqrt{\phi(du,dv)}\sqrt{\phi(\d u, \d v)}}\nonumber
\end{eqnarray}
where
$$
\phi(du, dv) = \langle d\overline{r}, d\overline{r}\rangle,\;\;
\phi(\d u, \d v) = \langle \d \overline{r}, \d
\overline{r}\rangle.
$$
In order to calculate the invariants $G,K,T$ we need to determine
the vectors $\dfrac{d\overline{\xi}}{ds},
\dfrac{d\overline{\nu}}{ds}$.

By means of (4.4.3) we obtain
\begin{eqnarray}
\dfrac{d\overline{\xi}}{ds}&=& \dfrac{d}{ds}\left(\dfrac{\d
\overline{r} }{\d s}\right) =
\overline{r}_{uu}\dfrac{du}{ds}\dfrac{\d u}{\d s} +
\overline{r}_{uv}\left(\dfrac{du}{ds}\dfrac{\d v}{\d s} + \dfrac{\d u
}{\d s}\dfrac{dv}{ds}\right)
\nonumber\\&&+\overline{r}_{vv}\dfrac{dv}{ds} \dfrac{\d v}{\d s} +
\overline{r}_{u}\dfrac{d}{ds}\left(\dfrac{\d u }{\d s}\right)
+\overline{r}_{v}\dfrac{d}{ds}\left(\dfrac{\delta\nu}{\delta s}\right).
\end{eqnarray}
The scalar mixt $\langle\overline{r}_{u}, \overline{r}_{v},
\overline{\nu})\rangle$ is equal to $\sqrt{\D}$. It results:
\begin{eqnarray}
\left\langle \overline{\xi}, \dfrac{d\overline{\xi}}{ds},
\overline{\nu}\right\rangle& =&
\sqrt{\D}\left[\dfrac{d}{ds}\left(\dfrac{\d v }{\d s}\right)
\dfrac{du}{ds} - \dfrac{d}{ds}\left(\dfrac{\d u }{\d s
}\right)\dfrac{dv}{ds}\right] \nonumber\\&&+ \langle
\overline{r}_{u},
\overline{r}_{uu},\overline{\nu}\rangle\dfrac{du}{ds}\left(\dfrac{\d
u}{\d s}\right)^{2}\nonumber+\\&&+ \langle \overline{r}_{u},
\overline{r}_{uv}, \overline{\nu}
\rangle\left(\dfrac{du}{ds}\dfrac{\d v}{\d s} +
\dfrac{dv}{ds}\dfrac{\d u}{\d s} \right)\dfrac{\d u}{\d s} + \\&&+
\langle \overline{r}_{u}, \overline{r}_{vv},
\overline{\nu}\rangle\dfrac{\d u}{\d s}\dfrac{dv}{ds}\dfrac{\d v}{\d
s }+\nonumber\\&&+\langle \overline{r}_{v}, \overline{r}_{uu},
\overline{\nu}\rangle \dfrac{du}{ds}\dfrac{\d u}{\d s}\dfrac{\d v}{\d
s } + \langle \overline{r}_{v}, \overline{r}_{uv},
\overline{\nu}\rangle\left(\dfrac{du}{ds}\dfrac{\d v}{\d s} +
\dfrac{dv}{ds} \dfrac{\d u}{\d s}\right)\dfrac{\d v}{\d
s}\nonumber\\&&+\langle \overline{r}_{v}, \overline{r}_{vv},
\overline{\nu}\rangle\dfrac{dv}{ds}\left(\dfrac{\d v}{\d s
}\right)^{2}.\nonumber
\end{eqnarray}
Taking into account (4.2.1), (4.2.2), (4.3.5) and (4.4.5) we have

\begin{proposition}
The geodesic curvature of the tangent versor field
$(C,\overline{\xi})$ on $S$, is expressed as follows:
\begin{eqnarray}
G(\d,d)&=& \sqrt{\D}\{\dfrac{d}{ds}\left(\dfrac{\d v}{\d
s}\right)\dfrac{\d u}{\d s} - \dfrac{d}{ds}\left(\dfrac{\d u}{\d s
}\right)\dfrac{\d v}{\d s} \nonumber\\&&+ \left( {2\brace 11}
\dfrac{du}{ds}+{2\brace12} \dfrac{dv}{ds}\right)\left(\dfrac{\d u}{\d
s}\right)^{2}\\&&\left(\left({ 1\brace 11}
 -{2\brace 12}
\right)\dfrac{du}{ds} - \left( {2\brace 22}
 - {1\brace 12}
\right)\dfrac{dv}{ds}\right)\dfrac{\d u}{\d s}\dfrac{\d v}{\d s}
\nonumber\\&& - \left({ 1\brace 21} \dfrac{du}{ds} + {1\brace 22}
\dfrac{dv}{ds}\right)\left(\dfrac{\d v}{\d s}\right)^{2}\}.\nonumber
\end{eqnarray}
\end{proposition}

Remark that the Christoffel symbols are expressed only by means of
the coefficients of the first fundamental form $\phi$ of surface
$S$ and their derivatives. It follows a very important result
obtained by Al Myller [34]:

\begin{theorem}
The geodesic curvature $G$ of a tangent versor field $(C,
\overline{\xi})$ on a surface $S$ is an intrisic invariant of $S$.
\end{theorem}

The invariant $G$ was named by Al. Myller the {\it deviation of
para\-llelism} of the tangent field $(C,\overline{\xi})$ on
surface $S$.

The expression (4.4.7) of $G$ can be simplified by introducing the
following notations
\begin{eqnarray}
&&\dfrac{D}{ds}\left(\dfrac{\d u^{i}}{\d s}\right) =
\dfrac{d}{ds}\left(\dfrac{\d u^{i}}{\d s}\right) +\sum_{j,k=1}^{2}{
i\brace jk}
 \dfrac{\d u^{j}}{\d s} \dfrac{d u^{k}}{ds}\nonumber\\&& (u= u^{1}; v =u^{2};
 i=1,2).
\end{eqnarray}
The operator (4.4.8) is the classical operator of covariant
derivative with respect to Levi-Civita connection.

Denoting by
\begin{equation}
\hspace*{10mm}G^{ij}(\d, d) = \sqrt{\D}\left\{\dfrac{D}{ds}\left(\dfrac{\d
u^{i}}{\d s }\right)\dfrac{\d u^{j}}{\d s} -
\dfrac{D}{ds}\left(\dfrac{\d u^{j} }{\d s}\right)\dfrac{\d u^{i}}{\d
s}\right\},(i,j=1,2)
\end{equation}
and remarking that $G^{ij}(\d,d) = -G^{ji}(\d,d)$ we have $G^{11}
= G^{22}=0$.

\begin{proposition}
The following formula holds:
\begin{equation}
G(\d, d) = G^{12}(\d,d).
\end{equation}
\end{proposition}

Indeed, the previous formula is a consequence of (4.4.7), (4.4.8)
and (4.4.9).

\begin{corollary}
The parallelism of tangent versor field $(C,\overline{\xi})$ to
$S$ along the curve $C$ is characterized by the differential
equation
\begin{equation}
G^{ij}(\d,d) = 0\;\; (i,j =1,2).
\end{equation}
\end{corollary}
In the following we introduce the notations \begin{equation} G(s)
= G(\d,d),\;\; K(s) = K(\d,d),\; T(s) = T(\d,d),
\end{equation}
since $\overline{\xi}(s) = \overline{r}_{u} \dfrac{\d u}{\d s} +
\overline{r}_{v}\dfrac{\d v}{\d s}$, $\overline{\a}(s) =
\overline{r}_{u}\dfrac{du}{ds}+\overline{r}_{v} \dfrac{dv}{ds}$.

\medskip

The second formula (4.3.5), by means of (4.4.5) gives us the
expression of normal curvature $K(\d,d)$ in the form
\begin{equation}
\hspace*{12mm}K(\d,d) = \dfrac{L du \d u +M(du \d v + dv \d u) + N dv \d
v}{\sqrt{Edu^{2} + 2F du dv + G dv^{2}}\sqrt{E\d u^{2} + 2F\d u \d
v + G\d v^{2}}}.
\end{equation}
It follows
\begin{equation}
K(\d,d) = K(d,\d)
\end{equation}

\begin{remark}\rm
The invariant $K(\d,d)$ can be written as follows
$$
K(\d,d) = \dfrac{\psi(\d,d)}{\sqrt{\phi(d,d)}\sqrt{\phi(\d,\d)}}
$$
where $\psi(\d,d)$ is the polar form of the second fundamental
form $\psi(d,d)$ of surface $S$.
\end{remark}

$K(\d,d)=0$ gives us the property of conjugation of $(C,
\overline{\xi})$ with $(C,\overline{\a})$.

The calculus of invariant $T(\d,d)$ can be made by means of
Weingarten formulae (4.2.3).

Since we have
\begin{eqnarray}
\dfrac{d\overline{\nu}}{ds} &=& \dfrac{1}{\D}\left\{[(FM -GL)
\overline{r}_{u} + (FL -
EM)\overline{r}_{v}]\dfrac{du}{ds}\nonumber\right.\\&&\left.+ \left[(FN
-GM)\overline{r}_{u} + (FM -
EN)\overline{r}_{v}\right]\dfrac{dv}{ds}\right\},
\end{eqnarray}
we deduce
\begin{eqnarray}
T(\d,d)&=& \dfrac{1}{\sqrt{\D}}\left\{\left[(FM - GL) \dfrac{du}{ds} +
(FN-GM)\dfrac{dv}{ds}\right]\dfrac{\d v}{\d s}\right.\nonumber\\&&\left.-
\left[(FL-EM)\dfrac{du}{ds} + (FM- EN)\dfrac{dv}{ds}\right]\dfrac{\d
u}{\d s}\right\}.
\end{eqnarray}
For simplicity we write $T(\d,d)$ from previous formula in the
following form
\begin{equation}
\hspace*{12mm}T(\d,d) =
\dfrac{1}{\sqrt{\D}\sqrt{\phi(d,d)}\sqrt{\phi(\d,\d)}}\left\|
\begin{array}{ccc}
E \d u\hspace{-0.5mm} +\hspace{-0.5mm} F\d v& F\d u\hspace{-0.5mm} +\hspace{-0.5mm}G \d v\\L du\hspace{-0.5mm} +\hspace{-0.5mm} Mdv & Mdu + N dv
\end{array} \right\|.
\end{equation}
Remember the expression of mean curvature $H$ and total curvature
$K_{t}$ of surface $S$:
\begin{equation}
H = \dfrac{EN - 2FM +GN}{2(EG -F^{2})},\;\; K_{t} = \dfrac{LN -
M^{2}}{EG-F^{2}}.
\end{equation}
It is not difficult to prove, by means of (4.4.16), the following
formula
\begin{equation}
T(\d, d) - T(d,\d) = 2\sqrt{\D} H\left(\dfrac{\d u}{\d
s}\dfrac{dv}{ds} - \dfrac{\d v}{\d s}\dfrac{du}{ds}\right).
\end{equation}
It has the following nice consequence:

\begin{theorem}
The geodesic torsion $T(\d,d)$ is symmetric with respect of
directions $\d$ and $d$, if and only if $S$ is a minimal surface.
\end{theorem}

Finally, $(C,\overline{\xi})$ is orthogonally conjugated with
$(C,\overline{\a})$ if and only if $T(\d,d) = 0$.

\begin{remark}\rm
The invariant $K(\d,d)$ was discovered by O. Mayer [20] and the
invariant $T(\d,d)$ was found by E. Bertoltti [3].
\end{remark}

\section{The Levi-Civita parallelism of vectors tangent to $S$}
\setcounter{theorem}{0}\setcounter{equation}{0}

The notion of Myller parallelism of vector field
$(C,\overline{V})$ tangent to $S$ along curve $C$, in the
associated Myller configuration $\mathfrak{M}_{t}(C,
\overline{\xi},\pi)$ is exactly the Levi-Civita parallelism of tangent
vector field $(C,\overline{V})$ to $S$ along the curve $C$. Indeed,
taking into account that the tangent versor field
$(C,\overline{\xi})$ has
\begin{equation}
\overline{\xi}(s) = \xi^{1}(s)\overline{r}_{u} +
\xi^{2}(s)\overline{r}_{v}
\end{equation}
and expression of the operator $\dfrac{D}{ds}$ is
\begin{equation}
\dfrac{D\xi^{i}}{ds} = \dfrac{d\xi^{i}}{ds} +
\sum_{j,k=1}^{2}\xi^{j}\left\{\begin{array}{cc} i\\jk
\end{array}\right\} \dfrac{du^{k}}{ds},\; (i=1,2)
\end{equation}
we have
$$
G^{ij} = \dfrac{D\xi^{i}}{ds}\xi^{j} -
\dfrac{D\xi^{j}}{ds}\xi^{i},\; (i=1,2).
$$
Thus the parallelism of versor $(C,\overline{\xi})$ along $C$ on
$S$ is expressed by $G^{ij} = 0.$ But these equations are
equivalent to
\begin{equation}
\dfrac{D(\l(s)\xi^{i}(s))}{ds} = 0\; (i=1,2),\; (\l(s)\neq 0).
\end{equation}
This is the definition of Levi-Civita parallelism of vectors
$\overline{V}(s)=\l(s)\overline{\xi}(s)$, $(\l(s)=\|V\|)$ tangent to $S$ along
$C$. Writing $\overline{V} = V^{1}(s) \overline{r}_{u} +
V^{2}(s)\overline{r}_{v}$ and putting
\begin{equation}
\dfrac{DV^{i}}{ds} =
\dfrac{dV^{i}}{ds}+\sum_{j,k=1}^{2}V^{j}{i\brace jk}
\dfrac{du^{k}}{ds}
\end{equation}
the Levi-Civita parallelism is expressed by
\begin{equation}
\dfrac{DV^{i}}{ds} = 0.
\end{equation}
In the case $\overline{V}(s)=\overline{\xi}(s)$ is a versor field,
parallelism of $(C,\overline{\xi})$ in the associated Myller
configurations $\mathfrak{M}_{t}(C, \overline{\xi},\pi)$ is called
the parallelism of directions tangent to $S$ along $C$. We can
express the condition of parallelism by means of invariants of the
field $(C,\overline{\xi})$, as follows:

1. {\it $(C,\overline{\xi})$ is parallel in the Levi-Civita sense on $S$
along $C$ iff $G(\d,d)=0.$}

2. {\it $(C, \overline{\xi})$ is parallel in the Levi-Civita sense on $S$
along $C$ iff the versor $\overline{\xi}(s^{\prime})$,
$\{|s^{\prime}-s|<\vp, \vp>0, s^{\prime}\in (s_{1},s_{2})\}$ is
parallel in the space $E_{3}$ with the normal plan $(P(s);
\overline{\xi}(s), \overline{\nu}(s))$.}

3. {\it A necessary and sufficient condition for the versor field $(C,
\overline{\xi})$ to be parallel on $S$ along $C$ in Levi-Civita sense
is that developing on a plane the ruled surface $E$ generated by
tangent planes field $(C,\pi)$ along $C$ to $E$- the directions
$(C,\overline{\xi})$ after developing to be parallel in Euclidean
sense.}

4. {\it The directions $(C,\overline{\xi})$ are parallel in the Levi-Civita
sense on $S$ along $C$ iff the versor field $(C,\overline{\xi}_{2})$ is normal to $S$.}

5. {\it $(C,\overline{\xi})$ is parallel on $S$ along $C$, iff $|K| =
K_{1}$.}

6. {\it If $(C,\xi)$ is parallel on $S$ along $C$, then $|K| = K_{1},
T=K_{2}$.}

7. {\it If the ruled surface $R(C,\overline{\xi})$ is not developable,
then $(C,\overline{\xi})$ is parallel in the Levi-Civita sense iff
$C$ is the striction line of surface $R(C,\overline{\xi})$.}

Taking into account the system of differential equations (4.5.5)
in the given initial conditions, we have assured the existence and
uniqueness of the Levi-Civita parallel versor fields on $S$ along
$C$.

Other results presented in Section 9, Chapter 2 can be
particularized here without difficulties.

A first application. The Tchebishev nets on $S$ are defined as a
net of $S$ for which the tangent lines to a family of curves of
net are parallel on $S$ along to every curve of another family of
curves of net and conversely.

\bigskip

\noindent{\bf A Bianchi result}

\begin{theorem}
In order that the net parameter $(u = u_{0}, v=v_{0})$ on $S$ to
be a Tchebishev is necessary and sufficient to have the following
conditions:
\begin{equation}
{1\brace 12}
 = {2\brace 12}
 = 0.
\end{equation}
\end{theorem}
Indeed, we have $G(\d,d) = G(d,\d) = 0$ if (4.5.6) holds.

But (4.5.6) is equivalent to the equations $\dfrac{\p E}{\p v} =
\dfrac{\p G}{\p u} = 0.$ So, with respect to a Tchebishev
parametrization of $S$ its arclength element $ds^{2} = \phi(d,d)$
is given by
$$
ds^{2} = E(u)du^{2} + 2F (u,v)du dv + G(v)dv^{2}.
$$

We finish this paragraph remarking that the parallelism of vectors
$(C, \overline{V})$ on $S$ or the concurrence of vectors
$(C,\overline{V})$ on $S$ can be studied using the corresponding
notions in configurations $\mathfrak{M}_{t}$ described in the
Chapter 3.

\section{The geometry of curves on a surface}
\setcounter{theorem}{0}\setcounter{equation}{0}

The geometric theory of curves $C$ embedded in a surface $S$ can
be derived from the geometry of tangent Myller configuration
$\mathfrak{M}_{t}(C,\overline{\a},\pi)$ in which
$\overline{\a}(s)$ is tangent versor to $C$ at point $P(s)\in C$
and $\pi(s)$ is tangent plane to $S$ at point $P$ for any $s\in
(s_{1},s_{2})$.

Since $\mathfrak{M}_{t}(c,\overline{\a},\pi)$ is geometrically
associated to the curve $C$ on $S$ we can define its Darboux frame
$\cal{R}_{D}$, determine the moving equations of $\cal{R}_{D}$ and
its invariants, as belonging to curve $C$ on surface $S$.

Applying the results established in Chapter 3, first of all we
have

\begin{theorem}
For a smooth curve $C$ embedded in a surface $S$, there exists a
Darboux frame $\cal{R}_{D} = (P(s); \overline{\a}(s),
\overline{\mu}^{*}(s), \overline{\nu}(s))$ and a system of
invariants $\k_{g}(s), \k_{n}(s)$ and $\tau_{g}(s)$, satisfying
the following moving equations
\begin{equation}
\dfrac{d\overline{r}}{ds} = \overline{\a}(s),\;\; \forall s\in
(s_{1},s_{2})
\end{equation}
and
\begin{eqnarray}
\dfrac{d\overline{\a}}{ds} &=& \k_{g}(s)\overline{\mu}^{*}
+\k_{n}(s)\overline{\nu},\nonumber\\\dfrac{d\overline{\mu}^{*}}{ds}
&=& -\k_{g}(s)\overline{\a} +
\tau_{g}(s)\overline{\nu},\\\dfrac{d\overline{\nu}}{ds}&=&
-\k_{n}(s)\overline{\a} - \tau_{g}(s)\overline{\mu}^{*},\;\;
\forall s\in (s_{1},s_{2}).\nonumber
\end{eqnarray}
\end{theorem}

The functions $\k_{g}(s), \k_{n}(s), \tau_{g}(s)$ are called the
{\it geodesic curvature}, the {\it normal curvature} and {\it
geodesic torsion} of curve $C$ at point $P(s)$ on surface $S,$
respectively.

Exactly as in Chapter 3, a fundamental theorem can be enounced and
proved.

The invariants $\k_{g}, \k_{n}$ and $\tau_{g}$ have the same
values along $C$ on any smooth surface $S^{\prime}$ which passes
throught curve $C$ and is tangent to surface $S$ along
$C$.

The geometric interpretations of these invariants and the cases of
curves $C$ for which $\k_{g}(s) = 0,$ or $\k_{n}(s) = 0$ or
$\tau_{g}(s)=0$ can be studied as in Chapter 3.

The curves $C$ on $S$ for which $\k_{g}(s) = 0$ are called (as in
Chapter 3) geodesics (or {\it autoparallel curves}) of $S$. If $C$
has the property $\k_{n}(s) = 0,$ $\forall s\in (s_{1},s_{2})$, it
is {\it asymptotic curve} of $S$ and the curve $C$ for which
$\tau_{g}(s) = 0,$ $\forall s\in(s_{1},s_{2})$ is the {\it
curvature line} of $S$.

The expressions of these invariants can be obtained from those of
the invariants $G(\d,d)$, $K(\d,d)$ and $T(\d,d)$ for
$\overline{\xi}(s) = \overline{\a}(s)$ given in Chapter 4:
\begin{equation}
\k_{g} = G(d,d),\; \k_{n}(s) = K(d,d),\;\; \tau_{g}(s) = T(d,d).
\end{equation}
So, we have
\begin{equation}
\hspace*{10mm}\begin{array}{lll}
\k_{g}&=&
\sqrt{\D}\left\{\dfrac{D}{ds}\left(\dfrac{du^{i}}{ds}\right)\dfrac{du^{j}}{ds}-\dfrac{D}{ds}
\left(\dfrac{du^{j}}{ds}\right)\dfrac{du^{i}}{ds}\right\},\ (i,j=1,2)\vspace*{2mm}\\
\k_{n}&=& \dfrac{L du^{2} + 2M du dv +N dv^{2}}{E du^{2} + 2F du dv
+ G dv^{2}},\vspace*{2mm}\\ \tau_{g}&=&  \dfrac{1}{\sqrt{\D}}
\dfrac{\left\|\begin{array}{cc}L du + Mdv& Mdu + N dv\vspace*{2mm}\\ Edu + Fdv &
Fdu + Gdv
\end{array}\right\|}{Edu^{2} + 2F du dv + Gdv^{2}}.
\end{array}
\end{equation}
Evidently, $\k_{g}$ is an intrinsic invariant of surface $S$ and
$\k_{n}$ is the ratio of fundamental forms $\psi$ and $\Phi$.

The Mark Krein formula (3.4.3), Chapter 3 gives now the
Gauss-Bonnet formula for a surface $S$.

\newpage

\section{The formulae of O. Mayer and E. Bortolotti}
\setcounter{theorem}{0}\setcounter{equation}{0}

The geometers O. Mayer and E. Bortolotti gave some new forms of
the invariants $K(\d,d)$ and $T(\d,d)$ which generalize the Euler
or Bonnet formulae from the geometry of surfaces in Euclidean
space $E_{3}$.

Let $S$ be a smooth surface in $E_{3}$ having the parametrization
given by curvature lines. Thus the coefficients $F$ and $M$ of the
first fundamental and the second fundamental form vanish.

We denote by $\t = \sphericalangle \(\overline{\a},
\dfrac{\overline{r}_{u}}{\sqrt{E}}\)$ and obtain
\begin{equation}
\cos \t = \dfrac{\sqrt{E}du}{\sqrt{Edu^{2} + G dv^{2}}},\;\; \sin
\t = \dfrac{\sqrt{G}dv}{\sqrt{Edu^{2} + G dv^{2}}}.
\end{equation}
The principal curvature are expressed by
\begin{equation}
\dfrac{1}{R_{1}} = \dfrac{L}{E},\;\; \dfrac{1}{R_{2}} = \dfrac{N}{G}.
\end{equation}

The mean curvature and total curvature are
\begin{eqnarray}
H&=&\dfrac{1}{2}\left(\dfrac{1}{R_{1}} + \dfrac{1}{R_{2}}\right) =
\dfrac{1}{2}\left(\dfrac{L}{E} +
\dfrac{N}{G}\right)\nonumber\\K_{t}&=&\dfrac{1}{R_{1}}\dfrac{1}{R_{2}}
= \dfrac{LN}{EG}.
\end{eqnarray}
Consider a tangent versor field $(C,\overline{\xi})$,
$\overline{\xi} = \overline{r}_{u}\dfrac{\d u}{\d s} +
\overline{r}_{v}\dfrac{\d v}{\d s}$ and let be $\s =
\sphericalangle
\(\overline{\xi},\dfrac{\overline{r}_{u}}{\sqrt{E}}\)$. One gets
\begin{equation}
\cos \s = \dfrac{\sqrt{E}\d u}{\sqrt{E\d u^{2} + G\d v^{2}}},\;
\sin \s = \dfrac{\sqrt{G}\d v}{\sqrt{E \d u^{2} + G\d v^{2}}}.
\end{equation}
Consequently, the expressions (4.4.13) and (4.4.17) of
the normal \linebreak curvature and geodesic torsion of
$(C,\overline{\xi})$ on $S$ are as follows
\begin{equation}
K(\d,d) = \dfrac{\cos \s \cos \t}{R_{1}} + \dfrac{\sin \s \sin
\t}{R_{2}}
\end{equation}
\begin{equation}
T(\d,d) = \dfrac{\cos \s\sin \t}{R_{2}} - \dfrac{\sin \s \cos
\t}{R}.
\end{equation}
The first formula was established by O. Mayer [20] and second
formula was given by E. Bortolotti [3].

In the case $\overline{\xi}(s) = \overline{\a}(s)$ these formulae
reduce to the known Euler and Bonnet formulas, respectively:
\begin{eqnarray}
\k_{n}&=& \dfrac{\cos^{2}\t}{R_{1}} +
\dfrac{\sin^{2}\t}{R_{2}}\nonumber\\\tau_{g} &=&
\dfrac{1}{2}\left(\dfrac{1}{R_{2}} - \dfrac{1}{R_{1}}\right)\sin
2\t.
\end{eqnarray}
For $\t =0$ or $\t = \dfrac{\pi}{2},\k_{n}$ is equal to
$\dfrac{1}{R_{1}}$ and $\dfrac{1}{R_{2}}$ respectively. For $\t =
\pm \dfrac{\pi}{4}$, $\tau_{g}$ takes the extremal values.
\begin{equation}
\dfrac{1}{T_{1}} = \dfrac{1}{2}\left(\dfrac{1}{R_{2}} -
\dfrac{1}{R_{1}}\right),\; \dfrac{1}{T_{2}} =
-\dfrac{1}{2}\left(\dfrac{1}{R_{2}} - \dfrac{1}{R_{1}}\right).
\end{equation}
Thus
\begin{equation}
T_{m} = \dfrac{1}{2}\left(\dfrac{1}{T_{1}} +
\dfrac{1}{T_{2}}\right),\; T_{t} = \dfrac{1}{T_{1}}\dfrac{1}{T_{2}}
\end{equation}
are the {\it mean torsion} of $S$ at point $P\in S$ and the {\it total
torsion} of $S$ at point $P\in S.$

For surfaces $S$, $T_{m}$ and $T_{t}$ have the following
properties:
\begin{equation}
T_{m} = 0,\;\; T_{t} = -\dfrac{1}{4}\left(\dfrac{1}{R_{1}} -
\dfrac{1}{R_{2}}\right)^{2}.
\end{equation}

\begin{remark}\rm
1. As we will see in the next chapter the nonholomorphic manifolds in
$E_{3}$ have a nonvanishing mean torsion $T_{m}$.

2. $T_{t}$ from (4.7.10) gives us the Bacaloglu curvature of
surfaces [35].
\end{remark}

Consider in plane $\pi(s)$ the cartesian orthonormal frame $(P(s);
\overline{i}_{1}, \overline{i}_{2})$, $i_{1} =
\dfrac{\overline{r}_{u}}{\sqrt{E}}, i_{2} =
\dfrac{\overline{r}_{v}}{\sqrt{G}}$ and the point $Q\in \pi(s)$
with the coordinates $(x,y)$, given by
$$
\overrightarrow{PQ} = x\overline{i}_{1} + y \overline{i}_{2},\;
\overrightarrow{PQ} = |\k_{n}|^{-1}\overline{\a}.
$$
But $\overline{\a} =\cos \t \overline{i}_{1} +\sin \t i_{2}$. So
we have the coordinates $(x,y)$ of point $Q$:
\begin{equation}
x = |\k_{n}|^{-1}\cos \t; \;\; y = |\k_{n}|^{-1}\sin \t.
\end{equation}

The locus of the points $Q$, when $\theta$ is variable in interval
$(0,2\pi)$ is obtained eliminating the variable $\t$ between the
formulae (4.7.7) and (4.7.11). One obtains a pair of conics:
\begin{equation}
\dfrac{x^{2}}{R_{1}} + \dfrac{y^{2}}{R_{2}} = \pm 1
\end{equation}
called the {\it Dupin indicatrix} of normal curvatures. This indicatrix
is important in the local study of surfaces $S$ in Euclidean
space.

Analogously we can introduce the Bonnet indicatrix. Consider in
the plane $\pi(s)$ tangent to $S$ at point $P(s)$ the frame
$(P(s), \overline{i}_{1}, \overline{i}_{2})$ and the point
$Q^{\prime}$ given by
$$
\overrightarrow{PQ}^{\prime} = |\tau_{g}|^{-1}\overline{\a} =
x\overline{i}_{1} +y \overline{i}_{2}.
$$
Then, the locus of the points $Q^{\prime}$, when $\t$ verifies
(4.7.7) and $x =|\tau_{g}|^{-1}\cos \t$, $y = |\tau_{g}|^{-1}\sin
\t$, defines the Bonnet indicatrix of geodesic torsions:
\begin{equation}
\left(\dfrac{1}{R_{2}} - \dfrac{1}{R_{1}}\right)xy = \pm 1
\end{equation}
which, in general, is formed by a pair of conjugated equilateral
hyperbolas.

Of course, the relations between  the indicatrix of Dupin and
Bonnet can be studied without difficulties.

Following the same way we can introduce the indicatrix of the
invariants $K(\d,d)$ and $T(\d,d)$.

So, consider the angles $\t = \sphericalangle
(\overline{\a},i_{1})$, $\s =
\sphericalangle(\overline{\xi},\overline{i}_{1})$ and $U\in
\pi(s)$. The point $U$ has the coordinates $x,y$ with respect to
frame $(P(s); \overline{i}_{1}, i_{2})$ given by
\begin{equation}
x= |K(\d, d)|^{-1}\cos \s,\;\; y = |K(\d,d)|^{-1} \sin \s.
\end{equation}
The locus of points $U$ is obtained from (4.7.5) and (4.7.14):
\begin{equation}
\dfrac{x\cos \t}{R_{1}} +\dfrac{y\sin \t}{R_{2}} = \pm 1.
\end{equation}
Therefore (4.7.15) is the indicatrix of the normal curvature
$K(\d,d)$ of versor field $(C,\overline{\xi}).$ It is a pair of
parallel straight lines.

Similarly, for
$$
x = |T(\d,d)|^{-1}\cos \s,\;\; y = |T(\d,d)|^{-1}\sin \s
$$
and the formula (4.7.6) we determine the indicatrix of geodesic
torsion $T_{g}(\d,d)$ of the versor field $(C,\overline{\xi}):$
$$
\dfrac{x\sin \t}{R_{2}} - \dfrac{y\cos \t}{R_{1}} = \pm 1.
$$
It is a pair of parallel straight lines.

Finally, we can prove without difficulties the following formula:
\begin{equation}
\hspace*{14mm}\k_{n}(\t)\k_{n}(\s)\hspace{-0.5mm}+\hspace{-0.5mm} \tau_{g}(\t) \tau_{g}(\s)\hspace{-0.5mm}=\hspace{-0.5mm}2HK(\s,\t)\cos
(\s-\t)\hspace{-0.5mm}-\hspace{-0.5mm}K_{t}\cos 2(\s-\t),
\end{equation}
where $\k_{n}(\t)=  \k_{n}(d,d),$ $\k_{n}(\s) = \k_{n}(\d,\d)$,
$K(\s,\t) = K(d,\d)$.

For $\s = \t,$ (i.e. $\overline{\xi} = \overline{\a}$) the
previous formulas leads to a known Baltrami-Enneper formula:
\begin{equation}
\k_{n}^{2} + \tau_{g}^{2} - 2H \k_{n} +K_{t} = 0,
\end{equation}
for every point $P(s)\in S.$

Along the asymptotic curves we have $\k_{n} = 0,$ and the previous
equations give us the Enneper formula
$$
\tau_{g}^{2} +K_{t} = 0,\; (\k_{n}=0).
$$

Along to the curvature lines, $\tau_{g}=0$ and one obtains from
(4.7.17)
$$
\k_{n}^{2} - 2H \k_{n} +K_{t} = 0,\;\; (\tau_{g} = 0).
$$

The considerations made in Chapter 4 of the present book allow
to affirm that the applications of the theory of Myller
configurations the geometry of surfaces in Euclidean space $E_{3}$
are interesting. Of course, the notion of Myller configuration can
be extended to the geometry of nonholomonic manifolds in $E_{3}$
which will be studied in next chapter. It can be applied to the
theory of versor fields in $E_{3}$ which has numerous applications
to Mechanics, Hydrodynamics (see the papers by Gh. Gheorghiev and
collaborators).

Moreover, the Myller configurations can be defined and
investigated in Riemannian spaces and applied to the geometry of
submanifolds of these spaces, [24], [25]. They can be studied in the Finsler, Lagrange or Hamilton spaces [58], [68], [70].

\chapter{Applications of theory of Myller configurations in the geometry of nonholonomic manifolds from $E_{3}$}

The second efficient applications of the theory of Myller
configurations can be done in the geometry of nonholonomic
manifolds $E_{3}^{2}$ in $E_{3}$. Some important results, obtained
in the geometry of manifolds $E_{3}^{2}$ by Issaly l$'$Abee, D.
Sintzov [40], Gh. Vr\u{a}nceanu [44], [45], Gr. Moisil [72], M. Haimovici [14], [15] Gh. Gheorghiev [10], [11], I. Popa [12],
G. Th. Gheorghiu [13], R. Miron [30], [64], [65], I. Creang\u{a} [4], [5], [6], A. Dobrescu [7], [8] and I. Vaisman [41], [42], can be unitary presented by means of
associated Myller configuration to a curve embedded in
$E_{3}^{2}$. Now, a number of new concepts appears, the mean
torsion, the total torsion, concurrence of tangent vector field.
The new formulae, as Mayer, Bortolotti-formulas, Dupin and Bonnet
indicatrices etc. will be studied, too.

\section{Moving frame in Euclidean spaces $E_{3}$}
\setcounter{theorem}{0}\setcounter{equation}{0}

Let $R = (P; I_{1}, I_{2}, I_{3})$ be an orthonormed positively
oriented frame in $E_{3}$. The application $P\in E_{3}\to R = (P;
I_{1}, I_{2}, I_{3})$ of class $C^{k}$, $k\geq 3$ is a moving frame
of class $C^{k}$ in $E_{3}$. If $\overline{r} =
\overrightarrow{OP} = x \overline{e}_{1} + y \overline{e}_{2} +
z\overline{e}_{3}$ is the vector of position of point $P$, and $P$
is the application point of the versors $I_{1}, I_{2}, I_{3}$,
thus the moving equation of $R$ can be expressed in the form (see
Spivak, vol. I [76] and Biujguhens [2]):
\begin{equation}
d\overline{r} = \o_{1}I_{1} + \o_{2}I_{2} + \o_{3}I_{3},
\end{equation}
where $\o_{i}(i=1,2,3)$ are independent 1 forms of class $C^{k-1}$
on $E_{3}$, and
\begin{equation}
dI_{i} = \sum_{j=1}^{3}\o_{ij}I_{j}, \;\; (i=1,2,3),
\end{equation}
with $\o_{ij}$, $(i,j = 1,2,3)$ are the rotation Ricci
coefficients of the frame $R$. They are 1-forms of class
$C^{k-1}$, satisfying the skewsymmetric conditions
\begin{equation}
\o_{ij} +\o_{ji} = 0,\;\; (i,j = 1,2,3).
\end{equation}

In the following, it is convenient to write the equations (5.1.2)
in the form
$$
\begin{array}{c}
dI_{1} = rI_{2} -qI_{3}\\dI_{2} = p I_{3} - r I_{1}\\dI_{3} = q
I_{1} - p I_{2}.
\end{array}\leqno{(5.1.2)^{\prime}}
$$

Thus, {\it thus structure equations} of the moving frame $R$ can
be obtained by exterior differentiating the equations (5.1.1) and
(5.1.2)$^{\prime}$ modulo the same system of equations (5.1.1),
(5.1.2)$^{\prime}$.

One obtains, without difficulties

\begin{theorem}
The structure equations of the moving frame $R$ are
\begin{eqnarray}
&&d\o_{1} = r \wedge \o_{2},-q\wedge \o_{3}, \;\; dp = r\wedge
q,\nonumber\\&&d\o_{2} = p\wedge \o_{3} - r \wedge \o_{1},\;\; dq
= p\wedge r ,\\&&d\o_{3} = q\wedge \o_{1} - p\wedge \o_{2},\;\; dr
= q\wedge p.\nonumber
\end{eqnarray}
\end{theorem}

In the vol. II of the book of Spivak [76], it is proved the theorem
of existence and uniqueness of the moving frames:

\begin{theorem}
Let $(p,q,r)$ be 1-forms of class $C^{k-1}$ ,$(k\geq 3)$ on $E_{3}$
which satisfy the second group of structure equation $(5.1.4)$,
thus:

$1^{\circ}.$ In a neighborhood of point $O\in E_{3}$ there is a
triple of vector fields $(I_{1}, I_{2}, I_{3})$, solutions of
equations (5.1.2)$^{\prime}$, which satisfy initial conditions
$(I_{1}(0), I_{2}(0), I_{3}(0))$-positively oriented, orthonormed
triple in $E_{3}$.

$2^{\circ}.$ In a neighborhood of point $O\in E_{3}$ there exists
a moving frame $R = (P; I_{1}, I_{2}, I_{3})$ orthonormed
positively oriented for which $\overline{r} = \overrightarrow{OP}$
is given by (5.1.1), where $\o_{i} (i=1,2,3)$ satisfy the first
group of equations (5.1.4). $(I_{1}, I_{2}, I_{3})$ are given in
$1^{\circ}$ and the initial conditions are verified: $(P_{0} =
\overrightarrow{OO}; I_{1}(0), I_{2}(0), I_{3}(0))$-the
orthonormed frame at point $O\in E_{3}$.\end{theorem}

Remarking that the 1-forms $\o_{1}, \o_{2}, \o_{3}$ are
independent we can express 1-forms $p,q,r$ with respect to
$\o_{1}, \o_{2}, \o_{3}$ in the following form:
\begin{eqnarray}
&&p = p_{1}\o_{1} + p_{2}\o_{2} + p_{3}\o_{3},\ q = q_{1}\o_{1} +
q_{2}\o_{2} + q_{3}\o_{3},\\&& r = r_{1}\o_{1} + r_{2}\o_{2} +
r_{3}\o_{3}.\nonumber
\end{eqnarray}

The coefficients $p_{i}, q_{i}, r_{i}$ are function of class
$C^{k-1}$ on $E_{3}$.

Of course we can write the structure equations (5.1.4) in terms of
these coefficients. Also, we can state the Theorem 5.1.2 by means of
coefficients of 1-forms $p,q,r$.

The 1-forms $p,q,r$ determine the rotation vector $\O$ of the
moving frame:
\begin{equation}
\O = p I_{1} + q I_{2} + r I_{3}.
\end{equation}
Its orthogonal projection on the plane $(P; I_{1}, I_{2})$ is
\begin{equation}
\overline{\t} = p I_{1} + q I_{2}.
\end{equation}
Let $C$ be a curve of class $C^{k}$, $k\geq 3$ in $E_{3}$, given
by
\begin{equation}
\overline{r} = \overline{r}(s), \;\; s\in(s_{1},s_{2}),
\end{equation}
where $s$ is natural parameter. If the origin $P$ of moving frame
describes the curve $C$ we have $R = (P; I_{1}(s), I_{2}(s),
I_{3}(s))$, where $P(s) = P(\overline{r}(s))$, $I_{j}(s) =
I_{j}(\overline{r}(s))$, $(j=1,2,3)$.

So, the tangent versor $\dfrac{d\overline{r}}{ds}$ to curve $C$ is:
\begin{equation}
\dfrac{d\overline{r}}{ds} = \dfrac{\o_{1}(s)}{ds}I_{1} +
\dfrac{\o_{2}(s)}{ds}I_{2} + \dfrac{\o_{3}(s)}{ds}I_{3}
\end{equation}
and $ds^{2}=\langle d\overline{r},d\overline{r}\rangle$ is:
\begin{equation}
ds^{2} = (\o_{1}(s))^{2} + (\o_{2}(s))^{2} + (\o_{3}(s))^{2},
\end{equation}
where $\o_{i}(s)$ are 1-forms $\o_{i}$ restricted to $C$.

The restrictions $p(s), q(s), r(s)$ of the 1-forms $p,q,r$ to $C$
give us:
\begin{equation}
\begin{array}{lll}
\dfrac{p(s)}{ds} &=& p_{1}(s)\dfrac{\o_{1}(s)}{ds} +
p_{2}(s)\dfrac{\o_{2}(s)}{ds}
+p_{3}(s)\dfrac{\o_{3}(s)}{ds},\\\spa\dfrac{q(s)}{ds} &=&
q_{1}(s)\dfrac{\o_{1}(s)}{ds} + q_{2}(s)\dfrac{\o_{2}(s)}{ds} +
q_{3}(s)\dfrac{\o_{3}(s)}{ds},\\\spa\dfrac{r(s)}{ds}&=&r_{1}(s)\dfrac{\o_{1}(s)}{ds}
+ r_{2}(s)\dfrac{\o_{2}(s)}{ds} +
r_{3}(s)\dfrac{\o_{3}(s)}{ds}.
\end{array}
\end{equation}

\section{Nonholonomic manifolds $E_{3}^{2}$}
\setcounter{theorem}{0}\setcounter{equation}{0}

\begin{definition}
A nonholonomic manifold $E_{3}^{2}$ on a domain $D\subset E_{3}$
is a nonintegrable distribution $\cal{D}$ of dimension 2, of class
$C^{k-1},$ $k\geq 3$.
\end{definition}

One can consider $\cal{D}$ as a plane field $\pi(P)$, $P\in D,$ the
application $P\to \pi(P)$ being of class $C^{k-1}$.

Also $\cal{D}$ can be presented as the plane field $\pi(P)$
orthogonal to a versor field $\overline{\nu}(P)$ normal to the
plane $\pi(P)$, $\forall P\in D.$

Assuming that $\overline{\nu}(P)$ is the versor of vector
$$\overline{V}(P) = X(x,y,z)\overline{e}_{1} +
Y(x,y,z)\overline{e}_{2} + Z(x,y,z)\overline{e}_{3}$$ and
$\overrightarrow{OP} = x \overline{e}_{1} + y \overline{e}_{2} +z
\overline{e}_{3}$ thus the nonholonomic manifold $E_{3}^{2}$ is
given by the Pfaff equation:
\begin{equation}
\o = X(x,y,z)dx + Y(x,y,z)dy + Z(x,y,z)dz = 0
\end{equation}
which is nonintegrable, i.e.
\begin{equation}
\o\wedge d\o \neq0.
\end{equation}
Consider a moving frame $\cal{R} = (P; I_{1}, I_{2}, I_{3})$ in
the space $E_{3}$ and the nonholonomic manifold $E_{3}^{2}$, on a
domain $D$ orthogonal to the versors field $I_{3}$. It is given by
the Pfaff equations
\begin{equation}
\o_{3} = \langle I_{3}, d\overline{r}\rangle = 0,\;\;\;
\mbox{on}\; D.
\end{equation}
By means of (5.1.4) we obtain
\begin{equation}
\o_{3}\wedge d\o_{3} = -(p_{1} +q_{2})\o_{1}\wedge \o_{2}\wedge
\o_{3}.
\end{equation}
So, the Pfaff equation $\o_{3} = 0$ is not integrable if and only
if we have
\begin{equation}
p_{1} + q_{2}\neq 0,\; \mbox{on}\; D.
\end{equation}

It is very known that:

In the case $p_{1} + q_{2} =0$ on the domain $D$, there are two
functions $h(x,y,z)\neq 0$ and $f(x,y,z)$ on $D$ with the property
\begin{equation}
h\o_{3} = df.
\end{equation}
Thus the equation $h\o_{3} = 0,$ have a general solution
\begin{equation}
f(x,y,z) = c,\;\;(c=const.),\;\; P(x,y,z)\in D.
\end{equation}
The manifold $E_{3}^{2}$, in this case is formed by a family of
surfaces, given by (5.2.7).

For simplicity we assume that the class of the manifold
$E_{3}^{2}$ is $C^{\infty}$ and the same property have the
geometric object fields or mappings which will be used in the
following parts of this chapter.

A smooth curve $C$ embedded in the nonholonomic manifolds
$E_{3}^{2}$, has the tangent vector $\overline{\a} =
\dfrac{d\overline{r}}{ds}$ given by (5.1.9) and $\o_{3} = 0.$

This is
\begin{equation}
\overline{\a}(s) = \dfrac{d\overline{r}}{ds}
=\dfrac{\o_{1}(s)}{ds}I_{1} + \dfrac{\o_{2}(s)}{ds}I_{2},\; \forall
P(\overline{r}(s))\in D.
\end{equation}
The square of arclength element, $ds$ is
\begin{equation}
ds^{2} = (\o_{1}(s))^{2} + (\o_{2}(s))^{2}.
\end{equation}
And the arclength of curve $C$ on the interval $[a,b]\subset
(s_{1},s_{2})$ is \begin{equation} s =
\int_{a}^{b}\sqrt{(\o_{1}(\s))^{2} + (\o_{2}(\s))^{2}}d\s.
\end{equation}
A tangent versor field $\overline{\xi}(s)$ at the same point
$P(s)$ to another curve $C^{\prime},$ $P\in C^{\prime}$ and $P\in
C,$ can be given in the same form (5.2.8):
\begin{equation}
\overline{\xi}(s) = \dfrac{\d\overline{r}}{\d s}  =
\dfrac{\widetilde{\o}_{1}(s)}{\d s}I_{1} +
\dfrac{\widetilde{\o}_{2}(s)}{\d s}I_{2},
\end{equation}
with $\d s^{2} = (\widetilde{\o}_{1}(s))^{2} +
(\widetilde{\o}_{2}(s))^{2}$.

Thus, $\overline{\xi}(s)$ is a versor field tangent to the
nonholonomic manifold $E_{3}^{2}$ along to the curve $C$.

To any tangent vector field $(C,\overline{\xi})$ to $E_{3}^{2}$
along the curve $C$ belonging to $E_{3}^{2}$ we uniquely
associated the tangent Myller configuration $\mathfrak{M}_{t}(C,
\overline{\xi}, I_{3})$. We say that: {\it the geometry of the
associated tangent Myller configuration $\mathfrak{M}_{t}(C,
\overline{\xi},I_{3})$ is the geometry of versor field $(C,
\overline{\xi})$ on $E_{3}^{2}$}. In particular the geometry of
configuration $\mathfrak{M}_{t}(C,\overline{\a},I_{3})$ is {\it
the geometry of curve $C$} in the nonholonomic manifold
$E_{3}^{2}$.

The Darboux frame of $(C,\overline{\xi})$ on $E_{3}^{2}$ is
$\cal{R}_{D} = (P; \overline{\xi}, \overline{\mu}, I_{3})$, with
$\overline{\mu} = I_{3}\times \overline{\xi}$, i.e.:
\begin{equation}
\overline{\mu} = \dfrac{\widetilde{\o}_{1}(s)}{\d s}I_{2} -
\dfrac{\widetilde{\o}_{2}(s)}{\d s}I_{1}.
\end{equation}
In Darboux frame $\cal{R}_{D}$, the tangent versor $\overline{\a}
= \dfrac{d\overline{r}}{ds}$ to $C$ can be expressed in the
following form:
\begin{equation}
\overline{\a}(s) = \cos \l(s) \overline{\xi}(s) + \sin
\l(s)\overline{\mu}(s).
\end{equation}
Taking into account the relations (5.2.8), (5.2.11) and (5.2.13),
one gets:
\begin{eqnarray}
\dfrac{\o_{1}}{ds}& =& \cos \l \dfrac{\widetilde{\o}_{1}}{\d s} -
\sin
\l \dfrac{\widetilde{\o}_{2}}{\d s},\nonumber\\
\dfrac{\o_{2}}{ds}&=& \sin \l \dfrac{\widetilde{\o}_{1}}{\d s} +
\cos \l \dfrac{\widetilde{\o}_{2}}{\d s}.
\end{eqnarray}
For $\l(s) = 0$, we have $\overline{\a}(s) = \overline{\xi}(s)$.

\section{The invariants $G,R,T$ of a tangent versor field $(C,\overline{\xi})$ in $E_{3}^{2}$}
\setcounter{theorem}{0}\setcounter{equation}{0}

The invariants $G,K,T$ of $(C,\overline{\xi})$ satisfy the
fundamental equations
\begin{equation}
\dfrac{d\overline{\xi}}{ds} = G \overline{\mu}+K I_{3},\
\dfrac{d\overline{\mu}}{ds} = -G \overline{\xi}+T I_{3},\
\dfrac{dI_{3}}{ds} = -K \overline{\xi} - T \overline{\mu}.
\end{equation}
Consequently,
\begin{eqnarray}
G(\d,d) &=& \left\langle \overline{\xi}, \dfrac{d\overline{\xi}}{ds},
I_{3} \right\rangle\nonumber\\
K(\d,d)&=& \left\langle\dfrac{d\overline{\xi}}{ds}, I_{3}\right\rangle =
-\left\langle \overline{\xi}, \dfrac{dI_{3}}{ds}\right\rangle\\T(\d,d)&=&
\left\langle \overline{\xi}, I_{3}, \dfrac{dI_{3}}{ds}\right\rangle.\nonumber
\end{eqnarray}

The proofs and the geometrical meanings are similar to those from
Chapter 3.

By means of expression (5.2.11) of versors $\overline{\xi}(s)$ we
deduce:
\begin{eqnarray}
\dfrac{d\overline{\xi}}{ds}&=& \left[
\dfrac{d}{ds}\left(\dfrac{\widetilde{\o}_{1}}{\d s}\right) -
\dfrac{r}{ds}\dfrac{\widetilde{\o}_{2}}{\d s} \right]I_{1} + \left[
\dfrac{d}{ds}\left(\dfrac{\widetilde{\o}_{2}}{\d s}\right) +
r\dfrac{\widetilde{\o}_{1}}{\d s} \right]I_{2} + \nonumber\\&&+
\left[ \dfrac{p}{ds}\dfrac{\widetilde{\o}_{2}}{\d s} - \dfrac{q}{ds}
\dfrac{\widetilde{\o}_{1}}{\d s} \right]I_{3},
\end{eqnarray}
and
\begin{equation}
\dfrac{dI_{3}}{ds} = \left(q_{1}\dfrac{\o_{1}}{ds} +
q_{2}\dfrac{\o_{2}}{ds}\right)I_{1} - \left(p_{1}\dfrac{\o_{1}}{ds}
+ p_{2}\dfrac{\o_{2}}{ds} \right)I_{2}.
\end{equation}

The formulae (5.3.2) lead to the following expressions for {\it
geodesic curvature} $G(\d,d)$, {\it normal curvature} $K(\d,d)$
and {\it geodesic torsion} $T(\d,d)$ of the tangent versor field
$(C,\overline{\xi})$ on the nonholonomic manifold $E_{3}^{2}$:
\begin{equation}
\hspace*{8mm}G(\d,d) = \dfrac{\widetilde{\o}_{1}}{\d
s}\left[\dfrac{d}{ds}\left(\dfrac{\widetilde{\o}_{2}}{\d s }\right)
+ \dfrac{r}{ds}\dfrac{\widetilde{\o}_{1}}{\d s} \right] -
\dfrac{\widetilde{\o}_{2}}{\d
s}\left[\dfrac{d}{ds}\left(\dfrac{\widetilde{\o}_{1}}{\d s } -
\dfrac{r}{ds}\dfrac{\widetilde{\o}_{2}}{\d s} \right)\right],
\end{equation}
\begin{equation}
K(\d,d) = \dfrac{p\widetilde{\o}_{2} - q\widetilde{\o}_{1}}{\d
s\cdot d s}
\end{equation}
\begin{equation}
T(\d,d) = \dfrac{p\widetilde{\o}_{1} + q\widetilde{\o}_{2}}{\d s
ds}.
\end{equation}

All these formulae take very simple forms if we consider the
angles $\a = \sphericalangle (\overline{\a},I_{1})$, $\b =
\sphericalangle (\overline{\xi},I_{1})$ since we have
\begin{eqnarray}
&&\dfrac{\o_{1}}{ds} = \cos \a,\ \dfrac{\o_{2}}{ds} = \sin
\a,\nonumber\\&& \dfrac{\widetilde{\o}_{1}}{\d s} = \cos \b,\
\dfrac{\widetilde{\o}_{2}}{\d s} = \sin \b.
\end{eqnarray}
With notations
\begin{equation}
\hspace*{8mm}G(\d,d) = G(\b,\a); K(\d,d) = K(\b,\a), T(\d,d) = T(\b,\a)
\end{equation}
the formulae (5.3.2) can be written:
\begin{equation}
\hspace*{12mm}G(\b,\a) = \dfrac{d\overline{\b}}{ds} +r_{1}\cos \a +r_{2}\sin \a,
\end{equation}
\begin{equation}
\hspace*{14mm}K(\b,\a) = (p_{1}\cos \a + p_{2}\sin \a)\sin \b - (q_{1}\cos \a
+q_{2}\sin \a)\cos \b,
\end{equation}
\begin{equation}
\hspace*{14mm}T(\b,\a) = (p_{1}\cos \a +p_{2}\sin \a)\cos \b + (q_{1}\cos\a +
q_{2}\sin \a)\sin \b.
\end{equation}
Immediate consequences:

Setting
\begin{equation}
T_{m} = p_{1} +q_{2},
\end{equation}
(called the {\it mean torsion} of $E_{3}^{2}$ at point $P\in
E_{3}^{2}$),
\begin{equation}
H = p_{2} - q_{1}
\end{equation}
(called the {\it mean curvature} of $E_{3}^{2}$ at point $P\in
E_{3}^{2}$), from (5.3.11) and (5.3.12) we have:
\begin{eqnarray}
K(\b,\a) - K(\a,\b)&=& T_{m}\cos (\a-\b),\nonumber\\T(\b,\a)-
T(\a,\b) &=& H \cos (\a-\b).
\end{eqnarray}

\begin{theorem}
The following properties hold:

1. Assuming $\b-\a\neq \pm \dfrac{\pi}{2}$, the normal curvature
$K(\b,\a)$ is symmetric with respect to the versor field
$(C,\overline{\xi}), (C,\overline{\a})$ at every point $P\in
E_{3}^{2}$, if and only if the mean torsion $T_{m}$ of $E_{3}^{2}$ vanishes
$(E_{3}^{2}$ is integrable$)$.

2. For $\b-\a \neq \pm \dfrac{\pi}{2}$, the geodesic torsion
$T(\b,\a)$ is symmetric with respect to the versor fields
$(C,\overline{\xi})$, $(C,\overline{\a})$, if and only if the
nonholonomic manifold $E_{3}^{2}$ has null mean curvature
$($$E_{3}^{2}$ is minimal$)$.\end{theorem}

\section{Parallelism, conjugation and orthonormal conjugation}
\setcounter{theorem}{0}\setcounter{equation}{0}

The notion of conjugation of versor field $(C,\overline{\xi})$
with tangent field $(C,\overline{\a})$ on the nonholonomic
manifold $E_{3}^{2}$ is that studied for the associated Myller
configuration $\mathfrak{M}_t$,

So that $(C,\overline{\xi})$ are conjugated with
$(C,\overline{\a})$ on $E_{3}^{2}$ if and only if $K(\d,d)=0$ i.e.
\begin{equation}
(p_{1}\o_{1} +p_{2}\o_{2})\widetilde{\o}_{2} - (q_{1}\o_{1}
+q_{2}\o_{2})\widetilde{\o}_{1} = 0
\end{equation}
or, by means of (5.3.11):
\begin{equation}
(p_{1}\cos\a + p_{2}\sin \a)\sin \b - (q_{1}\cos\a +q_{2}\sin
\a)\cos \b = 0.
\end{equation}

All propositions established in section 3.3, Chapter 3, can be
applied.

The notion of orthogonal conjugation of $(C, \overline{\xi})$ with
$(C,\overline{\a})$ is given by $T(\d,d) = 0$ or by
\begin{equation}
(p_{1}\o_{1} + p_{2}\o_{2})\widetilde{\o}_{1} + (q_{1}\o_{1} +
q_{2}\o_{2})\widetilde{\o}_{2} = 0
\end{equation}
or, by means of (5.3.11), it is characterized by
\begin{equation}
(p_{1}\cos \a + p_{2}\sin \a)\cos \b + (q_{1}\cos \a + q_{2}\sin
\a)\sin \b = 0.
\end{equation}

The relation of conjugation is symmetric iff $E_{3}^{2}$ is
integrable $(T_{m} = 0)$ and that of orthogonal conjugation is
symmetric iff $E_{3}^{2}$ is minimal (i.e. $H=0$).

Now we study the case of tangent versor field $(C,\overline{\xi})$
parallel along $C$, in $E_{3}^{2}$. Applying the theory of
parallelism in $\mathfrak{M}_{t}(C,\overline{\xi},\pi)$, from
Chapter 3 we obtain:

\begin{theorem}
The versors $(C,\overline{\xi})$, tangent to the manifold
$E_{3}^{2}$ along the curve $C$ is parallel in the Levi-Civita
sense if and only if the following system of equations holds
\begin{equation}
\dfrac{\widetilde{\o}_{1}}{\d
s}\left[\dfrac{d}{ds}\left(\dfrac{\widetilde{\o}_{2}}{\d s }\right)
+\dfrac{r}{ds}\dfrac{\widetilde{\o}_{1}}{\d s} \right] -
\dfrac{\widetilde{\o}_{2}}{\d
s}\left[\dfrac{d}{ds}\left(\dfrac{\widetilde{\o}_{1}}{\d s }\right)
- \dfrac{r}{ds}\dfrac{\widetilde{\o}_{2}}{\d s} \right] = 0.
\end{equation}
\end{theorem}
But, the previous system is equivalent to the system:
\begin{eqnarray}
\dfrac{d}{ds}\left(\dfrac{\widetilde{\o}_{1}}{\d s}\right) -
\dfrac{r}{ds}\dfrac{\widetilde{\o}_{2}}{\d s} &=&
h(s)\dfrac{\widetilde{\o}_{1}}{\d
s}\nonumber\\\dfrac{d}{ds}\left(\dfrac{\widetilde{\o}_{2}}{\d
s}\right) +\dfrac{r}{ds}\dfrac{\widetilde{\o}_{1}}{\d s} &=&
h(s)\dfrac{\widetilde{\o}_{2}}{\d s},\;\; \forall h(s)\neq 0.
\end{eqnarray}
Since, the tangent vector field $h(s)\overline{\xi}(s)$ has the
same direction with the versor $\overline{\xi}(s)$, from (5.2.11) we
obtain $G(\widetilde{\d},d) = h^{2}G(\d,d)$, $\widetilde{\d}$
being the direction of tangent vector $h(s)\overline{\xi}$.
Therefore, the equations $G(\d,d)=0$ is invariant with respect to
the applications $\overline{\xi}(s)\to h(s)\overline{\xi}(s)$.
Thus the equations (5.4.3) or (5.4.4) characterize the parallelism of
directions $(C, h(s) \overline{\xi}(s))$ tangent to $E_{3}^{2}$
along $C$.

All properties of the parallelism of versors $(C,\xi)$ studied for
the configuration $\mathfrak{M}_{t}$ in sections 2.8, Chapter 2 are
applied here.

Theorem of existence and uniqueness of parallel of tangent versors
$(C,\overline{\xi})$ on $E_{3}^{2}$ can be formulated exactly as
for this notion in $\mathfrak{M}_{t}$.

But a such kind of theorem can be obtained by means of the
following equation, given by (5.3.10):
\begin{equation}
G(\b,\a) \equiv \dfrac{d\b}{ds} + r_{1}(s)\cos \a +r_{2}(s)\sin \a
= 0.
\end{equation}

The parallelism  of vector field $(C,\overline{V})$ tangents to
$E_{3}^{2}$ along the curve $C$ can be studied using the form
$$
\overline{V}(s)= \|\overline{V}(s)\|\overline{\xi}(s)
$$
where $\overline{\xi}$ is the versor of $\overline{V}(s)$.

Also, we can start from the definition of Levi-Civita parallelism
of tangent vector field $(C,\overline{V})$, expressed by the
property
\begin{equation}
\dfrac{d\overline{V}}{ds} = h(s)I_{3}.
\end{equation}
Setting
\begin{equation}
\overline{V}(s) = V^{1}(s)I_{1} + V^{2}(s)I_{2}
\end{equation}
and using (5.4.8) we obtain the system of differential  equations
\begin{equation}
\dfrac{dV^{1}}{ds} - \dfrac{r}{ds}V^{2}=0,\ \
\dfrac{dV^{2}}{ds}+\dfrac{r}{ds}V^{1}=0.
\end{equation}

All properties of parallelism in Levi-Civita sense of tangent
vectors $(C,\overline{V})$ studied in section 2.8, Chapter 2 for
Myller configuration are valid. For instance

\begin{theorem}
The Levi-Civita parallelism of tangent vectors $(C,\overline{V})$
in the manifold $E_{3}^{2}$ preserves the lengths and angle of
vectors.
\end{theorem}

The concurrence in Myller sense of tangent vector fields
$(C,\overline{\xi})$ is characterized by Theorem 2.9.2 Chapter 2, by
the following equations
\begin{equation}
\dfrac{d}{ds}\left(\dfrac{c_{2}}{G}\right) = c_{1},
\end{equation}
where $c_{1} = \langle \overline{\a},\overline{\xi}\rangle$,
$c_{2} = \langle \overline{\a}, \overline{\mu}\rangle,$ $G =
G(\d,d)\neq 0.$

Consequence: the concurrence of tangent versor field
$(C,\overline{\xi})$ in $E_{3}^{2}$ for $G\neq 0$ is characterized
by
\begin{equation}
\dfrac{d}{ds}\left[\dfrac{\widetilde{\o}_{1}\o_{2} -
\widetilde{\o}_{2}\o_{1}}{\d s\cdot ds}\cdot \dfrac{1}{G}\right] =
\dfrac{\widetilde{\o}_{1}\o_{1} + \widetilde{\o}_{2}\o_{2}}{\d sds}
\end{equation}
or by:
\begin{equation}
\dfrac{d}{ds}\left[\dfrac{1}{G}\sin (\a - \b)\right] =\cos (\a-\b).
\end{equation}
The properties of concurrence in Myller sense can be taken from
section 2.8, Chapter 2.

\section{Theory of curves in $E_{3}^{2}$}
\setcounter{theorem}{0}\setcounter{equation}{0}

To a curve $C$ in the nonholonomic manifold $E_{3}^{2}$ one can
uniquely associate a Myller configuration $\mathfrak{M}_{t} =
\mathfrak{M}_{t}(C, \overline{\a}, \pi)$ where $\pi(s)$ is the
orthogonal to the normal versor $I_{3}(s)$ of $E_{3}^{2}$ at every
point $P\in C.$

Thus the geometry of curves in $E_{3}^{2}$ is the geometry of
associated Myller configurations $\mathfrak{M}_{t}$. It is
obtained as a particular case taking $\overline{\xi}(s) =
\overline{\a}(s)$ from the previous sections of this chapter.

The Darboux frame of the curve $C$, in $E_{3}^{2}$ is given by
$$
\cal{R}_{D} = \{P(s); \overline{\a}(s), \overline{\mu}^{*}(s),
I_{3}(s)\},\;\; \overline{\mu}^{*} = I_{3}\times \overline{\a}
$$
and the fundamental equations of $C$ in $E_{3}^{2}$ are as
follows:
\begin{equation}
\dfrac{d\overline{r}}{ds} = \overline{\a}(s)
\end{equation}
and
\begin{eqnarray}
\dfrac{d\overline{\a}}{ds} &=& \k_{g}(s)\overline{\mu}^{*} +
\k_{n}(s)I_{3},\\\dfrac{d\overline{\mu}^{*}}{ds} &=&
-\k_{g}(s)\overline{\a} +\tau_{g}(s)I_{3},\\\dfrac{dI_{3}}{ds}&=&
-\k_{n}(s)\overline{\a} - \tau_{g}(s)\overline{\mu}^{*}.
\end{eqnarray}
The invariant $\k_{g}(s)$ is the {\it geodesic curvature} of $C$
at point $P\in C,$ $\k_{n}(s)$ is the {\it normal curvature} of
curve $C$ at $P\in C$ and $\tau_{g}(s)$ it {\it geodesic torsion}
of $C$ at $P\in C$ in $E_{3}^{2}$.

The geometrical interpretations of these invariants are exactly
those described in the section 3.2, Chapter 3. Also, a fundamental
theorem can be enounced as in section 3.2, Theorem 3.2.2, Chapter
3.

The calculus of $\k_{g}$, $\k_{n}$ and $\tau_{g}$ is obtained by
the formulae (5.3.2), for $\overline{\xi}=\overline{\a}$.

We have:
\begin{eqnarray}
\k_{g}&=&\left\langle \overline{\a}, \dfrac{d\overline{\a}}{ds},
I_{3}\right\rangle,\\\k_{n}&=&
\left\langle\dfrac{d\overline{\a}}{ds},I_{3}\right\rangle =
-\left\langle\overline{\a},\dfrac{dI_{3}}{ds}\right\rangle,\\\tau_{g} &=&
\left\langle\overline{\a}, I_{3}, \dfrac{dI_{3}}{ds}\right\rangle.
\end{eqnarray}
By using the expressions (3.1.4), Ch. 3 we
obtain
\begin{equation}
\k_{g} = \k \sin \v^{*},\; \k_{n}=\k\cos \v^{*},\; \tau_{g} =
\tau+\dfrac{d\v^{*}}{ds}
\end{equation}
with $\v^{*} = \sphericalangle (\overline{\a}_2,I_3)$, $\overline{\a}_{2}$ being the versor of
principal normal of curve $C$ at point $P$.

A theorem of Meusnier can be formulates as in section 3.1, Chapter
3.

The line $C$ is {\it geodesic} (or {\it autoparallel}) {\it line} of the
nonholonomic manifold $E_{3}^{2}$ if $\k_{g}(s) = 0$, $\forall
s\in (s_{1},s_{2})$.

\begin{theorem}
Along a geodesic $C$ of the nonholonomic manifold $E_{3}^{2}$ the
normal curvature $\k_{n}$ is equal to $\pm \k$ and geodesic
torsion $\tau_{g}$ is equal to the torsion $\tau$ of $C$.

The differential equations of geodesics are as follows
\begin{equation}
\k_{g} \equiv
\dfrac{\o_{1}}{ds}\left[\dfrac{d}{ds}\left(\dfrac{\o_{2}}{ds}\right)
+ r_{1} \right] -
\dfrac{\o_{3}}{ds}\left[\dfrac{d}{ds}\left(\dfrac{\o_{1}}{ds}\right)
- r_{2} \right]=0.
\end{equation}
\end{theorem}

This equations are equivalent to the system of differential equations
\begin{eqnarray}
\dfrac{d}{ds}\left(\dfrac{\o_{1}}{ds}\right) -r_{2}(s)&=&
h(s)\dfrac{\o_{1}}{ds},\\\dfrac{d}{ds}\left(\dfrac{\o_{2}}{ds}\right)
+ r_{1}(s) &=& h(s)\dfrac{\o_{2}}{ds},
\end{eqnarray}
where $h(s)\neq 0$ is an arbitrary smooth function.

If we consider $\s = \sphericalangle (I_{3}, \overline{\nu}_{3})$,
where $\overline{\nu}_{3}$ the versor of binormal is of the
field $(C,I_{3})$ and $\chi_{2}$ is the torsion of $(C, I_{3})$,
then the $\k_{g}(s)$, the geodesic curvature of the field
$(C,\overline{\a})$ is obtained by the formulas of Gh. Gheorghiev:
\begin{equation}
\k_{g} = \dfrac{d\s}{ds}+\chi_{2}.
\end{equation}
It follows that: the geodesic lines of the nonholonmic manifold
$E_{3}^{2}$ are characterized by the differential equations:
\begin{equation}
\dfrac{d\s}{ds} + \chi_{2}(s) = 0. \end{equation}

By means of this equations we can prove a theorem of existence of
uniqueness of geodesic on the nonholonomic manifold $E_{3}^{2}$,
when $\s_{0} = \sphericalangle (I_{3}(s_{0}),
\overline{\nu}_{3}(s_{0}))$ is a priori given.

\section{The fundamental forms of $E_{3}^{2}$}
\setcounter{theorem}{0}\setcounter{equation}{0}

The first, second and third fundamental forms $\phi, \psi$ and
$\chi$ of the nonholonomic manifold $E_{3}^{2}$ are defined as in
the theory of surfaces in the Euclidean space $E_{3}.$

The tangent vector $d\overline{r}$ to a curve $C$ at point $P\in
C$ in $E_{3}^{2}$ is given by:
\begin{equation}
d\overline{r} = \o_{1}(s)I_{1} + \o_{2}(s)I_{2},\ \omega_3=0.
\end{equation}

The first fundamental form of $E_{3}^{2}$ at point $P\in
E_{3}^{2}$ is defined by:
\begin{equation}
\phi = \langle d\overline{r},d\overline{r}\rangle =
(\o_{1}(s))^{2} + (\o_{2}(s))^{2},\ \omega_3=0.
\end{equation}
Clearly, $\phi$ has geometric meaning and it is a quadratic positive
defined form.

The arclength of $C$ is determined by the formula (5.2.10) and the
arclength element $ds$ is given by:
\begin{equation}
ds^{2} = \phi = (\o^{1})^{2} + (\o^{2})^{2},\ \omega_3=0.
\end{equation}

The angle of two tangent vectors $d\overline{r}$ and $\d
\overline{r} =\widetilde{\o}_{1} I_{1} + \widetilde{\o}_{2}I_{2}$
is expressed by
\begin{equation}
\cos \t = \dfrac{\langle \d \overline{r}, d\overline{r}\rangle}{\d
sds} = \dfrac{\widetilde{\o}_{1}\o_{1} +
\widetilde{\o}_{2}\o_{2}}{\sqrt{(\o_{1})^{2}+(\o_{2})^{2}}\sqrt{(\widetilde{\o}_{1})^{2}
+ (\widetilde{\o}_{2})^{2} }}.
\end{equation}

The second fundamental from $\psi$ of $E_{3}^{2}$ at point $P$ is
defined by
\begin{equation}
\psi =-\langle d\overline{r}, dI_{3}\rangle = p\o_{2} - q \o_{1} =
p_{2}\o_{2}^{2} +(p_{1} - q_{2})\o_{1}\o_{2} - q_{1}\o_{1}^{2}.
\end{equation}
The form $\psi$ has a geometric meaning. It is not symmetric.

The third fundamental form $\chi$ of $E_{3}^{2}$ at point $P$ is
defined by
\begin{equation}
\chi = \langle d\overline{r},\overline{\t}\rangle,\ \omega_3=0,
\end{equation}
where $\overline{\t}$ is given by (5.1.7).

As a consequence, we have
\begin{equation}
\chi  = p\o_{1}+q\o_{2} = p_{1}\o_{1}^{2}
+(p_{2}+q_{1})\o_{1}\o_{2} + q_{2}\o_{2}^{2}
\end{equation}
$\chi$ has a geometrical meaning and it is not symmetric.

One can introduce a fourth fundamental form of $E_{3}^{2}$, by
$$
\Theta = \langle dI_{3}, dI_{3}\rangle = p^{2} + q^{2}\; (mod\;
\o_{3}).
$$
But $\Theta$ linearly depends by the forms $\phi, \psi, \chi$.
Indeed, we have
$$
\Theta = M\psi - K_{t}\phi - T_{m}\chi,
$$
where $M$ is mean curvature, $T_{m}$ is mean torsion  and $K_{t}$
is total curvature of $E_{3}^{2}$ at point $P.$ The expression of
$K_{t}$ is
\begin{equation}
K_{t} = p_{1}q_{2} - p_{2}q_{1}.
\end{equation}

The formulae (5.5.6) and (5.5.7) and the fundamental forms $\phi, \psi, \chi$
allow to express the normal curvature and geodesic torsion of a
curve $C$ at a point $P\in C$ as follows
\begin{equation}
\k_{n} = \dfrac{\psi}{\phi} = \dfrac{p\o_{2} - q\o_{1}}{\o_{1}^{2}
+\o_{2}^{2}} = \dfrac{p_{2}\o_{2}^{2} + (p_{1} - q_{2})\o_{1}\o_{2}
- q_{1}\o_{1}^{2}}{\o_{1}^{2} + \o_{2}^{2}}
\end{equation}
and
\begin{equation}
\tau_{g} = \dfrac{\chi}{\phi} = \dfrac{p\o_{1} +q\o_{2}}{\o_{1}^{2}
+ \o_{2}^{2}} = \dfrac{p_{1}\o_{1}^{2} + (p_{2} +q_{1})\o_{1}\o_{2}
+ q_{2}\o_{2}^{2}}{\o_{1}^{2} + \o_{2}^{2}}.
\end{equation}
The cases when $\k_{n} = 0$ and $\tau_{g} =0$ are important. They
will be investigated in the next section.

\section{Asymptotic lines. Curvature lines}
\setcounter{theorem}{0}\setcounter{equation}{0}

An asymptotic tangent direction $d\overline{r}$ to $E^2_3$ at a point $P\in
E_{3}^{2}$ is defined by the property $\k_{n}(s) = 0.$ A line $C$
in $E_{3}^{2}$ for which the tangent directions $d\overline{r}$
are asymptotic directions is called an {\it asymptotic line} of the
nonholonomic manifold $E_{3}^{2}$.

The asymptotic directions are characterized by the following
equation of degree 2.
\begin{equation}
p_{2}\o_{2}^{2}  + (p_{1} - q_{2})\o_{1}\o_{2} - q_{1}\o_{1}^{2}
=0.
\end{equation}
The realisant of this equation is the invariant
\begin{equation}
K_{g} = K_{t} - \dfrac{1}{4}T_{m}^{2},
\end{equation}
called the {\it gaussian curvature} of $E_{3}^{2}$ at point $P\in
E_{3}^{2}$.

We have

\begin{theorem}
At every point $P\in E_{3}^{2}$ there are two asymptotic
directions:

- real if $K_{g}<0$

- imaginary if $K_{g}>0$

- coincident if $K_{g}=0$
\end{theorem}

The point $P\in E_{3}^{2}$ is called {\it planar} if the
asymptotic directions of $E_{3}^{2}$ at $P$ are nondeterminated.

A planar point is characterized by the equations
\begin{equation}
p_{1} - q_{2} = 0,\;\; p_{2} = q_{1} = 0.
\end{equation}
The point $P\in E_{3}^{2}$ is called {\it elliptic} if the
asymptotic direction of $E_{3}^{2}$ at $P$ are imaginary and $P$
is a {\it hyperbolic} point if the asymptotic directions of
$E_{3}^{2}$ at $P$ are real.

Of course, if $P$ is a hyperbolic point of $E_{3}^{2}$ then,
exists two asymptotic line through the point $P$, tangent to the
asymptotic directions, solutions of the equations (5.7.1).

The curvature direction $d\overline{r}$ at a point $P\in E_{3}$ is
defined by the property $\tau_{g}(s) = 0.$ A line $C$ in
$E_{3}^{2}$ for which the tangent directions $d\overline{r}$ are
the curvature directions is called a {\it curvature line} of the
nonholonomic manifold $E_{3}^{2}$.

The curvature directions are characterized by the following second
order equations.
\begin{equation}
p_{1}\o_{1}^{2} + (p_{2}+q_{1})\o_{1}\o_{2} + q_{2}\o_{2}^{2} =0.
\end{equation}

The realisant of this equations is
\begin{equation}
T_{t} = K_{t} - \dfrac{1}{4}M^{2}.
\end{equation}
We have

\begin{theorem}
At every point $P\in E_{3}^{2}$ there are two curvature directions

real, if $T_{t}<0$

imaginary, if $T_{t} >0$

coincident, if $T_{t} = 0$
\end{theorem}

The curvature lines on $E_{3}^{2}$ are obtained by integrating the
equations (5.7.1) in the case $T_{t}\leq 0$.

\begin{remark}
In the case of surfaces $(T_{m}=0)$, the curvature lines coincides
with the lines of extremal normal curvature.
\end{remark}

\section{The extremal values of $\k_{n}$. Euler formula. Dupin indicatrix}
\setcounter{theorem}{0}\setcounter{equation}{0}

The extremal values at a point $P\in E_{3}^{2}$, of the normal
curvature $\k_{n} = \dfrac{\psi}{\phi}$ when $(\o_{1},\o_{2})$ are
variables are given by $\phi \dfrac{\p \psi}{\p \o_{i}} -
\psi\dfrac{\p\phi}{\p \o_{i}} = 0,$ $(i=1,2)$ which are equivalent
to the equations
\begin{equation}
\dfrac{\p\psi}{\p w_{i}} - \k_{n}\dfrac{\p\phi}{\p\o_{i}} = 0,\;\; (i=1,2).
\end{equation}

Taking into account the form (5.6.5) of $\psi$ and (5.6.3) of $\phi$,
the system (5.8.1) can be written:
\begin{eqnarray}
-(q_{1} + \k_{n})\o_{1} + \dfrac{1}{2}(p_{1}
-q_{2})\o_{2}&=&0\\\dfrac{1}{2}(p_{1} - q_{1})\o_{1}
+(p_{2}-\k_{n})\o_{2} &=& 0.
\end{eqnarray}

It is a homogeneous and linear system in $(\o_{1},\o_{2})$-which
gives the directions $(\o_{1},\o_{2})$ of extremal values of
normal curvature.

But, the previous system has solutions if and only if the following equations hold
$$
\left|\begin{array}{ccc} -(q_{1} +\k_{n})&\dfrac{1}{2}(p_{1} -
q_{2})\\\dfrac{1}{2}(p_{1}-q_{2})&p_{2} - \k_{n}
\end{array}\right|=0
$$
equivalent to the following equations of second order in $\k_{n}$:
\begin{equation}
\k_{n}^{2} - (p_{2} - q_{1})\k_{n} +p_{1}q_{2} - p_{2}q_{1} -
\dfrac{1}{4}(p_{1} + q_{2})^{2} = 0.
\end{equation}
This equation has two real solutions $\dfrac{1}{R_{1}}$ and
$\dfrac{1}{R_{2}}$ because its realisant is
\begin{equation}
\D = (p_{2} +q_{1})^{2} + (p_{1} -q_{2})^{2}.
\end{equation}
We have $\D\geq 0.$ Therefore the solutions
$\dfrac{1}{R_{1}},\dfrac{1}{R_{2}}$ are real, different or equal,
$\dfrac{1}{R_{1}}, \dfrac{1}{R_{2}}$ are the extremal values of
normal curvature $\k_{n}$.

Let us denote \begin{equation} H = \dfrac{1}{R_{1}}+\dfrac{1}{R_{2}}
=p_{2} - q_{1}
\end{equation}
called the mean curvature of the nonholonomic manifold $E_{3}^{2}$
at point $P$, and
\begin{equation}
K_{g} = \dfrac{1}{R_{1}}\dfrac{1}{R_{2}} = p_{1}q_{2} - p_{2}q_{1} -
\dfrac{1}{4}(p_{1} +q_{2}).
\end{equation}
called the Gaussian curvature of $E_{3}^{2}$ at point $P$.

Substituting $\dfrac{1}{R_{1}},$ $\dfrac{1}{R_{2}}$ in the equations
(5.8.2) we have the directions of extremal values of the normal
curvature-called the {\it principal directions} of $E_{3}^{2}$ at
point $P\in E_{3}^{2}$. These directions are obtained from (5.8.2)
for $\k_{n} = \dfrac{1}{R_{1}}$, $\k_{n}=\dfrac{1}{R_{2}}$. Thus,
one obtains the following equations which determine the principal
directions:
\begin{equation}
(p_{1} - q_{2})\o_{1}^{2} +2(p_{2}+q_{1})\o_{1}\o_{2} - (p_{1} -
q_{2})\o_{2}^{2}=0.
\end{equation}
Let $(\o_{1}, \o_{2}), (\widetilde{\o}_{1},\widetilde{\o}_{2})$
the solution of (5.8.8). Then $d\overline{r} = \o_{1}I_{1} +
\o_{2}I_{2}$ and $\d \overline{r} = \widetilde{\o}_{1}I_{1} +
\widetilde{\o}_{2}I_{2}$ are the principal directions on
$E_{3}^{2}$ at point $P.$

The principal directions $d\overline{r},\d\overline{r}$ of $E_{3}^{2}$ in
every point $P$ are real and orthogonal.

Indeed, because $\D>0$, and if $\dfrac{1}{R_{1}}\neq
\dfrac{1}{R_{2}}$ we have $\langle \d \overline{r},
d\overline{r}\rangle = 0$.

The curves on $E_{3}^{2}$ tangent to $\d \overline{r}$ and
$d\overline{r}$ are called the {\it lines of extremal normal curvature}.
So, we have

\begin{theorem}
At every point $P\in E_{3}^{2}$ there are two real and orthogonal
lines of extremal normal curvature.\end{theorem}

Assuming that the frame $\mathcal{R} = (P; I_{1}, I_{2},I_{3})$
has the vectors $I_{1}, I_{2}$ in the principal directions $\d
\overline{r}, d\overline{r}$ respectively, thus the equation (5.8.8)
implies
\begin{equation}
p_{1} - q_{2} = 0
\end{equation}
and the extremal values of normal curvature $\k_{n}$ are given by
\begin{equation}
\k_{n}^{2} - (p_{2}-q_{1})\k_{n}-p_{2}q_{1} =0.
\end{equation}
We have
\begin{equation}
\dfrac{1}{R_{1}} = p_{2},\;\;\dfrac{1}{R_{2}} = -q_{1}.
\end{equation}

Denoting by $\a$ the angle between the versor $I_{1}$ and the
versor $\overline{\xi}$ of an arbitrary direction at $P$, tangent
to $E_{3}^{2}$ we can write the normal curvature $\k_{n}$ from
(5.6.8) in the form
\begin{equation}
 \k_{n} = \dfrac{1}{R_{1}}\cos^{2}\a +
 \dfrac{1}{R_{2}}\sin^{2}\a,\;\; \a =
 \sphericalangle(\overline{\xi},I_{1}).
\end{equation}

The formula (5.8.12) is called the {\it Euler formula} for normal
curvatures on the nonholonomic manifold $E_{3}^{2}$.

Consider the tangent vector $\overrightarrow{PQ} =
|\k_{n}|^{-1/2}\overline{\xi}$, i.e.
$$
\overrightarrow{PQ} = |\k_{n}|^{-\dfrac{1}{2}}(\cos \a I_{1} +\sin
\a I_{2} ).
$$
The cartesian coordinate $(x,y)$ of the point $Q$, with respect to
the frame $(P; I_{1}, I_{2})$ are given by
\begin{equation}
x = |\k_{n}|^{-\dfrac{1}{2}}\cos \a, \;\; y =
|\k_{n}|^{-\dfrac{1}{2}}\sin \a.
\end{equation}

Thus, eliminating the parameter $\a$ from (5.8.12) and (5.8.13) are
gets the geometric locus of the point $Q$ (in the tangent plan
$(P;I_{1},I_{2})$):
\begin{equation}
\dfrac{x^{2}}{R_{1}} + \dfrac{y^{2}}{R_{2}} = \pm 1
\end{equation}
called {\it the Dupin$'$s indicatrix} of the normal curvature of
$E_{3}^{2}$ at point $P$.

It is formed by a pair of conics, having the axis in the principal
directions $\d \overline{r}, d \overline{r}$ and the invariants:
$H$-the mean curvature and $K_{g}$ - the gaussian curvature of
$E_{3}^{2}$ at $P$.

The Dupin indicatrix is formed by two ellipses, one real and
another imaginary, if $K_{g}>0.$ It is formed by a pair of
conjugate hyperbolae  if $K_{g}<0$, whose asymptotic lines are
tangent to the asymptotic lines of $E_{3}^{2}$ at point $P$. It
follows that the asymptotic directions of $E_{3}^{2}$ at $P$ are
symmetric with respect to the principal directions of $E_{3}^{2}$
at $P.$

In the case $K_{g} = 0$, the Dupin indicatrix of $E_{3}^{2}$ at
$P$ is formed by a pair of parallel straight lines - one real and
another imaginary.

\section{The extremal values of $\tau_{g}$. Bonnet formula. Bonnet indicatrix}
\setcounter{theorem}{0}\setcounter{equation}{0}

The extremal values of the geodesic torsion $\tau_{g} =
\dfrac{\chi}{\phi}$ at the point $P$, on $E_{3}^{2}$ can be
obtained following a similar way as in the previous paragraph.

Such that, the extremal values of $\tau_{g}$ are given by the
system of equations
$$
\dfrac{\p \chi}{\p \o_{i}} - \tau_{g}\dfrac{\p \phi}{\p \o_{i}} = 0,
(i=1,2).
$$
Or, taking into account the expressions of $\phi = (\o_{1})^{2}
+ (\o_{2})^{2}$ and $\chi = p\o_{1} + q\o_{2} = (p_{1}\o_{1}
+p_{2}\o_{2})\o_{1} +(q_{1}\o_{1}  + q_{2}\o_{2})\o_{2}$, we have:
\begin{eqnarray}
(p_{1}-\tau_{g})\o_{1} + \dfrac{1}{2}(p_{2}+q_{1})\o_{2} =
0\\\dfrac{1}{2}(p_{2} +q_{1})\o_{1} + (q_{2}-\tau_{g})\o_{2} =
0.\nonumber
\end{eqnarray}
But, these imply
\begin{equation}
\left|
\begin{array}{ccc}
p_{1} - \tau_{g}& \dfrac{1}{2}(p_{2}
+q_{1})\\\dfrac{1}{2}(p_{2}+q_{1})&q_{2} - \tau_{g}
\end{array}\right|=0
\end{equation}
or expanded:
\begin{equation}
\tau_{g}^{2} - (p_{1}+q_{2})\tau_{g} +(p_{1}q_{2} - p_{2}q_{1} -
\dfrac{1}{4}(p_{2} - q_{1})^{2}) = 0.
\end{equation}
The solutions $\dfrac{1}{T_{1}}, \dfrac{1}{T_{2}}$ of this equations
are the extremal values of geodesic torsion $\tau_{g}$. They are
real. The realisant of (5.9.3) is
\begin{equation}
\D_{1} = (p_{1} - q_{2})^{2} +(p_{2}+q_{1})^{2}
\end{equation}
Therefore $\dfrac{1}{T_{1}},$ $\dfrac{1}{T_{2}}$ are real, distinct
or coincident.

We denote
\begin{equation}
T_{m} = \dfrac{1}{T_{1}} + \dfrac{1}{T_{2}} = p_{1} +q_{2},
\end{equation}
{\it the mean torsion} of $E_{3}^{2}$ at point $P.$ And
\begin{equation}
T_{t} = \dfrac{1}{T_{1}}\dfrac{1}{T_{2}} = p_{1}q_{2} - p_{2}q_{1} -
\dfrac{1}{4}(p_{2} -q_{1})^{2} =K_{g} - \dfrac{1}{4}H^{2}
\end{equation}
is called {\it the total torsion} of $E_{3}^{2}$ at $P.$

But, as we know the condition of integrability of the Pfaf
equation $\o_{3}=0$ is
$$
d\o_{3} = 0,\; \mbox{modulo}\; \o_{3},
$$
which is equivalent to $T_{m} = p_{1} +q_{2} = 0.$

\begin{theorem}
The nonholonomic manifold $E_{3}^{2}$ of mean torsion $T_{m}$
equal to zero is a family of surfaces in $E_{3}$.
\end{theorem}
The directions $\d \overline{r} = \widetilde{\o}_{1}I_{1} +
\widetilde{\o}_{2}I_{2}$ and $d\overline{r} = \o_{1}I_{1} +
\o_{2}I_{2}$ of extremal geodesic torsion, corresponding to the
extremal values $\dfrac{1}{T_{1}}$ and $\dfrac{1}{T_{2}}$ of
$\tau_{g}$ are given by the equation obtained from (5.9.2) by
eliminating $\tau_{g}$ i.e.
\begin{equation}
(p_{2} + q_{1})\o_{1}^{2} - (p_{1}-q_{2})\o_{1}\o_{2} -
(p_{2}+q_{1})\o_{2}^{2} = 0.
\end{equation}

Thus $\d \overline{r}$ and $d\overline{r}$ are real and orthogonal
at every point $P\in E_{3}^{2}$ for which $\dfrac{1}{T_{1}}\neq
\dfrac{1}{T_{2}}$. $\d \overline{r}$ and $d \overline{r}$ are
called the direction of the {\it extremal geodesic torsion} at
point $P\in E_{3}^{2}.$

\begin{theorem}
At every point $P\in E_{3}^{2}$, where $\dfrac{1}{T_{1}}\neq
\dfrac{1}{T_{2}}$ there exist two real and orthogonal directions of
extremal geodesic torsion.
\end{theorem}

Now, assuming that, at point $P$, the frame $\cal{R} = (P;I_{1},
I_{2}, I_{3})$ has the vectors $I_{1}$ and $I_{2}$ in the
directions $\d \overline{r}$, $d\overline{r}$ of the extremal
geodesic torsion, respectively, from (5.9.7) we deduce
\begin{equation}
p_{2} + q_{1} =0.
\end{equation}

Thus, the equation (5.9.3) takes the form
\begin{equation}
\tau_{g}^{2} - (p_{1}+q_{2})\tau_{g} +p_{1}q_{2} = 0.
\end{equation}
Its solutions are
\begin{equation}
\dfrac{1}{T_{1}} = p_{1},\;\; \dfrac{1}{T_{2}} = q_{2}
\end{equation}
and, setting $\b = \sphericalangle (\overline{\xi}, I_{1})$, from
$\tau_{g} = \dfrac{\chi}{\phi}$ one obtains [65] the so called {\it
Bonnet formula} giving $\tau_{g}$ for the nonholonomic manifold
$E_{3}^{2}$:
\begin{equation}
\tau_{g}=\dfrac{\cos^{2}\b}{T_{1}}+\dfrac{\sin^{2}\b}{T_{2}},\; \b =
\sphericalangle (\overline{\xi},I_{1}).
\end{equation}

If $T_{m} = 0,$ (5.9.11) reduce to the very known formula from
theory of surfaces:
\begin{equation}
\tau_{g} = \dfrac{1}{T_{1}}\cos 2\b.
\end{equation}

By means of the formula (5.9.9) we can determine an indicatrix of
geodesic torsions.

In the tangent plane at point $P$ we take the cartesian frame $(P;
I_{1}, I_{2})$, $I_{1}$ having the direction of extremal geodesic
torsion $\d \overline{r}$ corresponding to $\dfrac{1}{T_{1}}$ and
$I_{2}$ having the same direction with $d\overline{r}$
corresponding to $\dfrac{1}{T_{2}}$. Consider the point $Q^{\prime}
\in \pi(s)$ given by
\begin{equation}
\overrightarrow{PQ^{\prime}} = |\tau_{g}|^{-1}\overline{\xi} =
|\tau_{g}|^{-1}(\cos \b I_{1} + \sin \b I_{2})
\end{equation}
with $\b = \sphericalangle (\overline{\xi},I_{1})$. The cartesian
coordinate $(x,y)$ at point $Q^{\prime}$:
$\overrightarrow{PQ^{\prime}} = xI_{1} +y I_{2}$, by means of
(5.9.13) give us
\begin{equation}
x = |\tau_{g}|^{-1}\cos\b,\;\; y = |\tau_{g}|^{-1}\sin \b.
\end{equation}

Eliminating the parameter $\b$ from (5.9.14) and (5.9.11) one obtains
the locus of point $Q^{\prime}$:
\begin{equation}
\dfrac{x^{2}}{T_{1}} + \dfrac{y^{2}}{T_{2}} = \pm 1,
\end{equation}
called the {\it indicatrix of Bonnet} for the geodesic torsions at
point $P\in E_{3}^{2}$.

It consists of two conics having the axis the tangents in the
directions of extremal values of geodesic torsion. The invariants
of these conics are $T_{m}$-the mean torsion and $T_{t}$-the total
torsion of $E_{3}^{2}$ at point $P.$

The Bonnet indicatrix is formed by two ellipses, one real and
another imaginary, if and only if the total torsion $T_{t}$ is positive. It
is composed by two conjugated hyperbolas if $T_{t}<0.$ In this
case the curvature lines of $E_{3}^{2}$ at $P$ are tangent to the
asymptotics of these conjugated hyperbolas. Finally, if $T_{t} =
0,$ the Bonnet indicatrix (5.9.15) is formed by two pair of parallel
straightlines, one being real and another imaginary.

\noindent {\bf Remark 5.9.1}
1. If $T_{t}<0,$ the directions of asymptotics of the hyperbolas
$(5.9.15)$ are symmetric to the directions of the extremal geodesic
torsion.

2. The direction of extremal geodesic torsion are the bisectrices
of the principal directions.

3. In the case of surfaces, $T_{m} =0,$ the Bonnet indicatrix is
formed by two equilateral conjugates hyperbolas.

Assuming that the Dupin indicatrix is formed by two conjugate
hyperbolas and the Bonnet indicatrix is formed  by two conjugate
hyperbolas, we present here only one hyperbolas from every of this
indicatrices as follows (Fig. 1):

\newpage

\includegraphics{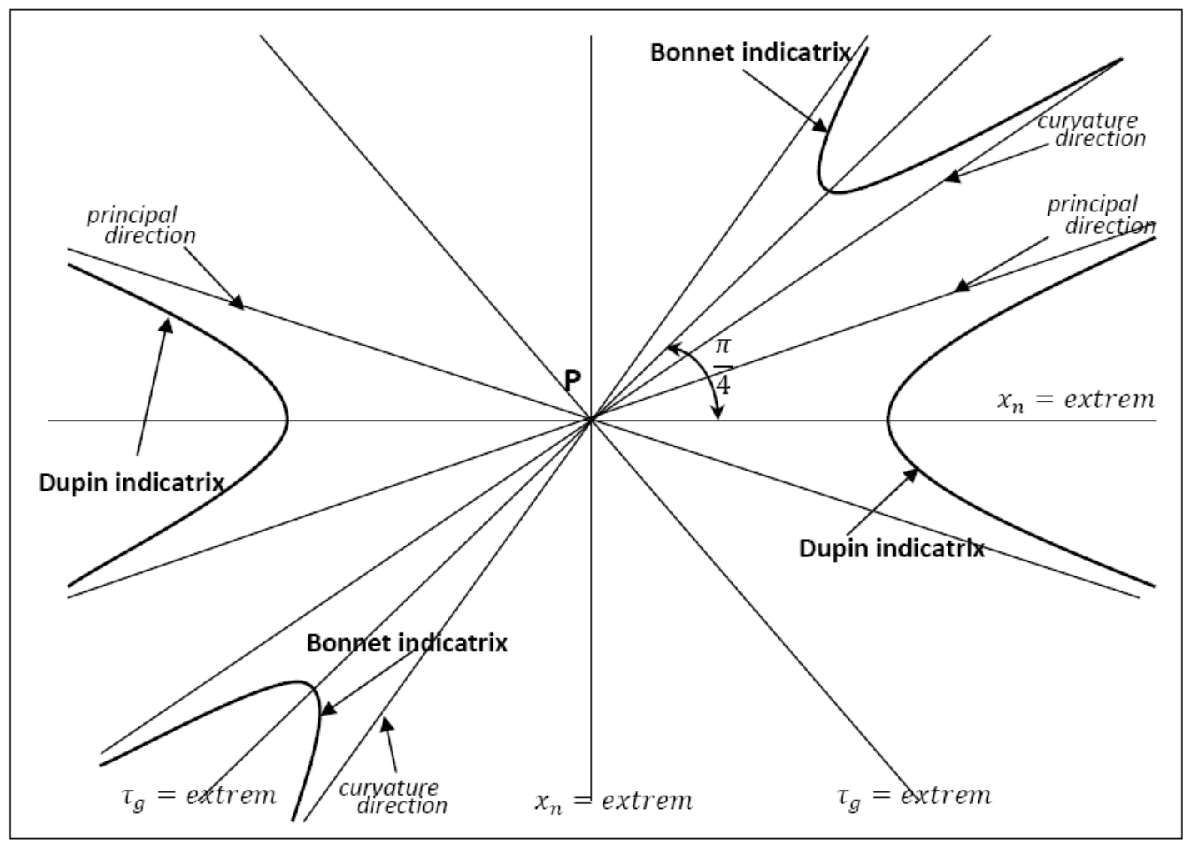}
\centerline{\footnotesize{Fig. 1}}

\section{The circle of normal curvatures and geodesic torsions}
\setcounter{theorem}{0}\setcounter{equation}{0}

Consider a curve $C\subset E_{3}^{2}$ and its tangent versor
$\overline{\a}$ at point $P\in C.$ Denoting by $\a =
\sphericalangle (\overline{\a},I_{1})$ and using the formulae
(5.3.8) we obtain \begin{eqnarray} \k_{n}&=& p_{2}\sin^{2}\a +
(p_{1} - q_{2})\sin \a \cos \a -
q_{1}c\cos^{2}\a\nonumber\\\tau_{g}&=& p_{1}\cos^{2}\a +(p_{2} +
q_{1})sin \a \cos \a +q_{2}\sin^{2}\a.
\end{eqnarray}
Eliminating the parameter $\a$ from (5.10.1) we deduce
\begin{equation}
\k_{n}^{2} +\tau_{g}^{2} -H \k_{n} - T_{m}\tau_{g} - K_{g} = 0
\end{equation}
where $H = p_{2}-q_{1}$ is the mean curvature of $E_{3}^{2}$ at
point $P$, $T_{m} = p_{1} + q_{2}$ is the mean torsion of
$E_{3}^{2}$ at $P$ and $K_{g} = p_{1}q_{2} - p_{2}q_{1}$ is the
Gaussian curvature of $E_{3}^{2}$ at point $P$. Therefore we have
\begin{theorem}
1. The normal curvature $\k_{n}$ and the geodesic torsion
$\tau_{g}$ at a point $P\in E_{3}^{2}$ satisfy the equations
$(5.10.2)$.

2. In the plan of variable $(\k_{n}, \tau_{g})$ the equation
$(5.10.2)$ represents a circle with the center $\(\dfrac{H}{2},
\dfrac{T_{m}}{2}\)$ and of the radius given by
\begin{equation}
R^{2} = \dfrac{1}{4}(H^{2}+T_{m}^{2})-K_{g}.
\end{equation}
\end{theorem}
Evidently, this circle is real if we have
$$
K_{g} <\dfrac{1}{4}(H^{2} + T_{m}^{2}).
$$
It is a complex circle if $K_{g}>\dfrac{1}{4}(H^{2} + T_m^{2})$.
This circle, for the nonholonomic manifolds $E_{3}^{2}$ was
introduced by the author [62]. It has been studied by Izu Vaisman
in [41], [42].

In the case of surfaces, $T_{m} = 0$ and (5.10.3) gives us
$$
R^{2} = \dfrac{1}{4}H^{2} - K_{g} = \dfrac{1}{4}\(\dfrac{1}{R_{1}} +
\dfrac{1}{R_{2}}\)^{2} -\dfrac{1}{R_{1}R_{2}}
=\dfrac{1}{4}\(\dfrac{1}{R_{1}} - \dfrac{1}{R_{2}}\)^{2}.
$$
So, for surfaces the circle of normal curvatures and geodesic torsions is real.

If $T_{m} = 0$, the total torsion is of the form
\begin{equation}
T_{t} = -\dfrac{1}{4}\left(\dfrac{1}{R_{1}} -
\dfrac{1}{R_{2}}\right)^{2}.
\end{equation}

This formula was established by Alexandru Myller [35]. He proved
that $T_{t}$ is just ``the curvature of Sophy Germain and Emanuel
Bacaloglu'' ([35]).

\section{The nonholonomic plane and sphere}
\setcounter{theorem}{0}\setcounter{equation}{0}

The notions of nonholonomic plane  and nonholonomic sphere have
been defined by Gr. Moisil who expressed the Pfaff equations, provided the existence of these manifolds, and R. Miron [64] studied them. Gh.
Vr\u{a}nceanu, in the book [44] has studied the notion of non
holonomic quadrics.

The Pfaff equation of the nonholonomic plane $\Pi_{3}^{2}$ given
by Gr. Moisil [30], [71] is as follows
\begin{equation}
\hspace*{14mm}\o\equiv (q_{0}z - r_{0}y +a)dx +(r_{0}x - p_{0}z +b)dy +(p_{0}y -
q_{0}x +c)dz = 0,
\end{equation}
where $p_{0}, q_{0}, r_{0}, a,b,c$ are real constants verifying
the nonholonomy condition: $ap_{0} + bq_{0} + cr_{0}\neq 0$.

While, the nonholonomic sphere $\Sigma_{3}^{2}$ has the equation
given by Gr. Moisil [30]:
\begin{eqnarray}
\hspace*{5mm}\o&\equiv & [2x (ax+by+cz) - a(x^{2} + y^{2}+z^{2}) + \mu x
+q_{0}z - r_{0}y + h] dx \nonumber\\&&+[2y(ax+by+cz)- b(x^{2} +
y^{2}+z^{2}) + \mu y +r_{0}x - p_{0}z + k]dy \nonumber\\&& +
[2z(ax+by+cz) - c(x^{2} + y^{2}+z^{2}) + \mu z +p_{0}y -
q_{0}x+l]dz.
\end{eqnarray}
$\mu, p_{0}, q_{0}, r_{0}, a,b,c, h,k,l$ being constants which
verify the condition $\o\wedge d\o \neq 0.$

In this section we investigate the geometrical properties of these
special manifolds.

1. {\it The nonholonomic plane $\Pi_{3}^{2}$}

The manifold $\Pi_{3}^{2}$ is defined by the property $\psi\equiv
0$, $\psi$ being the second fundamental form (5.6.3) of $E_{3}^{2}$.

From the formula (5.6.3) is follows
\begin{equation}
p_{2} = q_{1} = p_{1} - q_{2} = 0.
\end{equation}
The nonholonomic plane $\Pi_{3}^{2}$ has the invariants $H,
T_{m},\ldots$ given by
\begin{equation}
H=0, T_{m}\neq 0, T_{t} = \dfrac{1}{4}T_{m}^{2}, K_{t} = T_{t},
K_{g} = 0.
\end{equation}
Conversely, (5.11.4) imply (5.11.3).

Therefore, the properties (5.11.4), characterize the nonholonomic
plane $\Pi_{3}^{2}$.

Thus, the following theorems hold:

\begin{theorem}
1. The following line on $\Pi_{3}^{2}$ are nondetermined:

-The lines tangent to the principal directions.

- The lines tangent to the directions of extremal geodesic
torsion.

- The asymptotic lines.
\end{theorem}

\begin{theorem}
Also, we have

- The lines of curvature coincide with the minimal lines,
$\langle \d \overline{r}, \d \overline{r}\rangle=0$. The normal
curvature $\k_{n}\equiv 0.$

- The geodesic of $\Pi_{3}^{2}$ are straight lines.
\end{theorem}

Consequences: the nonholonomic manifold $\Pi_{3}^{2}$ is totally
geodesic.

Conversely, a total geodesic manifold $\Pi_{3}^{2}$ is a
nonholonomic plane.

If $T_{m} = 0,$ $\Pi_{3}^{2}$ is a family of ordinary planes from
$E_{3}^{2}.$

\medskip
2. {\it The nonholonomic sphere, $\Sigma_{3}^{2}$}.

A nonholonomic manifold $E_{3}^{2}$ is a nonholonomic sphere
$\Sigma_{3}^{2}$ if the second fundamental from $\psi$ is
proportional to the first fundamental form $\Phi.$

By means of (5.6.3) it follows that $\Sigma_{3}^{2}$ is
characterized by the following conditions
\begin{equation}
p_{1} - q_{2} = 0,\;\; p_{2} +q_{1} = 0.
\end{equation}

It follows that we have
\begin{equation}
\psi = \dfrac{1}{2}H\phi,\ \chi = \dfrac{1}{2}T_{m}\phi,\ \Theta^{2} =
K_{t}\phi.
\end{equation}
Thus, we can say

\begin{proposition}
1. The normal curvature $\k_{n}$, at a point $P\in \Sigma^{2}_{3}$
is the same $\(=\dfrac{1}{2}H\)$ in all tangent direction at point
$P$.

2. The geodesic torsion $\tau_{g}$ at point $P\in \Sigma_{3}^{2}$
is the same in all tangent direction at $P.$
\end{proposition}

\begin{proposition}
1. The principal directions at point $P\in E_{3}^{2}$ are
nondetermined.

2. The directions of extremal geodesic torsion at $P$ are non
determinated, too.
\end{proposition}

\begin{proposition}
At every point $P\in \Sigma_{3}^{2}$ the following relations hold:
\begin{equation}
K_{g} = \dfrac{1}{4}H^{2},\;\; T_{t} = \dfrac{1}{4}T_{m}^{2},\;\;
K_{t} = \dfrac{1}{4}(H^{2} +T_{m}^{2}),
\end{equation}
\end{proposition}
Conversely, if the first two relations (5.11.7) are verified at any
point $P\in E_{3}^{2}$, then $E_{3}^{2}$ is a nonholonomic sphere.

Indeed, (5.11.7) are the consequence of (5.11.5). Conversely, from
$K_{g} = \dfrac{1}{4}H^{2}$, $T_{t} = \dfrac{1}{4}T_{m}^{2}$ we
deduce (5.11.5).

\begin{proposition}
The nonholonomic sphere $\Sigma_{3}^{2}$ has the properties
\begin{equation}
K_{g}>0, T_{t} >0, K_{t} >0.
\end{equation}
\end{proposition}

\begin{proposition}
1. If $H\neq 0$, the asymptotic lines of $\Sigma_{3}^{2}$ are
imaginary.

2. The curvature lines of $\Sigma_{3}^{2}$, $(T_{m}\neq 0)$ are
imaginary.
\end{proposition}

The geodesic of $\Sigma_{3}^{2}$ are given by the equation (5.5.7)
or by the equation (5.5.11).

\begin{theorem}
The geodesics of nonholonomic sphere $\Sigma_{3}^{2}$ cannot be
the plane curves.
\end{theorem}

\begin{proof}
By means of theorem 5.5.1, along a geodesic $C$ we have $\k_{n} =
|\k|,$ $\tau_{g} = \tau$, $\k$ and $\tau$ being curvature and
torsion of $C$. From (11.6) we get
$$
\k_{n} = \dfrac{1}{2}H,\;\; \tau_{g} = \dfrac{1}{2}T_{m}.
$$
So, the torsion $\tau$ of $C$ is different from zero. It can not be
a plane curve.

We stop here the theory of nonholonomic manifolds $E_{3}^{2}$ in
the Euclidean space $E_{3}$, remarking that the nonholonomic quadric was investigated by G. Vr\u{a}nceanu and A. Dobrescu [8], [44], [45].
\end{proof}
\newpage
\thispagestyle{empty}

\backmatter
 
\tableofcontents
%\addcontentsline{toc}{chapter}{Index} \balance \printindex
\newpage
{---}
\vskip4cm
\begin{figure}[h]
\begin{center}
\includegraphics[width=7cm,height=7cm]{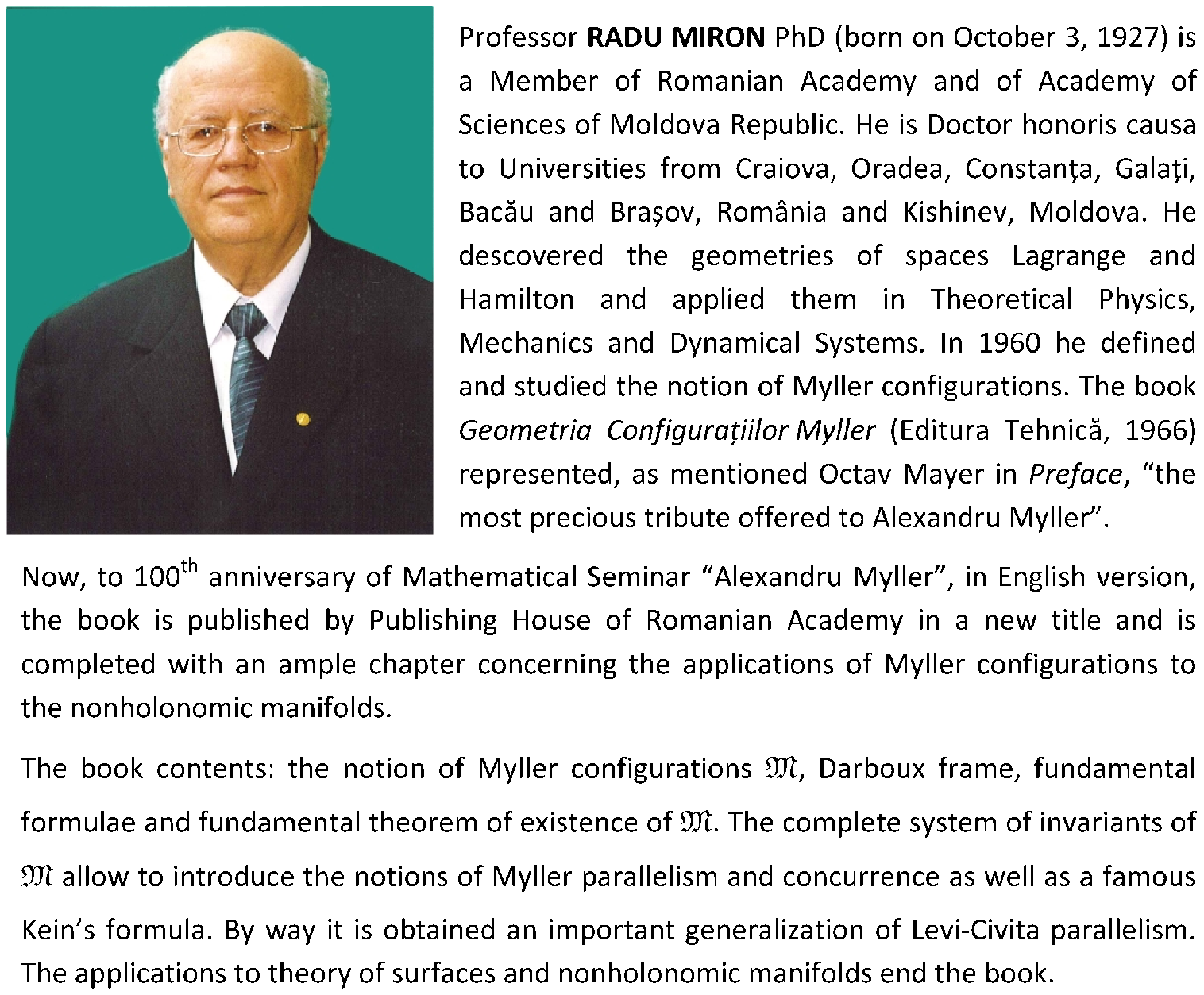}
\end{center}
\end{figure}
\end{document}